\documentclass[11pt,reqno]{amsart}
\usepackage[letterpaper,margin=1in]{geometry}
\usepackage{lmodern}
\usepackage{amsmath,amssymb,amsthm,mathtools,amscd}
\usepackage{graphicx}
\usepackage{tikz-cd}
\usepackage{mathrsfs}
\usepackage{microtype}
\usepackage{needspace}
\usepackage[colorlinks=false,pdfborder={0 0 0}]{hyperref}

\newtheorem{theorem}{Theorem}[section]
\newtheorem{proposition}[theorem]{Proposition}
\newtheorem{lemma}[theorem]{Lemma}
\newtheorem{corollary}[theorem]{Corollary}

\theoremstyle{definition}
\newtheorem{definition}[theorem]{Definition}
\newtheorem{assumption}[theorem]{Assumption}
\newtheorem{example}[theorem]{Example}

\theoremstyle{remark}
\newtheorem{remark}[theorem]{Remark}

\numberwithin{equation}{section}

\title[MIXED WEIGHT FORMAL LINEARITY AND COHOMOLOGY JUMP LOCI]{MIXED WEIGHT FORMAL LINEARITY AND \texorpdfstring{$\ell$}{l}-ADIC COHOMOLOGY JUMP LOCI}
\author{NINGYI LI}

\begin{document}

\begin{abstract}
We prove a mixed weight formal linearity theorem for diagonal power endomorphisms $\Phi$ of formal multiplicative groups over a complete discrete valuation ring of mixed characteristic $(0,\ell)$. On the coordinate ring, every prime ideal stable under the pullback $\Phi^\#$ and not containing $\ell$ defines a $\Phi$-stable formal subtorus. Consequently, a flat reduced equidimensional closed formal subscheme $Z$ becomes a finite union of torsion translates of $\Phi$-stable formal subtori over the valuation ring of a finite extension of the coefficient field if and only if there are integers $a>b\ge0$ for which the scheme theoretic images satisfy $\Phi^a(Z)\subset\Phi^b(Z)$. For an arithmetic constructible complex and an arithmetic quotient whose rational dual is mixed with strictly positive weights, the rigid form of this criterion gives the quasilinearity conjectured by Esnault and Kerz in the first paragraph of \cite[Remark~6.3]{EK21}. More precisely, after a finite extension of the coefficient field, both the jump loci for (compactly supported) cohomology and the Fourier--Mellin nonflat loci are finite unions of torsion translates of rigid subtori. All these loci carry canonical reduced rigid analytic structures.
\end{abstract}

\maketitle
\vfill
\begin{figure}[!h]
  \centering
  \includegraphics[width=.77\textwidth]{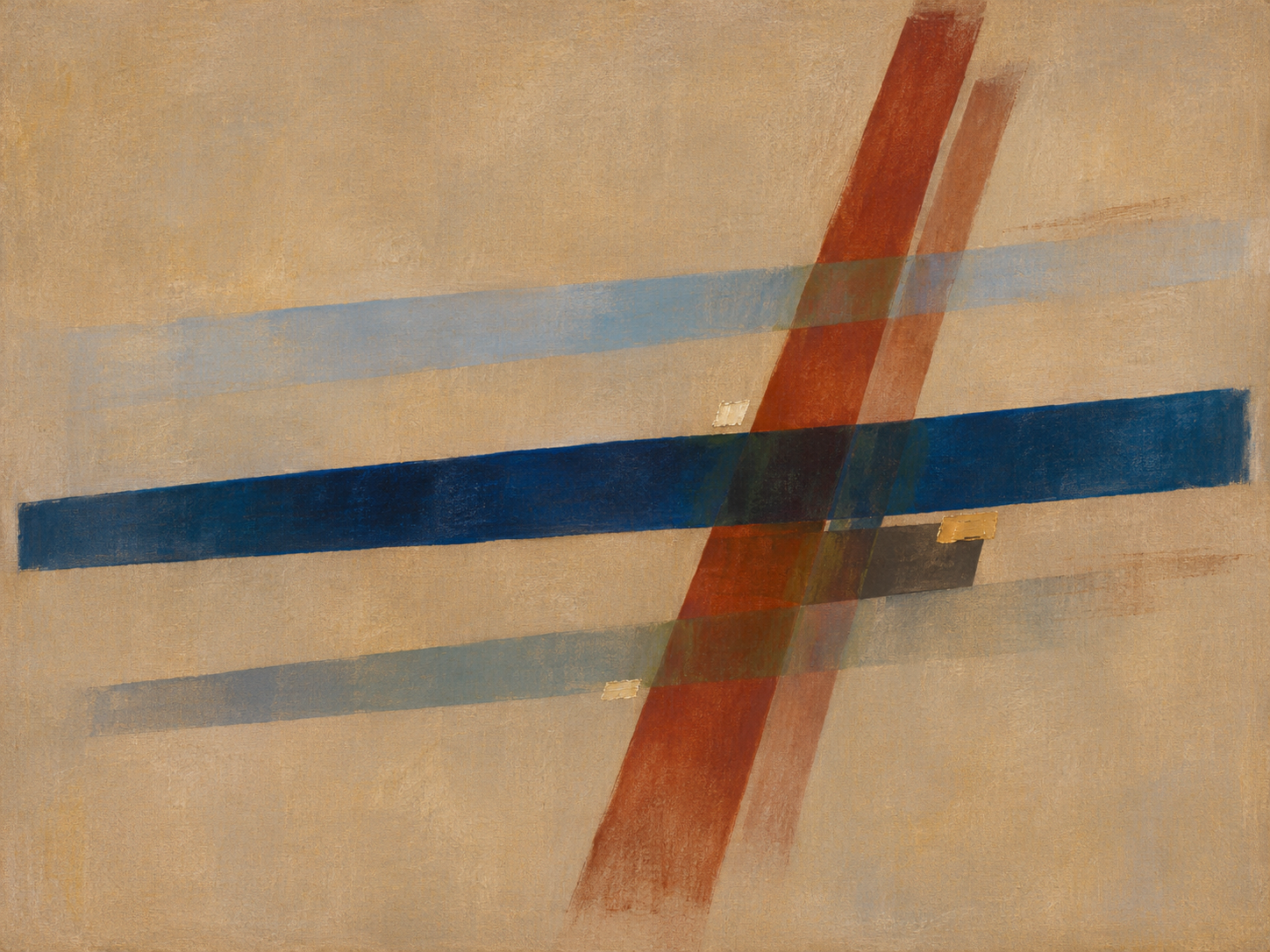}
  \caption{\textit{Geometric Rhythms in Blue and Rust}}
\end{figure}
\vfill
\clearpage
\tableofcontents

\section{Introduction}

Power endomorphisms impose strong restrictions on closed subschemes of formal multiplicative groups. De Jong established a formal linearity theorem for the scalar power map \cite{Jon00}. We prove a mixed weight extension for diagonal power endomorphisms $\Phi$ whose coordinate exponents are arbitrary positive powers of $\ell$.

For a flat reduced equidimensional closed formal subscheme $Z$, the formal criterion compares the images $\Phi^a(Z)$ and $\Phi^b(Z)$ for $a>b\ge0$. The inclusion $\Phi^a(Z)\subset\Phi^b(Z)$ holds for some pair of integers $a>b\ge0$ if and only if $Z$ becomes a finite union of torsion translates of $\Phi$-stable formal subtori over the valuation ring of a finite extension of the coefficient field. We then prove a rigid quasilinearity criterion for those closed subspaces of a rigid character torus that are defined by radical $\ell$-saturated ideals in the integral model. We apply this criterion after passing to a lattice that splits along the weights and then descend quasilinearity along the induced finite isogeny.

For a smooth complex quasiprojective variety, the cohomology jump loci of rank one local systems are finite unions of torsion translates of subtori \cite{BW15}, extending Simpson's theorem for smooth projective varieties \cite{Sim93}. Esnault and Kerz prove an $\ell$-adic analogue for an arithmetic complex and an arithmetic quotient whose rational dual is pure of nonzero weight \cite[Theorem~6.2]{EK21}. When the quotient is arithmetic and its dual is mixed with nonzero weights, our rigid criterion shows that both the jump loci for (compactly supported) cohomology and the nonflat loci of the Fourier--Mellin transform are quasilinear, proving in particular the conjecture in the first paragraph of \cite[Remark~6.3]{EK21}.

\medskip
\noindent\textbf{Mixed weight formal linearity.}
Let $\Lambda$ be a complete discrete valuation ring of mixed characteristic $(0,\ell)$ with fraction field $K$. Let $R=\Lambda[[x_1,\dots,x_n]]$ and $\mathfrak X=\operatorname{Spf}(R)$, and write $q_i=1+x_i$. For positive integers $t_1,\dots,t_n$, let $\Phi:\mathfrak X\to\mathfrak X$ be the endomorphism induced by the $\Lambda$-algebra homomorphism $\Phi^\#:R\to R$, $q_i\mapsto q_i^{\ell^{t_i}}$. For $m\in\mathbf Z_\ell^n$, the character $q^m=\prod_i q_i^{m_i}$ then satisfies $\Phi^\#(q^m)=q^{Um}$, where $U=\operatorname{diag}(\ell^{t_1},\ldots,\ell^{t_n})$ acts on $\mathbf Z_\ell^n$.

A primitive sublattice $M\subset\mathbf Z_\ell^n$ defines the formal subtorus on which every character $q^m$ with $m\in M$ is trivial. Its ideal is $I_M=(q^m-1\mid m\in M)$. Proposition~\ref{prop:stable-primes-binomial} shows that every prime ideal $\mathfrak p\subset R$ satisfying $\ell\notin\mathfrak p$ and $\Phi^\#(\mathfrak p)\subset\mathfrak p$ is $I_M$ for a unique primitive sublattice $M$ with $U(M)\subset M$.

\begin{theorem}\label{thm:intro-formal}
Let $Z\subset\mathfrak X$ be a reduced equidimensional closed formal subscheme that is flat over $\Lambda$.
\begin{enumerate}
\item[(1)] The inclusion $\Phi^a(Z)\subset\Phi^b(Z)$ holds for some integers $a>b\ge0$ if and only if there is a finite extension $L/K$, with valuation ring $\mathcal O_L$, such that $Z_{\mathcal O_L}$ is a finite scheme theoretic union of torsion translates of $\Phi_{\mathcal O_L}$-stable formal subtori of $\mathfrak X_{\mathcal O_L}$.
\item[(2)] Moreover, $Z$ is a $\Phi$-stable formal subtorus if and only if the defining ideal of $Z$ is prime and $\Phi^c(Z)\subset Z$ for some integer $c\ge1$.
\end{enumerate}
\end{theorem}

Proposition~\ref{prop:formal-dynamical-criterion} proves Theorem~\ref{thm:intro-formal}. Both the classification of stable prime ideals and the formal criterion also hold when $\Phi^\#$ acts on $\Lambda$ by an arbitrary automorphism $\tau$ (Proposition~\ref{prop:semilinear-primes} and Corollary~\ref{cor:semilinear-formal}).

\medskip
\noindent\textbf{Mixed weight rigid linearity.}
Let $\pi$ be a finite free $\mathbf Z_\ell$-module. Its completed group algebra $R_\pi=\Lambda[[\pi]]$ defines the formal character torus $\mathfrak X_\pi=\operatorname{Spf}(R_\pi)$. A basis of $\pi$ identifies $\mathfrak X_\pi$ with a split formal torus as in Theorem~\ref{thm:intro-formal}. Let $X_\pi=(\mathfrak X_\pi)_\eta$ be its rigid generic fibre \cite[Sections~7.1--7.2]{Jon95}. For every finite extension $L/K$, the $L$-points of $X_\pi$ are the continuous characters $\pi\to 1+\mathfrak m_L$, where $\mathfrak m_L$ is the maximal ideal of $\mathcal O_L$. An endomorphism $\Phi$ of $\mathfrak X_\pi$ over $\Lambda$ induces an endomorphism $\Phi:X_\pi\to X_\pi$ on the rigid generic fibre.

To prove quasilinearity of the jump loci for (compactly supported) cohomology and the Fourier--Mellin nonflat loci, we need a rigid analogue of Theorem~\ref{thm:intro-formal}. A closed rigid subspace of $X_\pi$ is \emph{standard} if it has the form $Z_J=(\operatorname{Spf}(R_\pi/J))_\eta$ for a radical $\ell$-saturated ideal $J\subset R_\pi$. Proposition~\ref{prop:standard-fm-loci} shows that the jump and nonflat loci are standard when endowed with their reduced structures.

For radical $\ell$-saturated ideals $I,J\subset R_\pi$, one has $Z_J\subseteq Z_I$ if and only if $I\subseteq J$ (Lemma~\ref{lem:faithfulness}). Moreover, taking the rigid generic fibre commutes with scheme theoretic images under a diagonal power endomorphism $\Phi$ (Lemma~\ref{lem:rigid-images}). Thus, for $\mathfrak Z_J=\operatorname{Spf}(R_\pi/J)$, the inclusion $\Phi^a(\mathfrak Z_J)\subseteq\Phi^b(\mathfrak Z_J)$ in Theorem~\ref{thm:intro-formal}(1) is equivalent to $\Phi^a(Z_J)\subseteq\Phi^b(Z_J)$ on $X_\pi$.

A \emph{rigid subtorus} is the generic fibre of a formal subtorus. A reduced closed rigid subspace is \emph{rigid quasilinear} if, after a finite extension of $K$, it is a finite closed analytic union of torsion translates of rigid subtori. It is \emph{$\Phi$-stable rigid quasilinear} if the subtori in such a decomposition can be chosen stable under $\Phi$.

\begin{theorem}\label{thm:intro-rigid}
Fix a basis $e_1,\dots,e_n$ of $\pi$ and positive integers $t_1,\dots,t_n$. The $\mathbf Z_\ell$-linear map $U:\pi\to\pi$, $e_i\mapsto\ell^{t_i}e_i$, induces the endomorphism $\Phi:X_\pi\to X_\pi$, $\chi\mapsto\chi\circ U$. Let $Z_J\subset X_\pi$ be a standard closed rigid subspace with $R_\pi/J$ equidimensional.
\begin{enumerate}
\item[(1)] $Z_J$ is $\Phi$-stable rigid quasilinear if and only if $\Phi^a(Z_J)\subset\Phi^b(Z_J)$ as closed rigid subspaces for some integers $a>b\ge0$.
\item[(2)] Moreover, $Z_J$ is a $\Phi$-stable rigid subtorus if and only if it is irreducible and $\Phi^c(Z_J)\subset Z_J$ for some integer $c\ge1$.
\end{enumerate}
\end{theorem}

We prove Theorem~\ref{thm:intro-rigid} in Theorem~\ref{thm:rigid-dynamical}. Let $Z\subset X_\pi$ be a standard closed rigid subspace. When its integral model has components of different dimensions, Corollary~\ref{cor:fixed-dim-components} applies the criterion to the rigid generic fibre of the union of the components of each fixed dimension.

\needspace{4\baselineskip}
For a finite index inclusion $\pi\subset\pi'$, restriction of characters gives a finite isogeny $f:X_{\pi'}\to X_\pi$. Choose a basis of $\pi'$ and write $q'_1,\ldots,q'_n$ for the corresponding multiplicative coordinates. For positive integers $t_1,\ldots,t_n$, let $\Phi':X_{\pi'}\to X_{\pi'}$ be the endomorphism whose pullback satisfies $(\Phi')^\#(q'_i)=(q'_i)^{\ell^{t_i}}$ for each $i$. If $f^{-1}(Z)$ satisfies the criterion of Corollary~\ref{cor:fixed-dim-components} for $\Phi'$, then $f^{-1}(Z)$ is $\Phi'$-stable rigid quasilinear and $Z$ is rigid quasilinear (Corollary~\ref{cor:isogenous-rigid-criterion}).

\medskip
\noindent\textbf{Arithmetic jump loci.}
For the arithmetic applications, let $E/\mathbf Q_\ell$ be a finite extension and take $\Lambda=\mathcal O_E$. Let $F$ be an algebraically closed field of characteristic different from $\ell$, and let $X/F$ be a separated, connected scheme of finite type. Denote by $\Gamma_\ell(X)$ the maximal torsion free quotient of the maximal abelian pro-$\ell$ quotient of $\pi_1(X)$. Fix a surjection $u:\Gamma_\ell(X)\twoheadrightarrow\pi$ onto a finite free $\mathbf Z_\ell$-module, and let $V_\pi=\pi\otimes_{\mathbf Z_\ell}\mathbf Q_\ell$.

For a model $X_{F_0}$ of $X$ over a finitely generated subfield $F_0\subset F$, let $\overline F_0$ be the algebraic closure of $F_0$ in $F$ and let $p:X\to X_{\overline F_0}=X_{F_0}\times_{F_0}\overline F_0$ be the projection. The group $G_{F_0}=\mathrm{Gal}(\overline F_0/F_0)$ acts on $X_{\overline F_0}$ through the second factor and hence on $\Gamma_\ell(X)$ via $p_*:\Gamma_\ell(X)\xrightarrow{\sim}\Gamma_\ell(X_{\overline F_0})$.

The quotient $u$ is \emph{arithmetic} if, for some such model, there is an open subgroup $H\subset G_{F_0}$ such that $h(\ker(u))=\ker(u)$ for all $h\in H$. Equivalently, $H$ preserves the image of $V_\pi^\vee\hookrightarrow H^1(X_{\overline F_0},\mathbf Q_\ell)$, $\lambda\mapsto\lambda\circ u$.

A bounded constructible complex $\mathcal F\in D_c^b(X,\overline{\mathbf Q}_\ell)$ is \emph{arithmetic} if, for some such model, there exists $\mathcal F_{\overline F_0}\in D_c^b(X_{\overline F_0},\overline{\mathbf Q}_\ell)$ such that $p^*\mathcal F_{\overline F_0}\simeq\mathcal F$ and $h^*\mathcal F_{\overline F_0}\simeq\mathcal F_{\overline F_0}$ for every $h\in G_{F_0}$. If $u$ and $\mathcal F$ are arithmetic, they admit a common model with $H=G_{F_0}$ (Remark~\ref{rem:common-descent}).

After a finite extension of $E$, choose a model of the arithmetic complex $\mathcal F$ in $D_c^b(X,E)$ and denote it again by $\mathcal F$ \cite[Sections~1.1.2--1.1.3]{Del80}. For a rigid point $\chi\in X_\pi$ with residue field $\kappa(\chi)$, let $\mathcal L_\chi$ be the rank one local system on $X$ defined by the character $\pi_1(X)\to\Gamma_\ell(X)\xrightarrow{u}\pi\xrightarrow{\chi}\kappa(\chi)^\times$. For $i\in\mathbf Z$ and $m\ge1$, the jump loci for (compactly supported) cohomology are
\[
\begin{aligned}
\Sigma_m^i(X,\mathcal F)
 &=\bigl\{\chi\in X_\pi\mid \dim_{\kappa(\chi)}H^i(X,\mathcal F\otimes_E^{\mathbf L}\mathcal L_\chi)\ge m\bigr\},\\
\Sigma_m^{i,c}(X,\mathcal F)
 &=\bigl\{\chi\in X_\pi\mid \dim_{\kappa(\chi)}H_c^i(X,\mathcal F\otimes_E^{\mathbf L}\mathcal L_\chi)\ge m\bigr\}.
\end{aligned}
\]
Let $\mathrm{FM}_\pi(X,\mathcal F)$ denote the rational Fourier--Mellin transform of $\mathcal F$ over $R_\pi[1/\ell]$ (Definition~\ref{def:fm-transform}(2)). The nonflat locus $\mathrm{NF}^i(X,\mathcal F)$ consists of those rigid points $\chi\in X_\pi$ for which $H^i(\mathrm{FM}_\pi(X,\mathcal F))$ is not flat at the maximal ideal defined by $\chi$. Proposition~\ref{prop:standard-fm-loci} equips all three loci with canonical reduced closed rigid structures defined by radical $\ell$-saturated ideals in $\mathcal O_E[[\pi]]$.

\begin{theorem}\label{thm:intro-arithmetic}
Suppose that $u:\Gamma_\ell(X)\twoheadrightarrow\pi$ and $\mathcal F$ are arithmetic relative to a common descent and that $V_\pi^\vee$ is mixed with all weights strictly positive. For $i\in\mathbf Z$ and $m\ge1$, let $Z$ be any of the three reduced closed rigid subspaces $\Sigma_m^i(X,\mathcal F)_{\mathrm{red}}$, $\Sigma_m^{i,c}(X,\mathcal F)_{\mathrm{red}}$, and $\mathrm{NF}^i(X,\mathcal F)_{\mathrm{red}}$.
\begin{enumerate}
\item[(1)] $Z$ is rigid quasilinear on $X_\pi$.
\item[(2)] Let $f:X_{\pi'}\to X_\pi$ be the finite étale isogeny of Proposition~\ref{prop:weight-adapted-torus}(2), and let $\Phi_w:X_{\pi'}\to X_{\pi'}$ be the rigid generic fibre of the mixed power endomorphism of Proposition~\ref{prop:weight-adapted-torus}(3). Then $f^{-1}(Z)$ is $\Phi_w$-stable rigid quasilinear on $X_{\pi'}$.
\end{enumerate}
\end{theorem}

To construct $f$ and $\Phi_w$, choose a lift $\sigma$ of a power of geometric Frobenius as in Lemma~\ref{lem:frob-lift}. The action of $\sigma$ on $V_\pi$ is semisimple, with weight decomposition $V_\pi=\bigoplus_r V_r$, where $V_r$ has weight $-r$. Since this decomposition need not split the lattice $\pi$, we use Proposition~\ref{prop:weight-adapted-torus}(1) to choose a $\sigma$-stable lattice $\pi'$ containing $\pi$ with finite index and satisfying $\pi'=\bigoplus_r(\pi'\cap V_r)$. Restriction of characters along $\pi\hookrightarrow\pi'$ gives $f$, and multiplication by $\ell^r$ on each summand $\pi'\cap V_r$ induces $\Phi_w$. If $V_\pi^\vee$ is pure of weight $w>0$, there is only one weight summand, so we can take $\pi'=\pi$ and $f=\mathrm{id}_{X_\pi}$ (Corollary~\ref{cor:pure-weight}).

\medskip
\noindent\textbf{Relation to earlier work.}
When $t_1=\cdots=t_n$, de Jong uses logarithmic differentials and Cartier duality to prove that every prime ideal $\mathfrak p\subset R$ satisfying $\ell\notin\mathfrak p$ and $\Phi^\#(\mathfrak p)\subseteq\mathfrak p$ is $I_M$ for a primitive sublattice $M\subset\mathbf Z_\ell^n$ \cite[Proposition~2.2]{Jon00}. In the proof of \cite[Lemma~2.1(2)]{Jon00}, he drops the stability hypothesis and evaluates a fraction $h$ at the identity to deduce $h\in K$ from $dh=0$. Remark~\ref{rem:const-field-step} shows that this implication fails without stability. For a stable prime $\mathfrak p$, our classification identifies $R/\mathfrak p$ with a formal power series ring over $\Lambda$. Every fixed fraction then lies in $K$ because its differential vanishes, so the fixed field is $K^\tau$ when the power map acts on $\Lambda$ by $\tau$ (Corollaries~\ref{cor:fixed-field} and~\ref{cor:semilinear-fixed-field}).

As written, the proof of \cite[Theorem~6.2, p.~20]{EK21} contains an error in the reduction to the maximal quotient. Replacing the prescribed arithmetic quotient $\pi$ by $\Gamma_\ell(X)$ while keeping $X$ fixed need not preserve purity. Indeed, in Example~\ref{ex:zero-weight}, $V_\pi^\vee$ is pure of weight two, whereas $H^1(X,\mathbf Q_\ell)$ has weights zero and two. To prove quasilinearity for the given arithmetic quotient $\pi$, we use Lemma~\ref{lem:frob-lift} to choose a lift of a power of geometric Frobenius acting semisimply on $V_\pi$ itself. When $V_\pi^\vee$ is pure of positive weight, Frobenius invariance allows us to apply the rigid criterion directly on $X_\pi$ and conclude that, after a finite extension of $E$, both the cohomology jump loci and the Fourier--Mellin nonflat loci are finite unions of torsion translates of rigid subtori (Corollary~\ref{cor:pure-weight}). Taking $\overline{\mathbf Q}_\ell$-points recovers both conclusions of \cite[Theorem~6.2]{EK21}.

\medskip
\noindent\textbf{Organization.}
Section~\ref{sec:formal-tori} develops the character lattice description of formal tori and their torsion translates. Sections~\ref{sec:formal-rigidity} and~\ref{sec:rigid-spaces} prove the formal and rigid criteria, respectively. Section~\ref{sec:frobenius-loci} derives the inclusions in Proposition~\ref{prop:frob-gives-phi-criterion} from Frobenius invariance and constructs the finite isogeny $f:X_{\pi'}\to X_\pi$ of Theorem~\ref{thm:intro-arithmetic}(2). Section~\ref{sec:jump-loci} applies the rigid criterion and Frobenius stability to cohomology jump loci and Fourier--Mellin nonflat loci.

\section{Preliminaries on formal tori}\label{sec:formal-tori}

Let $\Lambda$ be a complete discrete valuation ring of mixed characteristic $(0,\ell)$ with fraction field $K$ and uniformizer $\varpi$. The inclusion $\mathbf Z\hookrightarrow\Lambda$ extends to $\mathbf Z_\ell\hookrightarrow\Lambda$ by completeness. We give $\Lambda[[x_1,\dots,x_n]]$ the $(\varpi,x_1,\dots,x_n)$-adic topology. An ideal $J$ in a $\Lambda$-algebra is called \emph{$\ell$-saturated} if $\ell f\in J$ implies $f\in J$. Since $\ell$ is a unit times a positive power of $\varpi$, $\ell$-saturation of $J$ is equivalent to flatness of the quotient over $\Lambda$.

\subsection{The character lattice dictionary}\label{subsec:lattice-dictionary}
Let $\pi$ be a finite free $\mathbf Z_\ell$-module. Following \cite[Section~3.1]{EK21}, define its completed group algebra and the associated formal group by
\[
R_\pi=\Lambda[[\pi]]=\varprojlim_{N,c}(\Lambda/\varpi^c)[\pi/\ell^N\pi],
\qquad \mathfrak X_\pi=\mathrm{Spf}(R_\pi).
\]
Write $[\gamma]$ for the image of $\gamma$ under the natural inclusion $\pi\hookrightarrow R_\pi^\times$. At level $(N,c)$, $[\gamma]$ is the basis element indexed by $\gamma\bmod\ell^N\pi$. The comultiplication and counit are
\[
\Delta:R_\pi\to R_\pi\mathbin{\widehat\otimes}_\Lambda R_\pi,\quad [\gamma]\mapsto[\gamma]\mathbin{\widehat\otimes}[\gamma],
\qquad
\varepsilon:R_\pi\to\Lambda,\quad [\gamma]\mapsto1.
\]

A basis $e_1,\ldots,e_b$ of $\pi$ gives an isomorphism $R_\pi\xrightarrow{\sim}\Lambda[[x_1,\ldots,x_b]]$, $[e_i]\mapsto q_i:=1+x_i$. Indeed, at level $(N,c)$ this assignment identifies $(\Lambda/\varpi^c)[\pi/\ell^N\pi]$ with the quotient of $\Lambda[[x_1,\ldots,x_b]]$ by $\varpi^c$ and all $(1+x_i)^{\ell^N}-1$. This quotient is Artinian local since the relations reduce to $x_i^{\ell^N}=0$ modulo $\varpi$. Conversely, the binomial expansion gives $(1+x_i)^{\ell^N}-1\to0$ in the $(\varpi,x_1,\ldots,x_b)$-adic topology as $N\to\infty$. Thus the ideals defining these quotients are cofinal with the powers of $(\varpi,x_1,\ldots,x_b)$, and taking inverse limits gives the asserted isomorphism. For $a\in\mathbf Z_\ell$, define
\[
(1+x)^a=\sum_{k\ge0}\binom{a}{k}x^k\in\mathbf Z_\ell[[x]]\subset\Lambda[[x]].
\]
In these coordinates $[\gamma]=q^m=\prod_i(1+x_i)^{m_i}$ if $\gamma=\sum_i m_i e_i$.

\begin{definition}\label{def:split-formal-torus}
A \emph{split formal torus} over $\Lambda$ is a formal group scheme isomorphic to $\widehat{\mathbf{G}}_m^b = \mathfrak{X}_{\mathbf{Z}_\ell^b}$ for some integer $b \ge 0$. A \emph{formal subtorus} of a split formal torus is a closed formal subgroup that is itself isomorphic to a split formal torus.
\end{definition}

\begin{definition}\label{def:formal-torus-hom}
A \emph{homomorphism of split formal tori} over $\Lambda$ is a morphism $f:\mathfrak X\to\mathfrak Y$ over $\mathrm{Spf}(\Lambda)$ satisfying
\[
f\circ m_{\mathfrak X}=m_{\mathfrak Y}\circ(f\times_\Lambda f),
\qquad
f\circ e_{\mathfrak X}=e_{\mathfrak Y},
\]
where $m$ and $e$ denote multiplication and the identity section.
\end{definition}

Equivalently, the pullback $f^\#:\mathcal O(\mathfrak Y)\to\mathcal O(\mathfrak X)$ is a continuous $\Lambda$-algebra homomorphism satisfying
\[
\Delta\circ f^\#=(f^\#\mathbin{\widehat\otimes}_\Lambda f^\#)\circ\Delta,
\qquad
\varepsilon\circ f^\#=\varepsilon.
\]

\begin{lemma}\label{lem:grouplike}
Every unit $u\in R_\pi$ satisfying $\Delta(u)=u\mathbin{\widehat{\otimes}}u$ and $\varepsilon(u)=1$ is of the form $u=[\gamma]$ for a unique $\gamma\in\pi$.
\end{lemma}

\begin{proof}
For each $N\ge1$, write the image of $u$ in $\Lambda[\pi/\ell^N\pi]$ as $\sum_{g\in\pi/\ell^N\pi}a_g[g]$. Comparing the coefficients of $[g]\otimes[h]$ in $\Delta(u)=u\otimes u$ gives $a_g^2=a_g$ and $a_g a_h=0$ for $g\ne h$. The ring $\Lambda$ is local and $\sum_g a_g=1$, so exactly one coefficient is $1$ and the others are zero. Let $g_N\in\pi/\ell^N\pi$ be the unique element with $a_{g_N}=1$. By compatibility of the images of $u$, the reduction of $g_{N+1}$ modulo $\ell^N\pi$ is $g_N$. Thus $(g_N)_N$ defines a unique $\gamma\in\pi$, and $u=[\gamma]$ since both have image $[g_N]$ for every $N$.
\end{proof}

Use multiplicative coordinates $q_i=1+x_i$ on $\widehat{\mathbf G}_m^n$ and $q'_j=1+y_j$ on $\widehat{\mathbf G}_m^r$.

\begin{proposition}\label{prop:torus-morphisms}
There is a bijection between matrices $U\in\mathrm M_{n\times r}(\mathbf Z_\ell)$ and homomorphisms $f:\widehat{\mathbf G}_m^n\to\widehat{\mathbf G}_m^r$ over $\Lambda$, via the formula
\begin{equation}\label{eq:exponent-matrix}
f^\#((q')^m)=q^{Um}\qquad(m\in\mathbf Z_\ell^r).
\end{equation}
\end{proposition}

\begin{proof}
For $1\le j\le r$, the element $f^\#(q'_j)$ satisfies the hypotheses of Lemma~\ref{lem:grouplike} since $f^\#$ preserves comultiplication and the counit. Hence there are unique $u_{1j},\ldots,u_{nj}\in\mathbf Z_\ell$ such that $f^\#(q'_j)=\prod_{i=1}^n q_i^{u_{ij}}$. These exponents form the $j$-th column of $U=(u_{ij})$. The identity~\eqref{eq:exponent-matrix} follows first for integral exponents by multiplication and then for $m\in\mathbf Z_\ell^r$ by continuity.

Conversely, let $U=(u_{ij})$. The continuous homomorphism $\Lambda[[y_1,\ldots,y_r]]\to\Lambda[[x_1,\ldots,x_n]]$, $y_j\mapsto\prod_i q_i^{u_{ij}}-1$, exists because these images lie in $(x_1,\ldots,x_n)$. The identities $\Delta(q_i^a)=q_i^a\mathbin{\widehat\otimes}q_i^a$ and $\varepsilon(q_i^a)=1$ hold for $a\in\mathbf Z_\ell$ by continuity, so this homomorphism preserves comultiplication and the counit. It therefore defines $f$, uniquely determined by the images of the $y_j$.
\end{proof}

If $f:\widehat{\mathbf G}_m^n\to\widehat{\mathbf G}_m^r$ and $g:\widehat{\mathbf G}_m^r\to\widehat{\mathbf G}_m^s$ have matrices $U$ and $V$, respectively, then $g\circ f$ has matrix $UV$. Indeed, $(g\circ f)^\#=f^\#\circ g^\#$ acts on exponents first by $V$ and then by $U$.

\begin{definition}\label{def:character-lattice}
A \emph{character} of a split formal torus $\mathfrak X$ over $\Lambda$ is a homomorphism $\mathfrak X\to\widehat{\mathbf G}_m$ over $\Lambda$. The \emph{character lattice} $X^*(\mathfrak X)$ is the $\mathbf Z_\ell$-module of characters of $\mathfrak X$. Addition is pointwise multiplication. For $a\in\mathbf Z_\ell$, scalar multiplication satisfies $(a\chi)^\#(z)=(\chi^\#(z))^a$, where $z$ is the multiplicative coordinate on $\widehat{\mathbf G}_m$.
\end{definition}

A homomorphism $f:\mathfrak X\to\mathfrak Y$ induces the pullback $f^*:X^*(\mathfrak Y)\to X^*(\mathfrak X)$. In the coordinates of Proposition~\ref{prop:torus-morphisms}, let $U$ be the matrix of $f$. For each character $\chi:\widehat{\mathbf G}_m^r\to\widehat{\mathbf G}_m$, there is a unique $m\in\mathbf Z_\ell^r$ such that $\chi^\#(z)=(q')^m$. By~\eqref{eq:exponent-matrix},
\[
(\chi\circ f)^\#(z)=f^\#\bigl((q')^m\bigr)=q^{Um}.
\]

\begin{proposition}\label{prop:lattice-equivalence}
The character functor $X^*(-)$ defines an anti-equivalence between the category of split formal tori over $\Lambda$ and the category of finite free $\mathbf Z_\ell$-modules.
\end{proposition}

\begin{proof}
For every finite free $\mathbf Z_\ell$-module $\pi$, there is a natural isomorphism $\pi\xrightarrow{\sim}X^*(\mathfrak X_\pi)$, $\gamma\mapsto\chi_\gamma$, where $\chi_\gamma^\#(z)=[\gamma]$. Indeed, in a basis of $\pi$, $\chi_\gamma$ is the monomial character with exponent vector $\gamma$. Proposition~\ref{prop:torus-morphisms} with $r=1$ shows that every character occurs exactly once. Multiplication adds exponents, and raising a monomial to the power $a\in\mathbf Z_\ell$ multiplies its exponent by $a$. This proves essential surjectivity.

For full faithfulness, let $\pi_1$ and $\pi_2$ be finite free $\mathbf Z_\ell$-modules of ranks $n$ and $r$, and choose bases. Proposition~\ref{prop:torus-morphisms} identifies homomorphisms $\mathfrak X_{\pi_1}\to\mathfrak X_{\pi_2}$ with matrices in $\mathrm M_{n\times r}(\mathbf Z_\ell)$, which also represent linear maps $\pi_2\to\pi_1$. The pullback of characters by the homomorphism with matrix $U$ acts on exponent vectors by $m\mapsto Um$. Thus pullback gives a bijection between the two sets of morphisms, proving full faithfulness.
\end{proof}

Recall that a submodule $M\subset\pi$ is \emph{primitive} if $\pi/M$ is torsion free.

\begin{lemma}[{\cite[Proposition~1.4.6 and Lemma~1.4.11]{Sch23}}]\label{lem:primitivity}
The following are equivalent:
\begin{enumerate}
\item[(1)] $M$ is primitive.
\item[(2)] $M=(M\otimes_{\mathbf Z_\ell}\mathbf Q_\ell)\cap\pi$ inside $\pi\otimes_{\mathbf Z_\ell}\mathbf Q_\ell$.
\item[(3)] $M$ is a direct summand of $\pi$.
\end{enumerate}
\end{lemma}

\begin{samepage}
\begin{lemma}\label{lem:primitive-subtori}
\leavevmode
\begin{enumerate}
\item[(1)] The assignments
\[
M\longmapsto T_M\coloneqq\bigcap_{\gamma\in M}\ker(\chi_\gamma),
\qquad
T\longmapsto\{\gamma\in\pi\mid\chi_\gamma|_T=1\}
\]
give mutually inverse order reversing bijections between the primitive submodules of $\pi$ and the formal subtori of $\mathfrak X_\pi$.
\item[(2)] The defining ideal of $T_M$ is $I_M=([\gamma]-1\mid\gamma\in M)$. If $m_1,\ldots,m_r$ is a basis of $M$, then $I_M=([m_1]-1,\ldots,[m_r]-1)$.
\end{enumerate}
\end{lemma}
\end{samepage}

\begin{proof}
Extend a basis $m_1,\dots,m_r$ of $M$ to a basis of $\pi$. In the corresponding coordinates, the equations $[m_i]=1$ for $1\le i\le r$ cut out a coordinate subtorus. Every character in $M$ is trivial on this subtorus, so it is $T_M$. This gives (2). The quotient map $\pi\to\pi/M$ identifies $\mathfrak X_{\pi/M}$ with this coordinate subtorus. Thus for $\gamma\in\pi$ one has $\chi_\gamma|_{T_M}=1$ if and only if $\gamma\in M$.

Conversely, let $\mathfrak X_{\pi'}\hookrightarrow\mathfrak X_\pi$ be a formal subtorus. By Proposition~\ref{prop:lattice-equivalence}, restriction of characters induces a linear map $v:\pi\to\pi'$. Let $k=\Lambda/(\varpi)$ and let $\bar v:\pi/\ell\pi\to\pi'/\ell\pi'$ be the reduction of $v$ modulo $\ell$. The closed immersion induces a surjection on cotangent spaces at the identity in the special fibres. Since the linear term of $(1+x)^a-1$ is $ax$, this surjection identifies with
\[
(\pi/\ell\pi)\otimes_{\mathbf F_\ell}k
 \xrightarrow{\bar v\otimes 1}
(\pi'/\ell\pi')\otimes_{\mathbf F_\ell}k.
\]
Faithful flatness of $k$ over $\mathbf F_\ell$ implies that $\bar v$ is surjective. Hence $v$ is surjective by Nakayama's lemma over $\mathbf Z_\ell$. Its kernel $M$ is primitive and $\pi/M\cong\pi'$, which identifies the given subtorus with $T_M$. Taking common kernels reverses inclusions.
\end{proof}

\begin{corollary}\label{cor:restriction-subtorus}
Let $M\subset\pi$ be primitive. Under the anti-equivalence of Proposition~\ref{prop:lattice-equivalence}, the quotient $\pi\twoheadrightarrow\pi/M$ and the inclusion $M\hookrightarrow\pi$ induce the exact sequence of formal groups
\[
1\longrightarrow\mathfrak X_{\pi/M}\longrightarrow\mathfrak X_\pi
\xrightarrow{\mathrm{res}_M}\mathfrak X_M\longrightarrow1.
\]
\end{corollary}

\begin{proof}
Choose a direct complement of $M$ in $\pi$. Then $\mathfrak X_\pi\cong\mathfrak X_{\pi/M}\times\mathfrak X_M$, and $\mathrm{res}_M$ is projection onto the second factor.
\end{proof}

\subsection{Torsion translates of formal subtori}\label{subsec:images-translations}

\begin{definition}\label{def:formal-image}
Let $\Phi$ be an endomorphism of $\mathfrak X=\mathrm{Spf}(R)$, where $R$ is a complete Noetherian adic ring, and let $Z\subset\mathfrak X$ be a closed formal subscheme.
\begin{enumerate}
\item[(1)] The \emph{scheme theoretic image} $\Phi(Z)$ is the smallest closed formal subscheme of $\mathfrak X$ through which $\Phi|_Z:Z\to\mathfrak X$ factors as
\[
Z\longrightarrow\Phi(Z)\hookrightarrow\mathfrak X.
\]
\item[(2)] The \emph{scheme theoretic inverse image} $\Phi^{-1}(Z)=Z\times_{\mathfrak X,\Phi}\mathfrak X$ is a closed formal subscheme of $\mathfrak X$, with Cartesian square
\[
\begin{tikzcd}[column sep=1.8em]
\Phi^{-1}(Z) \arrow[r] \arrow[d, hook] & Z \arrow[d, hook] \\
\mathfrak X \arrow[r, "\Phi"'] & \mathfrak X\mathrlap{.}
\end{tikzcd}
\]
\item[(3)] The \emph{scheme theoretic union} of finitely many closed formal subschemes of $\mathfrak X$ is defined by the intersection of their defining ideals in $R$.
\end{enumerate}
\end{definition}

\begin{samepage}
\begin{lemma}\label{lem:image-ideals}
Let $Z=\mathrm{Spf}(R/I)$.
\begin{enumerate}
\item[(1)] The image $\Phi(Z)$ has defining ideal $(\Phi^\#)^{-1}(I)$.
\item[(2)] The inverse image $\Phi^{-1}(Z)$ has defining ideal $\Phi^\#(I)R$.
\end{enumerate}
\end{lemma}
\end{samepage}

\begin{proof}
The restriction $\Phi|_Z:Z\to\mathfrak X$ has pullback $R\xrightarrow{\Phi^\#}R\twoheadrightarrow R/I$. Its kernel $(\Phi^\#)^{-1}(I)$ defines the scheme theoretic image. The fibre product $\Phi^{-1}(Z) = Z \times_{\mathfrak{X},\Phi} \mathfrak{X}$ has coordinate ring $R/\Phi^\#(I)R$. Every ideal of a complete Noetherian adic ring is closed.
\end{proof}

A formal subtorus $T$ is \emph{$\Phi$-stable} if the restriction of $\Phi$ to $T$ factors through $T$.

\begin{lemma}\label{lem:subtorus-stability}
Let $M\subset\pi$ be primitive and let $\Phi:\mathfrak X_\pi\to\mathfrak X_\pi$ be an endomorphism of formal groups over $\Lambda$. Pullback of characters induces $U:\pi\to\pi$, so $\chi_{U(\gamma)}=\chi_\gamma\circ\Phi$. The following conditions are equivalent:
\begin{enumerate}
\item[(1)] $T_M$ is $\Phi$-stable.
\item[(2)] $\Phi(T_M)\subset T_M$.
\item[(3)] $\Phi^\#(I_M)\subset I_M$.
\item[(4)] $U(M)\subset M$.
\end{enumerate}
\end{lemma}

\begin{proof}
Conditions (1) and (2) agree by the definition of the scheme theoretic image. By Lemma~\ref{lem:image-ideals}(1), the inclusion in (2) is equivalent to $I_M\subset(\Phi^\#)^{-1}(I_M)$, which is (3).

Take $\gamma\in M$. Since $\Phi^\#([\gamma]-1)=[U\gamma]-1$, one has $\Phi^\#([\gamma]-1)\in I_M$ if and only if $U\gamma\in M$. By Lemma~\ref{lem:primitive-subtori}(2), these generators suffice to test (3).
\end{proof}

For a finite extension $L/K$, let $\mathcal O_L$ be its valuation ring and $\mathfrak m_L$ its maximal ideal. There is a natural bijection $\mathfrak X_\pi(\mathcal O_L)\xrightarrow{\sim}\mathrm{Hom}_{\mathrm{cont}}(\pi,1+\mathfrak m_L)$, $\zeta\mapsto(\gamma\mapsto\zeta^\#([\gamma]))$. A point $\zeta\in\mathfrak X_\pi(\mathcal O_L)$ is torsion if and only if the corresponding character has finite image. When $\zeta$ is torsion, it has $\ell$-power order, since the image is a finite continuous quotient of the pro-$\ell$ group $\pi$. Write $t_\zeta:\mathfrak X_{\pi,\mathcal O_L}\to\mathfrak X_{\pi,\mathcal O_L}$ for translation by $\zeta$.

\begin{lemma}\label{lem:torsion-translate-ideal}
Let $M\subset\pi$ be primitive, and let $L/K$ be a finite extension.
\begin{enumerate}
\item[(1)] For a torsion point $\zeta\in\mathfrak X_\pi(\mathcal O_L)$, the defining ideal of $t_\zeta(T_M)$ in $\mathcal O_L[[\pi]]$ is $([\gamma]-\zeta(\gamma)\mid\gamma\in M)$.
\item[(2)] Conversely, any continuous character $a:M\to1+\mathfrak m_L$ with finite image extends to a torsion point $\zeta\in\mathfrak X_\pi(\mathcal O_L)$. All such extensions define the same translate, with ideal $([\gamma]-a(\gamma)\mid\gamma\in M)$.
\end{enumerate}
\end{lemma}

\begin{proof}
For (1), inverse translation induces $(t_{\zeta^{-1}})^\#:\mathcal O_L[[\pi]]\to\mathcal O_L[[\pi]]$, $[\gamma]\mapsto\zeta(\gamma)^{-1}[\gamma]$. For $\gamma\in M$, this gives $(t_{\zeta^{-1}})^\#([\gamma]-1)=\zeta(\gamma)^{-1}([\gamma]-\zeta(\gamma))$, with $\zeta(\gamma)^{-1}\in\mathcal O_L^\times$. For (2), choose a decomposition $\pi=M\oplus M'$ and extend $a$ to a character $\zeta:\pi\to1+\mathfrak m_L$ by taking $\zeta|_{M'}=1$. By (1), the defining ideal depends only on $a$.
\end{proof}

\begin{lemma}\label{lem:coeff-ext}
Fix an algebraic closure of $K$. For $N\ge0$, let $L=K(\mu_{\ell^N})$, a finite extension of $K$, with valuation ring $\mathcal O_L$. Every torsion point $\zeta$ of $\mathfrak X_\pi$ killed by $\ell^N$ is in $\mathfrak X_\pi(\mathcal O_L)$. For every formal subtorus $T\subset\mathfrak X_\pi$, the translate $t_\zeta(T_{\mathcal O_L})$ is defined over $\mathcal O_L$.
\end{lemma}

\begin{proof}
The character defining the point takes values in $\mu_{\ell^N}$, and every $\ell$-power root of unity reduces to $1$ in residue characteristic $\ell$. The statement for translates $t_\zeta(T_{\mathcal O_L})$ follows from Lemma~\ref{lem:torsion-translate-ideal}(1).
\end{proof}

\begin{definition}\label{def:formal-quasilinear}
Let $Z\subset\mathfrak X_\pi$ be a closed formal subscheme.
\begin{enumerate}
\item[(1)] The subscheme $Z$ is \emph{quasilinear} if there exist a finite extension $L/K$ and an integer $r\ge0$ such that
\[
Z_{\mathcal O_L}=\bigcup_{j=1}^r t_{\zeta_j}(T_j)
\]
as a scheme theoretic union, where each $\zeta_j\in\mathfrak X_\pi(\mathcal O_L)$ is a torsion point and each $T_j\subset(\mathfrak X_\pi)_{\mathcal O_L}$ is a formal subtorus. The subscheme $Z$ is \emph{linear} if such a decomposition can be chosen with $r=1$.
\item[(2)] For an endomorphism $\Phi$ of $\mathfrak X_\pi$ over $\Lambda$, the subscheme $Z$ is \emph{$\Phi$-stable quasilinear} if the decomposition in (1) can be chosen with every $T_j$ stable under $\Phi_{\mathcal O_L}$. The subscheme $Z$ is \emph{$\Phi$-stable linear} if such a decomposition can be chosen with $r=1$.
\end{enumerate}
\end{definition}

\begin{remark}\label{rem:phi-stable-convention}
When the exponent matrix $U$ is scalar, every primitive sublattice $M \subset \pi$ satisfies $U(M) \subset M$, so every quasilinear closed formal subscheme is $\Phi$-stable quasilinear, as in \cite[Definition~1.1]{Jon00}.
\end{remark}

\section{Mixed weight formal linearity}\label{sec:formal-rigidity}

Let $\mathfrak{X} = \widehat{\mathbf{G}}_m^n = \mathrm{Spf}(R)$, where $R = \Lambda[[x_1, \dots, x_n]]$ and $q_i = 1 + x_i$. By the \emph{irreducible components} of an affine formal scheme $\mathrm{Spf}(D)$, we mean the reduced closed formal subschemes defined by the minimal prime ideals of $D$; see \cite[Section~2.4]{Nic09}.

\subsection{Diagonal power endomorphisms}\label{subsec:diag-endo}
Let $t_1,\dots,t_n$ be positive integers and let $U=\mathrm{diag}(\ell^{t_1},\dots,\ell^{t_n})$ act on $\mathbf Z_\ell^n$. By Proposition~\ref{prop:torus-morphisms}, the corresponding endomorphism $\Phi:\mathfrak X\to\mathfrak X$ over $\Lambda$ has pullback $\Phi^\#:R\to R$ satisfying
\begin{equation}\label{eq:diag-phi}
\Phi^\#(q_i)=q_i^{\ell^{t_i}}\quad(1\le i\le n), \qquad \Phi^\#(q^m)=q^{U(m)}\quad(m\in\mathbf Z_\ell^n).
\end{equation}

\begin{remark}\label{rem:finite-locally-free}
The endomorphism $\Phi$ is finite locally free of degree $\ell^{t_1+\dots+t_n}$. Indeed, $R$ is free over $\Phi^\#(R)$ with basis $q^\beta=\prod_{i=1}^n q_i^{\beta_i}$ for $0\le\beta_i<\ell^{t_i}$. In one variable, set $y:=(1+x)^{\ell^t}-1$. The polynomial $(1+x)^{\ell^t}-(1+y)$ is monic of degree $\ell^t$ in $x$, and all its nonleading coefficients are in $(\varpi,y)\subset\Lambda[[y]]$. Weierstrass division gives a basis of powers of $x$ of degree less than $\ell^t$. Replacing $x^j$ by $(1+x)^j$ is a unitriangular change of basis. Applying Weierstrass division successively in the $n$ variables and replacing powers of $x_i$ by powers of $q_i$ gives the stated basis.
\end{remark}

Let $\mathfrak Z\subset\mathfrak X$ be a closed formal subscheme with defining ideal $I_{\mathfrak Z}$. By Definition~\ref{def:formal-image} and Lemma~\ref{lem:image-ideals}, the scheme theoretic image is defined by contraction and the inverse image by extension.
\begin{equation}\label{eq:image-inverse-ideals}
 I_{\Phi(\mathfrak Z)}=(\Phi^\#)^{-1}(I_{\mathfrak Z}),
 \qquad I_{\Phi^{-1}(\mathfrak Z)}=R\cdot\Phi^\#(I_{\mathfrak Z}).
\end{equation}
For a closed formal subscheme $\mathfrak Y\subset\mathfrak X$ and an integer $c\ge1$, it follows that $\Phi^c(\mathfrak Y)\subset\mathfrak Y$ if and only if $(\Phi^\#)^c(I_{\mathfrak Y})\subset I_{\mathfrak Y}$. For integers $a>b\ge0$, taking $\mathfrak Y=\Phi^b(\mathfrak Z)$ and $c=a-b$ gives
\[
\Phi^a(\mathfrak Z)\subset\Phi^b(\mathfrak Z)
\quad\Longleftrightarrow\quad
(\Phi^\#)^{a-b}\bigl(I_{\Phi^b(\mathfrak Z)}\bigr)\subset I_{\Phi^b(\mathfrak Z)}.
\]

For a submodule $M\subset\mathbf Z_\ell^n$, let $I_M=(q^m-1\mid m\in M)$. If $M$ is primitive, $V(I_M)=T_M$ by Lemma~\ref{lem:primitive-subtori}(2).

\begin{definition}\label{def:saturation}
The \emph{$\ell$-saturation} of a submodule $M$ of a finite free $\mathbf Z_\ell$-module $\pi$ is $M^{\mathrm{sat}}=(M\otimes_{\mathbf Z_\ell}\mathbf Q_\ell)\cap\pi$ inside $\pi\otimes_{\mathbf Z_\ell}\mathbf Q_\ell$.
\end{definition}

Equivalently, $m\in\pi$ is in $M^{\mathrm{sat}}$ if and only if $\ell^r m\in M$ for some $r\ge 0$. In particular, $M^{\mathrm{sat}}$ is primitive.

\begin{remark}\label{rem:sat-u-stable}
Let $U\in\mathrm{End}_{\mathbf Z_\ell}(\pi)$ and suppose that $U(M)\subset M$. Then $U(M^{\mathrm{sat}})\subset M^{\mathrm{sat}}$. Indeed, if $m\in M^{\mathrm{sat}}$, then $\ell^r m\in M$ for some $r\ge 0$. Hence $\ell^r U(m)=U(\ell^r m)\in M$, so $U(m)\in M^{\mathrm{sat}}$.
\end{remark}

For a submodule $M\subset\mathbf Z_\ell^n$, let $d=\operatorname{rank}(M)$. By Smith normal form there is a basis $e_1,\dots,e_n$ of $\mathbf Z_\ell^n$ in which
\begin{equation}\label{eq:smith-L}
M^{\mathrm{sat}}=\bigoplus_{i=1}^d\mathbf Z_\ell e_i,\qquad M=\bigoplus_{i=1}^d\ell^{a_i}\mathbf Z_\ell e_i,\qquad a_i\ge0.
\end{equation}
In Lemmas~\ref{lem:binomial-model-properties} and~\ref{lem:nonprimitive-torsion-union}, write $q_i=[e_i]$ and $x_i=q_i-1$ in this basis.

\begin{lemma}\label{lem:binomial-model-properties}
\leavevmode
\begin{enumerate}
\item[(1)] The quotient $R/I_M$ is finite free over $\Lambda[[x_{d+1},\dots,x_n]]$.
\item[(2)] The formal subscheme $V(I_M)$ is flat over $\Lambda$, reduced, and equidimensional.
\end{enumerate}
\end{lemma}

\begin{proof}
The relations $q_i^{\ell^{a_i}}-1$ for $1\le i\le d$ generate $I_M$. Indeed, $u^a-1$ is divisible by $u-1$ for $a\in\mathbf Z_\ell$ and $u\in1+(x_1,\ldots,x_n)$, and $uv-1=u(v-1)+(u-1)$. Each polynomial $(1+x_i)^{\ell^{a_i}}-1$ is monic in $x_i$ and reduces to $x_i^{\ell^{a_i}}$ modulo $\varpi$. Successive Weierstrass division shows that the monomials $\prod_{i=1}^d x_i^{b_i}$ with $0\le b_i<\ell^{a_i}$ form a basis of $R/I_M$ over $\Lambda[[x_{d+1},\dots,x_n]]$.

This proves (1) and flatness over $\Lambda$. After inverting $\ell$, each defining polynomial has invertible derivative $\ell^{a_i}q_i^{\ell^{a_i}-1}$ in the quotient. The polynomial presentation from (1) therefore makes $(R/I_M)[1/\ell]$ finite étale over the domain $\Lambda[[x_{d+1},\dots,x_n]][1/\ell]$, hence reduced. Flatness over $\Lambda$ gives an injection $R/I_M\hookrightarrow(R/I_M)[1/\ell]$, proving reducedness of $R/I_M$. Finally, every minimal prime of $R/I_M$ contracts to zero in $\Lambda[[x_{d+1},\dots,x_n]]$ by flatness. The resulting finite integral extensions of domains all have dimension $n-d+1$, proving equidimensionality.
\end{proof}

\begin{lemma}\label{lem:nonprimitive-torsion-union}
Let $L/K$ be a finite extension containing the $\ell^{a_i}$-th roots of unity for $1\le i\le d$. The decomposition of $V(I_M)_{\mathcal O_L}$ into irreducible components is
\[
V(I_M)_{\mathcal O_L}
=\overset{\mathrm{sch}}{\bigcup_{\zeta_i\in\mu_{\ell^{a_i}}(L)}}
t_{(\zeta_1,\ldots,\zeta_d,1,\ldots,1)}\bigl((T_{M^{\mathrm{sat}}})_{\mathcal O_L}\bigr).
\]
\end{lemma}

\begin{proof}
In the coordinates of~\eqref{eq:smith-L}, the subtorus $T_{M^{\mathrm{sat}}}$ is defined by $q_1=\cdots=q_d=1$. For $\zeta_i\in\mu_{\ell^{a_i}}(L)$, translation by $(\zeta_1,\ldots,\zeta_d,1,\ldots,1)$ multiplies the first $d$ coordinates by $\zeta_1,\ldots,\zeta_d$, respectively. Thus
\[
t_{(\zeta_1,\ldots,\zeta_d,1,\ldots,1)}\bigl((T_{M^{\mathrm{sat}}})_{\mathcal O_L}\bigr)
=V(q_1-\zeta_1,\ldots,q_d-\zeta_d).
\]
The defining ideal $(q_1-\zeta_1,\dots,q_d-\zeta_d)$ in $\mathcal O_L[[x_1,\dots,x_n]]$ is prime. Indeed, substituting $x_i=\zeta_i-1$ for $i\le d$ induces an isomorphism
\[
\mathcal O_L[[x_1,\dots,x_n]]/(q_1-\zeta_1,\dots,q_d-\zeta_d)
\xrightarrow{\sim}\mathcal O_L[[x_{d+1},\dots,x_n]].
\]
Each $\zeta_i-1$ is in the maximal ideal of $\mathcal O_L$, so the substitution is continuous.

After inverting $\ell$, the factors of each $q_i^{\ell^{a_i}}-1=\prod_{\zeta\in\mu_{\ell^{a_i}}(L)}(q_i-\zeta)$ are pairwise comaximal. The Chinese remainder theorem thus expresses the localized quotient by $I_M$ as the product of the localized coordinate rings of the displayed translates. The formal scheme $V(I_M)_{\mathcal O_L}$ is flat over $\mathcal O_L$ by Lemma~\ref{lem:binomial-model-properties}, applied over $\mathcal O_L$. The coordinate ring of each translate is also flat over $\mathcal O_L$. Contracting the localized ideal decomposition gives
\[
 I_M\mathcal O_L[[x_1,\ldots,x_n]]
 =\bigcap_{(\zeta_1,\ldots,\zeta_d)\in\prod_{i=1}^d\mu_{\ell^{a_i}}(L)}
       (q_1-\zeta_1,\ldots,q_d-\zeta_d).
\]
These primes are pairwise incomparable after localization and hence before localization. They therefore define precisely the irreducible components of $V(I_M)_{\mathcal O_L}$ and yield the stated scheme theoretic union.
\end{proof}

\begin{proposition}\label{prop:images-subtori}
Let $M \subset \mathbf{Z}_\ell^n$ be a primitive sublattice and let $T_M \subset \mathfrak{X}$ be the corresponding formal subtorus from Lemma~\ref{lem:primitive-subtori}(1).
\begin{enumerate}
\item[(1)] $\Phi(T_M) = T_{U^{-1}(M)}$, where $U^{-1}(M) = \{m \in \mathbf{Z}_\ell^n \mid U(m) \in M\}$ is a primitive sublattice of $\mathbf{Z}_\ell^n$.
\item[(2)] $\Phi^{-1}(T_M) = V(I_{U(M)})$. Over the valuation ring of a finite extension $L/K$, the decomposition of $\Phi^{-1}(T_M)$ into irreducible components is
\[
\bigl(\Phi^{-1}(T_M)\bigr)_{\mathcal O_L}
=\overset{\mathrm{sch}}{\bigcup_{i=1}^r}
t_{\zeta_i}\!\left(\bigl(T_{(U(M))^{\mathrm{sat}}}\bigr)_{\mathcal O_L}\right),
\]
where $\zeta_1,\ldots,\zeta_r\in\mathfrak X(\mathcal O_L)$ are torsion points.
\end{enumerate}
\end{proposition}

\begin{proof}
(1) The homomorphism $U$ induces the commutative diagram
\[
\begin{tikzcd}[column sep=1.8em]
\mathbf Z_\ell^n \arrow[r, "U"] \arrow[d, two heads] & \mathbf Z_\ell^n \arrow[d, two heads] \\
\mathbf Z_\ell^n/U^{-1}(M) \arrow[r, hook, "\bar U"'] & \mathbf Z_\ell^n/M\mathrlap{,}
\end{tikzcd}
\]
where the vertical arrows are the quotient maps. The map $\bar U$ is injective by definition of $U^{-1}(M)$ and has finite cokernel because $U$ does. The quotient $\mathbf Z_\ell^n/M$ is torsion free, so $U^{-1}(M)$ is primitive. By Lemma~\ref{lem:primitive-subtori}(2), the homomorphism of coordinate rings induced by $\bar U$ is
\[
\Lambda[[\mathbf Z_\ell^n/U^{-1}(M)]]\longrightarrow\Lambda[[\mathbf Z_\ell^n/M]].
\]
By Remark~\ref{rem:finite-locally-free}, this homomorphism is finite free in Smith normal form. It gives an injection $R/I_{U^{-1}(M)}\hookrightarrow R/I_M$ induced by $\Phi^\#$. Hence $(\Phi^\#)^{-1}(I_M)=I_{U^{-1}(M)}$, and \eqref{eq:image-inverse-ideals} gives $\Phi(T_M)=T_{U^{-1}(M)}$.

For (2), the formula for the scheme theoretic inverse image in \eqref{eq:image-inverse-ideals} gives
\[
I_{\Phi^{-1}(T_M)} = R \cdot \Phi^\#(I_M) = (q^{U(m)} - 1 \mid m \in M) = I_{U(M)},
\]
and the decomposition of $V(I_{U(M)})_{\mathcal O_L}$ follows from Lemma~\ref{lem:nonprimitive-torsion-union} applied to $U(M)$.
\end{proof}

\subsection{Classification of stable prime ideals}\label{subsec:stable-prime-ideals}
Let $\mathfrak n=(x_1,\dots,x_n)$ be the augmentation ideal and let $\mathfrak m=(\ell)+\mathfrak n$. Since $(\ell)$ is a power of $(\varpi)$ in $\Lambda$, the $\mathfrak m$-adic and $(\varpi,\mathfrak n)$-adic topologies agree. In particular, $R$ is $\mathfrak m$-adically complete.

\begin{lemma}\label{lem:m-adic-contraction}
One has
\[
\Phi^\#(\mathfrak m)\subset\mathfrak m,
\qquad
\Phi^\#(\mathfrak n)\subset\ell\mathfrak n+\mathfrak n^2=\mathfrak m\mathfrak n\subset\mathfrak m^2.
\]
In particular, $(\Phi^\#)^N(g)\to0$ in the $\mathfrak m$-adic topology for every $g\in\mathfrak n$.
\end{lemma}

\begin{proof}
Each $\Phi^\#(x_i)=(1+x_i)^{\ell^{t_i}}-1$ has linear term $\ell^{t_i}x_i$, and its higher degree terms lie in $(x_i)^2$. Thus $\Phi^\#(\mathfrak n)\subset\mathfrak m\mathfrak n$. Since $\Phi^\#(\ell)=\ell$, we also have $\Phi^\#(\mathfrak m)\subset\mathfrak m$. Induction gives $(\Phi^\#)^N(\mathfrak n)\subset\mathfrak m^N\mathfrak n$ for all $N\ge 1$.
\end{proof}

\begin{lemma}\label{lem:pass-identity}
Let $\mathfrak{p} \subset R$ be a prime ideal such that $\ell \notin \mathfrak{p}$ and $\Phi^\#(\mathfrak{p}) \subset \mathfrak{p}$. Then $\mathfrak{p} \subset \mathfrak{n} = (x_1, \dots, x_n)$, so $\mathrm{Spf}(R/\mathfrak{p}) \subset \mathfrak{X}$ passes through the identity section.
\end{lemma}

\begin{proof}
For $f\in\mathfrak p$, Lemma~\ref{lem:m-adic-contraction} gives $(\Phi^\#)^N(f)=f(0)+(\Phi^\#)^N(f-f(0))\to f(0)$. The ideal $\mathfrak p$ is closed and contains every iterate of $f$, so $f(0)\in\mathfrak p\cap\Lambda$. Since $\ell\notin\mathfrak p$ and $\Lambda$ is a discrete valuation ring of residue characteristic $\ell$, this intersection is zero. Thus $f\in\mathfrak n$.
\end{proof}

Let $\mathfrak p\subset R$ be a $\Phi^\#$-stable prime ideal with $\ell\notin\mathfrak p$. By Lemma~\ref{lem:pass-identity}, $\mathfrak p\subset\mathfrak n$, so the image of $\mathfrak n R[1/\ell]$ in $(R/\mathfrak p)[1/\ell]$ is the maximal ideal of the identity point. Completing $R[1/\ell]$ at $\mathfrak n R[1/\ell]$ gives $K[[x_1,\dots,x_n]]$, which contains the formal power series $\log(q_i)=\sum_{m\ge1}(-1)^{m+1}x_i^m/m$.

\Needspace{4\baselineskip}
\begin{lemma}\label{lem:ext-contract-identity}
Let $J\subset R$ be a radical ideal whose minimal primes are contained in $\mathfrak n$.
\begin{enumerate}
\item[(1)] The ideal $J$ is $\ell$-saturated.
\item[(2)] The extended ideal $J K[[x_1,\dots,x_n]]$ is radical.
\item[(3)] $J K[[x_1,\dots,x_n]]\cap R=J$.
\end{enumerate}
\end{lemma}

\begin{proof}
(1) Write $J=\bigcap_{\nu=1}^r\mathfrak p_\nu$ as the intersection of its minimal primes. Since $\mathfrak p_\nu\subset\mathfrak n$ and $\ell\notin\mathfrak n$, we have $\ell\notin\mathfrak p_\nu$ for every $\nu$. Thus $\ell f\in J$ implies $f\in\mathfrak p_\nu$ for all $\nu$, proving $\ell$-saturation.

(2) Since $\ell\notin\mathfrak n$, the local ring $(R[1/\ell])_{\mathfrak nR[1/\ell]}$ equals $R_{\mathfrak n}$ and has completion $K[[x_1,\dots,x_n]]$. The ring $R$ is excellent since it is complete Noetherian local, so $R_{\mathfrak n}/JR_{\mathfrak n}$ is excellent as well \cite[Tags~07QW and~07QU]{Stacks}. Since $J$ is radical, $R_{\mathfrak n}/JR_{\mathfrak n}$ is reduced. Consequently, the completion $K[[x_1,\dots,x_n]]/JK[[x_1,\dots,x_n]]$ is reduced \cite[Tags~0AH2 and~07QK]{Stacks}.

(3) Completion is faithfully flat over $R_{\mathfrak n}$, so $JK[[x_1,\dots,x_n]]\cap R_{\mathfrak n}=JR_{\mathfrak n}$. Since $\mathfrak p_\nu\subset\mathfrak n$, we have $\mathfrak p_\nu R_{\mathfrak n}\cap R=\mathfrak p_\nu$ for every $\nu$. Localization commutes with finite intersections, so $JR_{\mathfrak n}\cap R=J$.
\end{proof}

\begin{lemma}\label{lem:total-invariance}
Let $\mathfrak p\subset R$ be a prime ideal with $\Phi^\#(\mathfrak p)\subset\mathfrak p$. Then $(\Phi^\#)^{-1}(\mathfrak p)=\mathfrak p$. Consequently, $\Phi^\#$ induces injective endomorphisms of $R/\mathfrak p$ and $\mathrm{Frac}(R/\mathfrak p)$.
\end{lemma}

\begin{proof}
Stability gives $\mathfrak p\subset(\Phi^\#)^{-1}(\mathfrak p)$. The injection $R/(\Phi^\#)^{-1}(\mathfrak p)\hookrightarrow R/\mathfrak p$ induced by $\Phi^\#$ is finite by Remark~\ref{rem:finite-locally-free}. Hence the two quotients have the same dimension. A strict inclusion of primes would give a strict inequality between these dimensions. Hence $(\Phi^\#)^{-1}(\mathfrak p)=\mathfrak p$.
\end{proof}

\begin{lemma}\label{lem:eigenvalue-obstruction}
Let $\mathfrak p\subset R$ be a $\Phi^\#$-stable prime ideal with $\ell\notin\mathfrak p$, and let $A=R/\mathfrak p$. Let $\varphi$ denote the injective endomorphism of $\mathrm{Frac}(A)$ induced by $\Phi^\#$. If $w\in\mathrm{Frac}(A)^\times$ and $a\in\mathbf Z$ satisfy $\varphi(w)=\ell^a w$, then $a=0$.
\end{lemma}

\begin{proof}
Let $\widetilde A$ be the normalization of $A$ in $\mathrm{Frac}(A)$. The ring $A$ is excellent, so $\widetilde A$ is finite over $A$ \cite[Tag~07QV]{Stacks}. The map $\varphi$ carries integral equations over $A$ to integral equations over $\varphi(A)\subset A$, hence preserves $\widetilde A$. Since $A$ is finite over $\varphi(A)$ by Remark~\ref{rem:finite-locally-free}, $\widetilde A$ is finite over $\varphi(A)$ and hence over $\varphi(\widetilde A)$.

Since $\ell \notin \mathfrak{p}$ and $\ell$ lies in the maximal ideal of $A$, the set of height one prime ideals of $\widetilde{A}$ containing $\ell$ is finite and nonempty. Because $\varphi(\widetilde{A}) \subset \widetilde{A}$ is a finite integral extension of normal domains and $\varphi(\ell) = \ell$, the going down and incomparability theorems show that contraction along $\varphi$ preserves this finite set. Choosing a prime $\mathfrak{P}$ in this finite set that is periodic under contraction gives $(\varphi^N)^{-1}(\mathfrak{P}) = \mathfrak{P}$ for some integer $N \ge 1$. Let $v$ be the normalized discrete valuation of the discrete valuation ring $\widetilde{A}_{\mathfrak{P}}$. Then $v \circ \varphi^N = e v$ on $\mathrm{Frac}(A)^\times$ for some integer $e \ge 1$. Evaluating $v\circ\varphi^N=ev$ at $\ell$ and using $\varphi^N(\ell) = \ell$ with $v(\ell) > 0$ yields $e = 1$, so $v(\varphi^N(f)) = v(f)$ for all $f \in \mathrm{Frac}(A)^\times$. Applying $v$ to the iterated eigenvalue relation $\varphi^N(w) = \ell^{Na} w$ now gives
\[
v(w) = v(\varphi^N(w)) = v(w) + N a v(\ell),
\]
and since $N \ge 1$ and $v(\ell) > 0$, we conclude that $a = 0$.
\end{proof}

For a complete Noetherian local $\Lambda$-algebra $A$, the module $\Omega^{\mathrm{cont}}_{A/\Lambda}$ of continuous differentials is the completion of $\Omega^1_{A/\Lambda}$ with respect to the maximal ideal of $A$; see \cite[Chapter~0, Section~20.7]{Gro64}.

\begin{lemma}\label{lem:diff-alg-const}
Let $\mathfrak p\subset R$ be a $\Phi^\#$-stable prime ideal with $\ell\notin\mathfrak p$. Let $A=R/\mathfrak p$.
\begin{enumerate}
\item[(1)] After reordering the pairs $(q_i,t_i)$, the differentials $dq_j/q_j$ for $1\le j\le\dim A-1$ form a basis of $\Omega^{\mathrm{cont}}_{A/\Lambda}\otimes_A\mathrm{Frac}(A)$ over $\mathrm{Frac}(A)$.
\item[(2)] If $h\in\mathrm{Frac}(A)$ satisfies $dh=0$, then $h$ is algebraic over $K$.
\end{enumerate}
\end{lemma}

\begin{proof}
(1) Let $d := \dim A - 1$, and choose lifts $y_1,\dots,y_d\in\mathfrak n A$ of a system of parameters of $A/\varpi A$. Let $B_0=\Lambda[[y_1,\dots,y_d]]$. The continuous $\Lambda$-algebra homomorphism $B_0\to A$, given by evaluation at these lifts, is finite and local by completeness and the formal finiteness criterion \cite[Tag~0394]{Stacks}. Both $B_0$ and $A$ are integral domains of dimension $d + 1$, so the homomorphism $B_0\to A$ is injective.
Let $L=\mathrm{Frac}(B_0)$ and $F=\mathrm{Frac}(A)$. We show that $dy_1,\ldots,dy_d$ form a basis of $\Omega^{\mathrm{cont}}_{A/\Lambda}\otimes_A F$. For each $i$, let $P_i\in L[T]$ be the minimal polynomial of the image of $x_i$ in $F$. Since $F/L$ is separable, $P_i'(x_i)\ne0$. Differentiating $P_i(x_i)=0$ expresses $dx_i$ as an $F$-linear combination of $dy_1,\ldots,dy_d$, since the differential of each coefficient is an $L$-linear combination of the $dy_j$. The $dx_i$ generate $\Omega^{\mathrm{cont}}_{A/\Lambda}\otimes_A F$, so the $dy_j$ span it.

To prove linear independence, extend $\partial/\partial y_j$ from $L$ to a derivation $D_j$ of $F$ by separability. Since $D_j(B_0)\subseteq B_0$ and $A$ is finite over $B_0$, clearing denominators on a finite set of $B_0$-module generators gives $0\ne c_j\in A$ with $c_jD_j(A)\subseteq A$. Write $\mathfrak m_A$ for the maximal ideal of $A$. The derivation $c_jD_j:A\to A$ is continuous because it sends $\mathfrak m_A^m$ into $\mathfrak m_A^{m-1}$ for $m\ge1$. The resulting $F$-linear functional on $\Omega^{\mathrm{cont}}_{A/\Lambda}\otimes_A F$, divided by $c_j$, sends $dy_i$ to $\delta_{ij}$. Thus the $dy_j$ form a basis. Since the logarithmic differentials $dq_i/q_i = dx_i/q_i$ for $1 \le i \le n$ span $\Omega^{\mathrm{cont}}_{A/\Lambda} \otimes_A \mathrm{Frac}(A)$, we reorder the coordinates so that $dq_1/q_1, \dots, dq_{\dim A - 1}/q_{\dim A - 1}$ form a basis over $\mathrm{Frac}(A)$.

(2) Let $h \in \mathrm{Frac}(A)$ satisfy $dh = 0$, and let $T^m + a_1 T^{m-1} + \dots + a_m \in \mathrm{Frac}(B_0)[T]$ be the minimal polynomial of $h$ over $\mathrm{Frac}(B_0)$. Differentiating the relation $h^m + a_1 h^{m-1} + \dots + a_m = 0$ in $\Omega^{\mathrm{cont}}_{A/\Lambda} \otimes_A \mathrm{Frac}(A)$ and using $dh = 0$ yields
\[
h^{m-1} da_1 + h^{m-2} da_2 + \dots + da_m = 0.
\]
Expand each $da_k$ in the basis $dy_1,\dots,dy_d$ over $\mathrm{Frac}(B_0)$. By (1), the same differentials are independent over $\mathrm{Frac}(A)$. For $1\le r\le d$, the coefficient of $dy_r$ gives $\sum_{k=1}^m(\partial a_k/\partial y_r)h^{m-k}=0$. Since $1,h,\dots,h^{m-1}$ are linearly independent over $\mathrm{Frac}(B_0)$, every $\partial a_k/\partial y_r$ vanishes. Thus $da_k=0$ for every $k$. Embed $\mathrm{Frac}(B_0)$ into the iterated Laurent series field $K((y_1))\cdots((y_d))$. In characteristic zero a Laurent series with zero derivative in its last variable has only a constant term in that variable. Applying the vanishing of derivatives in the last variable successively to $y_d,\dots,y_1$ gives
\[
\ker\left(d : \mathrm{Frac}(B_0) \longrightarrow \bigoplus_{r=1}^{\dim A - 1} \mathrm{Frac}(B_0)\, dy_r\right) = K.
\]
Thus $a_k \in K$ for all $1 \le k \le m$, so $h$ is algebraic over $K$.
\end{proof}

\begin{lemma}\label{lem:fixed-elements-algebraic}
Let $A$ and $\varphi$ be as in Lemma~\ref{lem:eigenvalue-obstruction}. Every element of $\mathrm{Frac}(A)^\varphi$ is algebraic over $K$.
\end{lemma}

\begin{proof}
The endomorphism $\varphi$ induces a $\varphi$-semilinear endomorphism $\varphi^*$ of $\Omega^{\mathrm{cont}}_{A/\Lambda}\otimes_A\mathrm{Frac}(A)$, with
\[
\varphi^*(a\,db)=\varphi(a)\,d(\varphi(b)),
\qquad
\varphi^*\!\left(\frac{dq_j}{q_j}\right)
=\ell^{t_j}\frac{dq_j}{q_j}.
\]

For $f\in\mathrm{Frac}(A)^\varphi$, expand $df=\sum_j f_j\,dq_j/q_j$ in the basis of Lemma~\ref{lem:diff-alg-const}(1). Comparing coefficients in $\varphi^*(df)=df$ gives $\varphi(f_j)=\ell^{-t_j}f_j$. All $t_j$ are positive, so Lemma~\ref{lem:eigenvalue-obstruction} forces every $f_j$ to vanish. Now $df=0$, and Lemma~\ref{lem:diff-alg-const}(2) shows that $f$ is algebraic over $K$.
\end{proof}

\begin{lemma}\label{lem:log-diff-weight-pure}
Let $A$ and $\varphi$ be as in Lemma~\ref{lem:fixed-elements-algebraic}, and order the coordinates as in Lemma~\ref{lem:diff-alg-const}(1). For $\dim A\le i\le n$, one has
\[
\frac{dq_i}{q_i}=\sum_{\substack{1\le j\le\dim A-1\\ t_j=t_i}}w_j^{(i)}\frac{dq_j}{q_j},
\qquad w_j^{(i)}\in\mathrm{Frac}(A)^\varphi.
\]
In particular, all the coefficients $w_j^{(i)}$ are algebraic over $K$.
\end{lemma}

\begin{proof}
Write $dq_i/q_i=\sum_{j=1}^{\dim A-1}w_j^{(i)}\,dq_j/q_j$ in the basis of Lemma~\ref{lem:diff-alg-const}(1). Applying $\varphi^*$ and comparing coefficients of $dq_j/q_j$, using $\varphi^*(dq_k/q_k)=\ell^{t_k}\,dq_k/q_k$, gives $\varphi(w_j^{(i)})=\ell^{t_i-t_j}w_j^{(i)}$. A nonzero coefficient therefore has $t_i=t_j$ by Lemma~\ref{lem:eigenvalue-obstruction}. The same identity shows that this coefficient is fixed by $\varphi$, hence algebraic over $K$ by Lemma~\ref{lem:fixed-elements-algebraic}.
\end{proof}

\begin{lemma}[Constant logarithmic tangent spaces]\label{lem:const-log-tangent}
Let $P\subset(x_1,\dots,x_n)$ be a prime ideal of $R$. Let $A=R/P$, $F=\mathrm{Frac}(A)$, and $d=\dim A-1$. Let $W\subset K^n$ be an $(n-d)$-dimensional subspace. Suppose that the $F$-linear map
\[
F^n\longrightarrow\Omega^{\mathrm{cont}}_{A/\Lambda}\otimes_A F,
\qquad
(a_i)\longmapsto\sum_{i=1}^n a_i\,\frac{dq_i}{q_i},
\]
has kernel $W\otimes_K F$.
\begin{enumerate}
\item[(1)] For every basis $\omega^{(1)},\dots,\omega^{(n-d)}$ of $W$, the ideal $P K[[x_1,\dots,x_n]]$ is generated by the linear forms $z_k=\sum_{i=1}^n\omega_i^{(k)}\log(q_i)$ for $1\le k\le n-d$.
\item[(2)] The local ring of $A[1/\ell]$ at the identity is regular.
\end{enumerate}
\end{lemma}

\begin{proof}
Let $\widehat P=PK[[x_1,\dots,x_n]]$. For $f\in P$, the image of $df=\sum_i(q_i\partial f/\partial q_i)\,dq_i/q_i$ in $\Omega^{\mathrm{cont}}_{A/\Lambda}\otimes_A F$ vanishes. Thus the logarithmic gradient $(q_i\partial f/\partial q_i\bmod P)_i$ lies in $W\otimes_K F$. For $\lambda$ in the annihilator $W^\perp\subset(K^n)^\vee$, the $K$-derivation $\theta_\lambda:R[1/\ell]\to R[1/\ell]$, $f\mapsto\sum_i\lambda_i q_i\,\partial f/\partial q_i$, therefore sends $P$ into $PR[1/\ell]$.

In the coordinates $l_i=\log(q_i)$ the continuous extension of $\theta_\lambda$ is the constant coefficient derivation $\sum_i\lambda_i\,\partial/\partial l_i$. It preserves $\widehat P$ and annihilates each $z_k$. Extend $z_1,\dots,z_{n-d}$ to a linear coordinate system by forms $s_1,\dots,s_d$ dual to a basis of $W^\perp$. The derivations associated with this basis are precisely $\partial/\partial s_1,\dots,\partial/\partial s_d$.

For any ideal $J \subset K[[s, z]] = K[[s_1, \dots, s_d, z_1, \dots, z_{n-d}]]$ satisfying $\partial J/\partial s_j \subset J$ for all $1 \le j \le d$, Taylor expansion in $s$ gives
\begin{equation}\label{eq:taylor-ideal}
J = (J \cap K[[z]]) K[[s, z]].
\end{equation}
Indeed, for $f\in J$, each derivative $\partial_s^\alpha f$ is in $J$, so the Taylor expansion $f(s+u,z)=\sum_\alpha(\partial_s^\alpha f)(s,z)u^\alpha/\alpha!$ is in $JK[[s,z,u]]$. Setting $u=-s$ gives $f(0,z)\in J$. Expanding $(\partial_s^\alpha f)(s+u,z)$ and setting $u=-s$ gives $(\partial_s^\alpha f)(0,z)/\alpha!\in J\cap K[[z]]$ for every $\alpha$. These are the coefficients of $f$ as a power series in $s$. The ideal generated by $J\cap K[[z]]$ is closed, so it contains $f$.

Let $J_0=\widehat P\cap K[[z]]$. By \eqref{eq:taylor-ideal}, $K[[s,z]]/\widehat P$ is the power series ring over $K[[z]]/J_0$ in $d$ variables. Its dimension is $d$, since it is the completion of the local ring of $A[1/\ell]$ at the identity, of dimension $n-\mathrm{ht}(P)=d$ by catenarity. Thus $K[[z]]/J_0$ is zero dimensional. It is local and reduced by Lemma~\ref{lem:ext-contract-identity}(2), so it is a field. Hence $J_0=(z_1,\dots,z_{n-d})$, proving (1). The completed local ring at the identity is now $K[[s_1,\dots,s_d]]$, giving (2).
\end{proof}

\begin{lemma}\label{lem:equal-weight-binomial}
Fix an integer $t\ge 1$, and consider the $\Lambda$-algebra endomorphism $\Sigma_t^\#:R\to R$, $q_i\mapsto q_i^{\ell^t}$. Let $P \subset R$ be a prime ideal such that $\ell \notin P$ and $\Sigma_t^\#(P) \subset P$. If $K$ is algebraically closed in $\mathrm{Frac}(R/P)$, then there exists a unique primitive sublattice $M \subset \mathbf{Z}_\ell^n$ such that
\[
P = (q^m - 1 \mid m \in M) \subset R.
\]
\end{lemma}

\begin{proof}
Let $A=R/P$ and $F=\mathrm{Frac}(A)$.

\textit{Step 1. Geometric integrality and linear forms at the identity.} By Lemma~\ref{lem:total-invariance}, $\Sigma_t^\#$ induces an endomorphism $\sigma:F\to F$ fixing $K$. Lemma~\ref{lem:fixed-elements-algebraic} shows that every element of $F^\sigma$ is algebraic over $K$. Since $K$ is algebraically closed in $F$, we have
\begin{equation}\label{eq:fixed-field-F}
F^\sigma = K.
\end{equation}
Because $K$ has characteristic zero and is algebraically closed in $F$, the field extension $F/K$ is regular, so for every finite extension $K'/K$ the tensor product $K' \otimes_K F$ is a field. Since $K' \otimes_K A[1/\ell]$ embeds into $K' \otimes_K F$, it follows that
\begin{equation}\label{eq:geom-integral-A}
K' \otimes_K A[1/\ell] \text{ is an integral domain for every finite extension } K'/K.
\end{equation}

Let $d = \dim A - 1$ and $\mathfrak{a} = (x_1, \dots, x_n)A$. Lemma~\ref{lem:pass-identity} gives $P \subset (x_1, \dots, x_n)$. Reorder the coordinates so that $dq_1/q_1,\dots,dq_d/q_d$ form the basis of Lemma~\ref{lem:diff-alg-const}(1). For each $i>d$, Lemma~\ref{lem:log-diff-weight-pure} expresses $dq_i/q_i$ in this basis with coefficients $w_j^{(i)}\in F^\sigma=K$ by \eqref{eq:fixed-field-F}. Let $W\subset K^n$ be spanned by the vectors $e_i-\sum_{j=1}^d w_j^{(i)}e_j$ for $i>d$, where the $e_i$ are the standard basis vectors. These vectors form an $F$-basis of the kernel of the logarithmic differential map $F^n\to\Omega^{\mathrm{cont}}_{A/\Lambda}\otimes_A F$, $(a_i)\mapsto\sum_i a_i\,dq_i/q_i$. For $i>d$, write $z_{i-d}=\log(q_i)-\sum_{j=1}^d w_j^{(i)}\log(q_j)$. Lemma~\ref{lem:const-log-tangent} gives
\begin{equation}\label{eq:P-hat-linear}
\widehat{P} := P K[[x_1, \dots, x_n]] = (z_1, \dots, z_{n-d}),
\end{equation}
and shows that the local ring $A_{\mathfrak a}=(A[1/\ell])_{\mathfrak a A[1/\ell]}$ is regular.

\textit{Step 2. Regularity of $A[1/\ell]$ via the logarithmic Jacobian.} Let $h := n - d = \mathrm{ht}(P)$, choose generators $f_1, \dots, f_r$ of $P$, and let $J_0 \subset A$ be the ideal generated by the $h \times h$ minors of the logarithmic Jacobian matrix
\[
\left(q_j \frac{\partial f_a}{\partial q_j} \bmod P\right)_{1 \le a \le r,\, 1 \le j \le n},
\]
with the convention that the $0 \times 0$ minor for $h = 0$ is $1$. Let $C$ be a quotient of a formal power series ring over $\Lambda$, let $\mathfrak n_C\subset C[1/\ell]$ be a maximal ideal, and let $D=C/(\mathfrak n_C\cap C)$. Then $D$ is finite over $\Lambda$, and $\kappa(\mathfrak n_C)=D[1/\ell]$ is a finite extension of $K$. Indeed, $D$ is a complete Noetherian local domain in which $\ell \ne 0$, so formal Noether normalization provides a finite injective homomorphism $\Lambda[[T_1, \dots, T_e]] \hookrightarrow D$ with $e = \dim(D) - 1$. Inverting $\ell$ makes the field $D[1/\ell]$ integral over $\Lambda[[T_1, \dots, T_e]][1/\ell]$, which forces $\Lambda[[T_1, \dots, T_e]][1/\ell]$ to be a field and hence $e = 0$.

Now let $\mathfrak{m} \subset A[1/\ell]$ be any maximal ideal, let $\widetilde{\mathfrak{m}}$ be its preimage in $R[1/\ell]$, and let $\mathfrak{p} := \widetilde{\mathfrak{m}} \cap R$. Applying the finiteness property just proved to $R$ shows that $R/\mathfrak{p}$ is finite over $\Lambda$ of Krull dimension $1$ and that $L := \kappa(\mathfrak{m}) = (R/\mathfrak{p})[1/\ell]$ is a finite extension of $K$. By catenarity of the regular local ring $R$, we have $\mathrm{ht}(\mathfrak{p}) = n$, so $R[1/\ell]_{\widetilde{\mathfrak{m}}}$ is a regular local ring of dimension $n$ with residue field $L$. Let $a_j$ be the image of $q_j$ in $L$. Because $R/\mathfrak p$ is finite and local over $\Lambda$, each $a_j-1$ is in the maximal ideal of $\mathcal O_L$. After the finite separable extension $L/K$, the point has an $L$-rational lift given by these coordinates. Base change along $K\to L$ gives a finite étale morphism $\mathrm{Spec}(R[1/\ell]\otimes_K L)\to\mathrm{Spec}(R[1/\ell])$. At the chosen lift its map of local rings induces the identity on the residue field $L$, hence an isomorphism on completions. Over $\mathcal O_L$, translation by $a_j-1$ identifies the completion at the lift with $L[[q_1-a_1,\dots,q_n-a_n]]$. Thus the images of $q_j-a_j$ form a regular system of parameters, and multiplication by the units $q_j$ preserves the rank of the Jacobian. The Jacobian criterion for regularity in the complete regular local ring $\widehat{R[1/\ell]_{\widetilde{\mathfrak{m}}}}$ therefore shows that
\begin{equation}\label{eq:jacobian-criterion}
A[1/\ell]_{\mathfrak{m}} \text{ is regular} \quad \Longleftrightarrow \quad J_0 A[1/\ell] \not\subset \mathfrak{m}.
\end{equation}

Since $\Sigma_t^\#(P) \subset P$, we write $\Sigma_t^\#(f_a) = \sum_{b=1}^r h_{ab} f_b$ with $h_{ab} \in R$. Applying $q_j \partial/\partial q_j$ to this relation and reducing modulo $P$ yields
\[
\ell^t \sigma\left(\overline{q_j \frac{\partial f_a}{\partial q_j}}\right) = \sum_{b=1}^r \overline{h}_{ab}\, \overline{q_j \frac{\partial f_b}{\partial q_j}} \qquad \text{in } A.
\]
By the Cauchy--Binet formula for the $h \times h$ minors, we obtain $\ell^{th} \sigma(J_0) \subset J_0$. Let $J := (J_0 : \ell^\infty) = \{x \in A \mid \ell^m x \in J_0 \text{ for some } m \ge 0\}$ be the $\ell$-saturation of $J_0$ in $A$. Then $J$ is $\ell$-saturated and satisfies $\sigma(J) \subset J$.

Every proper $\ell$-saturated ideal $I\subset A$ preserved by $\sigma$ lies in $\mathfrak a$. Indeed, write $b\in I$ as $b=c+g$ with $c\in\Lambda$ and $g\in\mathfrak a$. By Lemma~\ref{lem:m-adic-contraction}, $\sigma^N(g)\to0$. The ideal $I$ is closed, so $c=\lim_N\sigma^N(b)\in I$. If $c\ne0$, then $c\Lambda$ contains a power of $\ell$, and saturation gives $1\in I$. Thus $c=0$ and $b\in\mathfrak a$.

Applying this containment in $\mathfrak a$ to the $\sigma$-stable $\ell$-saturated ideal $J$ shows that $J=A$, since the local ring $A_{\mathfrak a}$ is regular by \eqref{eq:P-hat-linear}, so \eqref{eq:jacobian-criterion} gives $J_0\not\subset\mathfrak a$ and hence $J\not\subset\mathfrak a$. It follows that $J_0A[1/\ell]=A[1/\ell]$, and a second application of \eqref{eq:jacobian-criterion} proves regularity at every maximal ideal of $A[1/\ell]$.

\textit{Step 3. Integrality of $A\mathbin{\widehat\otimes}_{\Lambda}A$.} Let $A_K := A[1/\ell]$, $B := A \mathbin{\widehat{\otimes}}_{\Lambda} A$, and $B_K := B[1/\ell]$. Because $A$ is $\ell$-torsion free and hence flat over $\Lambda$, and is pseudocompact as a complete local $\Lambda$-module, $A$ is topologically free over $\Lambda$ by \cite[Exp.~VII\textsubscript{B}, Corollary~0.3.8]{DG70}. For either of the two canonical homomorphisms $A\to B$, a pseudobasis of the other factor identifies $B$ with a direct product of copies of $A$. Thus both homomorphisms are flat by the same corollary in \cite{DG70}. In particular, $B$ is $\Lambda$-flat and both maps $A_K\to B_K$ are flat.

Consider the first projection $p_1 : \mathrm{Spec}(B_K) \to \mathrm{Spec}(A_K)$. For any closed point $y \in \mathrm{Spec}(A_K)$ with contraction $\mathfrak{q} \subset A$, the quotient $A/\mathfrak{q}$ is finite over $\Lambda$ with $(A/\mathfrak{q})[1/\ell] = \kappa(y)$. Because $A/\mathfrak{q}$ is a finite $\Lambda$-module, the completed tensor product $(A/\mathfrak{q}) \mathbin{\widehat{\otimes}}_{\Lambda} A$ coincides with the ordinary tensor product $(A/\mathfrak{q}) \otimes_{\Lambda} A$, and inverting $\ell$ yields
\[
B_K \otimes_{A_K} \kappa(y) \simeq (B/\mathfrak{q}B)[1/\ell] \simeq \bigl((A/\mathfrak{q}) \otimes_{\Lambda} A\bigr)[1/\ell] \simeq \kappa(y) \otimes_K A_K.
\]
By \eqref{eq:geom-integral-A}, this fibre $p_1^{-1}(y) = \mathrm{Spec}(\kappa(y) \otimes_K A_K)$ is integral for every closed point $y \in \mathrm{Spec}(A_K)$.

Fix a maximal ideal $\mathfrak{m} \subset B_K$. The residue field $\kappa(\mathfrak{m})$ is a finite extension of $K$ because $B$ is a quotient of a formal power series ring over $\Lambda$. Thus the image of $A_K$ in $\kappa(\mathfrak{m})$ is a finite dimensional integral domain over $K$ and hence a field, so $\mathfrak{m}_1 := \mathfrak{m} \cap A_K$ is a maximal ideal of $A_K$ whose residue field $\kappa(\mathfrak{m}_1)$ is a finite extension of $K$. In the flat local homomorphism $(A_K)_{\mathfrak{m}_1} \to (B_K)_{\mathfrak{m}}$, the base $(A_K)_{\mathfrak{m}_1}$ is regular and the closed fibre is a localization of $\kappa(\mathfrak{m}_1) \otimes_K A_K$, which is finite étale over the regular ring $A_K$ and hence regular. By the flat local criterion for regularity \cite[Tag~031E]{Stacks}, $(B_K)_{\mathfrak{m}}$ is regular.

Because $\mathrm{Spec}(B_K)$ is regular and Noetherian, its irreducible components are pairwise disjoint open and closed connected components. The diagonal $\Delta : \mathrm{Spec}(A_K) \hookrightarrow \mathrm{Spec}(B_K)$ is irreducible and lies in a single connected component $X_0$. For every closed point $y \in \mathrm{Spec}(A_K)$, the integral fibre $p_1^{-1}(y)$ meets $\Delta(\mathrm{Spec}(A_K))$ at $\Delta(y)$ and is contained in $X_0$. Every closed point of $\mathrm{Spec}(B_K)$ maps to a closed point of $\mathrm{Spec}(A_K)$ and therefore lies in $X_0$. Since every nonempty connected component of an affine Noetherian scheme contains a closed point, $\mathrm{Spec}(B_K) = X_0$ is connected and regular, so $B_K$ is an integral domain. As $B$ is $\Lambda$-flat, it embeds into $B_K$ and is also an integral domain.

\textit{Step 4. Formal subgroup structure of $\mathrm{Spf}(A)$.} The comultiplication and antipode of $R$ are
\[
\begin{aligned}
m^\# &: R\longrightarrow R\mathbin{\widehat\otimes}_\Lambda R,
& q_i&\longmapsto q_i\mathbin{\widehat\otimes}q_i,\\
\iota^\# &: R\longrightarrow R,
& q_i&\longmapsto q_i^{-1}.
\end{aligned}
\]
On the completed local rings at the identity, write $l_i=\log(q_i)$, with $l_i'$ denoting the logarithmic coordinates on the second factor. In these coordinates $m^\#(l_i)=l_i+l_i'$ and $\iota^\#(l_i)=-l_i$. By \eqref{eq:P-hat-linear}, the completed ideal $\widehat P$ is generated by the linear forms $z_1,\dots,z_h$. Since $z_k(l+l')=z_k(l)+z_k(l')$ and $z_k(-l)=-z_k(l)$, the comultiplication sends $P$ into the sum of the ideals generated by the $z_k(l)$ and the $z_k(l')$ in $K[[l,l']]$, and the antipode sends $P$ into $\widehat P$.
Let $\mathfrak b\subset B_K$ be the maximal ideal at the pair of identity points. By \eqref{eq:P-hat-linear}, there is an isomorphism
\[
\widehat{(B_K)_{\mathfrak b}} \simeq K[[l_1, \dots, l_n, l_1', \dots, l_n']]/\bigl(z_k(l),z_k(l')\mid 1\le k\le h\bigr).
\]
The quotient on the right is a power series ring in $2d$ variables over $K$. Thus $m^\#(P)$ vanishes in $\widehat{(B_K)_{\mathfrak b}}$ and $\iota^\#(P)$ vanishes in $\widehat{A_{\mathfrak a}}$. Because $B$ and $A$ are integral domains and the localizations $(B_K)_{\mathfrak{b}}$ and $A_{\mathfrak a}$ are Noetherian local domains, the Krull intersection theorem implies that both composites
\[
B \hookrightarrow (B_K)_{\mathfrak{b}} \hookrightarrow \widehat{(B_K)_{\mathfrak b}}, \qquad A = R/P \hookrightarrow A_{\mathfrak a} \hookrightarrow \widehat{A_{\mathfrak a}}
\]
are injective. As $B=(R\mathbin{\widehat\otimes}_{\Lambda}R)/(P\mathbin{\widehat\otimes}_{\Lambda}R+R\mathbin{\widehat\otimes}_{\Lambda}P)$, we obtain $m^\#(P)\subset P\mathbin{\widehat\otimes}_{\Lambda}R+R\mathbin{\widehat\otimes}_{\Lambda}P$ and $\iota^\#(P) \subset P$. Together with $P \subset (x_1, \dots, x_n)$, this shows that $P$ is a Hopf ideal of $R$, so $\mathrm{Spf}(A) \subset \widehat{\mathbf{G}}_m^n$ is a closed formal subgroup.

\textit{Step 5. Recovery of the primitive sublattice via Cartier duality.} By Subsection~\ref{subsec:lattice-dictionary}, there is an isomorphism $R\xrightarrow{\sim}\Lambda[[\mathbf Z_\ell^n]]$, $q^m\mapsto[m]$. For $r\ge 1$, let $R_r=\Lambda/(\varpi^r)$ with the discrete topology and $A_r=A/\varpi^r A$, so $R_1$ is the residue field. For a profinite group $H$, write $R_r[[H]]=\varprojlim_N R_r[H/N]$, where $N$ ranges over the open normal subgroups of $H$, and let $\Lambda[[H]]=\varprojlim_r R_r[[H]]$. Let $C_{R_r}(H)$ be the Hopf algebra of locally constant $R_r$-valued functions on $H$, with pointwise multiplication. Since $R_r$ is discrete and $H$ is compact, each such function factors through a finite quotient of $H$. By duality for finite group algebras, there is an isomorphism of Hopf algebras over $R_r$
\[
\mathrm{Hom}_{R_r,\mathrm{cont}}(R_r[[\mathbf Z_\ell^n]],R_r)
\xrightarrow{\sim} C_{R_r}(\mathbf Z_\ell^n),\qquad \lambda\longmapsto\bigl(g\mapsto\lambda([g])\bigr).
\]
Since $A$ is flat over $\Lambda$, each quotient $A_r$ is flat and pseudocompact over the Artinian local ring $R_r$, hence topologically free and projective in the category of pseudocompact $R_r$-modules \cite[Exp.~VII\textsubscript{B}, Corollary~0.3.8 and Section~2.2.2]{DG70}. By Cartier duality over $R_r$, the continuous dual $D_r := \mathrm{Hom}_{R_r,\mathrm{cont}}(A_r, R_r)$ is a flat commutative Hopf algebra over $R_r$, and the transpose of the quotient map $R_r[[\mathbf Z_\ell^n]]\twoheadrightarrow A_r$ is an injective homomorphism of Hopf algebras $D_r\hookrightarrow C_{R_r}(\mathbf Z_\ell^n)$. Since $A_r$ is projective, the quotient map admits a continuous $R_r$-linear section. Its transpose is a retraction of $D_r\hookrightarrow C_{R_r}(\mathbf Z_\ell^n)$. Identifying $D_r$ with its image inside $C_{R_r}(\mathbf{Z}_\ell^n)$, this $R_r$-linear splitting implies that
\begin{equation}\label{eq:Dr-split}
D_r \cap \varpi^j C_{R_r}(\mathbf{Z}_\ell^n) = \varpi^j D_r \qquad (0 \le j \le r),
\end{equation}
and compatibility of continuous duality with reduction modulo $\varpi$ identifies $D_r/\varpi D_r \simeq D_1 \hookrightarrow C_{R_1}(\mathbf{Z}_\ell^n)$.

On the special fibre $r = 1$, define an equivalence relation $\sim$ on $\mathbf{Z}_\ell^n$ by setting $g \sim h$ if and only if $f(g) = f(h)$ for all $f \in D_1$. Each $f \in D_1 \subset C_{R_1}(\mathbf{Z}_\ell^n)$ takes only finitely many values. Hence polynomial interpolation over $R_1$ shows that the characteristic function of each level set of $f$ is in $D_1$. Since $\mathbf{Z}_\ell^n$ is compact and $D_1$ separates distinct equivalence classes for $\sim$, $D_1$ contains the characteristic function of every clopen subset of $\mathbf{Z}_\ell^n$ saturated with respect to $\sim$, and hence $D_1$ consists of all locally constant $R_1$-valued functions on $\mathbf{Z}_\ell^n/{\sim}$. For $f\in D_1$, write $\Delta(f)=\sum_\nu a_\nu\otimes b_\nu$ with $a_\nu,b_\nu\in D_1$. Then
\[
f(g+h)=\sum_\nu a_\nu(g)b_\nu(h).
\]
Consequently, $g\sim g'$ and $h\sim h'$ imply $g+h\sim g'+h'$. The antipode sends $f$ to the function $g\mapsto f(-g)$, so $g\sim g'$ also implies $-g\sim-g'$. Thus the equivalence class $M := \{g \in \mathbf{Z}_\ell^n \mid g \sim 0\}$ is a closed subgroup of $\mathbf{Z}_\ell^n$, and
\begin{equation}\label{eq:D1-quotient}
D_1 = C_{R_1}(\mathbf{Z}_\ell^n/M) \subset C_{R_1}(\mathbf{Z}_\ell^n).
\end{equation}
Since $M$ is closed and stable under integer multiplication, it is a $\mathbf Z_\ell$-submodule.

For $r\ge1$, let $E_r=C_{R_r}(\mathbf Z_\ell^n/M)$ inside $C_{R_r}(\mathbf Z_\ell^n)$. Given a clopen subset $V\subset\mathbf Z_\ell^n/M$, its characteristic function in $D_1$ lifts uniquely to an idempotent in $D_r$, since $\varpi D_r$ is nilpotent. Uniqueness in the ambient function ring identifies the lift with $\mathbf 1_V$. Such functions span $E_r$, so $E_r\subset D_r$.

Conversely, if $d\in D_r$, choose $e\in E_r$ having the same reduction as $d$ in $D_1=E_1$. By \eqref{eq:Dr-split}, $d-e\in\varpi D_r$. Hence $D_r=E_r+\varpi D_r$. Iterating this equality and using $\varpi^r=0$ gives $D_r=E_r$.

Applying continuous Cartier duality over $R_r$ to $D_r=C_{R_r}(\mathbf Z_\ell^n/M)$ gives $A_r\simeq R_r[[\mathbf Z_\ell^n/M]]$. These isomorphisms commute with reduction from $r+1$ to $r$, so their inverse limit gives an isomorphism of completed group algebras
\[
A \simeq \varprojlim_r A_r \xrightarrow{\sim} \Lambda[[\mathbf{Z}_\ell^n/M]],\qquad q^m\bmod P\longmapsto[m+M].
\]
Since $[m+M]=1$ exactly when $m\in M$, we recover $M$ as $\{m\in\mathbf Z_\ell^n\mid q^m-1\in P\}$. In particular, $P$ determines $M$ uniquely.

Suppose that $\ell m\in M$ and let $u$ be the image of $q^m$ in $A$. Then $(u-1)(1+u+\cdots+u^{\ell-1})=u^\ell-1=0$. The augmentation $A\to\Lambda$ sends the second factor to $\ell\ne0$. Since $A$ is a domain, $u=1$, so $m\in M$. Thus $M$ is primitive. By Lemma~\ref{lem:primitive-subtori}(2), the kernel of $R\to\Lambda[[\mathbf Z_\ell^n/M]]$ is $I_M=(q^m-1\mid m\in M)$, and the displayed isomorphism identifies this kernel with $P$.
\end{proof}

\begin{proposition}[Classification of $\Phi^\#$-stable prime ideals]\label{prop:stable-primes-binomial}
Let $\mathfrak{p} \subset R$ be prime, with $\ell \notin \mathfrak{p}$ and $\Phi^\#(\mathfrak{p}) \subset \mathfrak{p}$. There is a unique primitive sublattice $M \subseteq \mathbf{Z}_\ell^n$ satisfying $U(M) \subseteq M$ such that
\[
\mathfrak{p} = I_M = (q^m - 1 \mid m \in M).
\]
\end{proposition}

\begin{proof}
Let $A=R/\mathfrak p$. Lemma~\ref{lem:pass-identity} gives $\mathfrak p\subset(x_1,\dots,x_n)$, and Lemma~\ref{lem:total-invariance} gives an injective endomorphism $\varphi$ of $\mathrm{Frac}(A)$.

Let $C$ be the relative algebraic closure of $K$ in $\mathrm{Frac}(A)$. Choose a formal Noether normalization $B_0=\Lambda[[y_1,\dots,y_{\dim A-1}]]\hookrightarrow A$ as in the proof of Lemma~\ref{lem:diff-alg-const}(1). The kernel of the continuous differential on $\mathrm{Frac}(B_0)$ is $K$. Thus $K$ is algebraically closed in $\mathrm{Frac}(B_0)$, and in characteristic zero this field is linearly disjoint over $K$ from every finite subextension $C_0/K$ of $C/K$. Hence
\[
[C_0:K]=[\mathrm{Frac}(B_0)C_0:\mathrm{Frac}(B_0)]
\le [\mathrm{Frac}(A):\mathrm{Frac}(B_0)].
\]
The upper bound is independent of $C_0$, so $[C:K]<\infty$. Since $\varphi$ fixes $K$, its restriction to $C$ is a $K$-automorphism. Choose $h\ge1$ such that $\varphi^h$ fixes $C$ pointwise.

Let $\mathcal O_C$ be the valuation ring of $C$. Since $\mathcal O_C$ is finite free over $\Lambda$, we have $S_C=\mathcal O_C[[x_1,\dots,x_n]]=R\otimes_\Lambda\mathcal O_C$. The multiplication map $\mu:S_C\to\mathrm{Frac}(A)$, $r\otimes c\mapsto(r\bmod\mathfrak p)c$, has prime kernel $P$, which contracts to $\mathfrak p$. Its image $B=A\mathcal O_C$ is finite over $A$ and lies in $\mathrm{Frac}(A)$, so $\dim B=\dim A$ and $\mathrm{Frac}(B)=\mathrm{Frac}(A)$.

Since $\varphi^h$ fixes $C$, the diagram
\[
\begin{tikzcd}[column sep=normal]
S_C \arrow[r, "(\Phi^\#)^h\otimes1"] \arrow[d, "\mu"'] & S_C \arrow[d, "\mu"] \\
\mathrm{Frac}(A) \arrow[r, "\varphi^h"'] & \mathrm{Frac}(A)
\end{tikzcd}
\]
commutes. Thus $(\Phi^\#)^h\otimes1$ preserves $P$. Its exponents are $h t_1,\dots,h t_n$, and $\ell\notin P$ since $P\cap R=\mathfrak p$. Lemma~\ref{lem:pass-identity} now gives $P\subset(x_1,\dots,x_n)S_C$.

Apply Lemma~\ref{lem:fixed-elements-algebraic} over $\mathcal O_C$ to $B=S_C/P$ and the endomorphism induced by $(\Phi^\#)^h\otimes1$. Every element of $\mathrm{Frac}(B)$ fixed by $\varphi^h$ is then algebraic over $C$. Since $C$ is algebraically closed in $\mathrm{Frac}(B)=\mathrm{Frac}(A)$ and $\varphi^h$ fixes $C$ pointwise, $\mathrm{Frac}(B)^{\varphi^h}=C$. Order the coordinates as in Lemma~\ref{lem:log-diff-weight-pure}. Since $h t_j=h t_i$ exactly when $t_j=t_i$, Lemma~\ref{lem:log-diff-weight-pure} gives $w_j^{(i)}\in C$ such that, for $\dim A\le i\le n$,
\[
\frac{dq_i}{q_i} = \sum_{\substack{1 \le j \le \dim A - 1 \\ t_j = t_i}} w_j^{(i)} \frac{dq_j}{q_j} \qquad \text{in } \Omega^{\mathrm{cont}}_{B/\mathcal{O}_C} \otimes_B \mathrm{Frac}(A).
\]
Consequently, the kernel of the map $\mathrm{Frac}(A)^n\to\Omega^{\mathrm{cont}}_{B/\mathcal O_C}\otimes_B\mathrm{Frac}(A)$, $(a_i)\mapsto\sum_i a_i\,dq_i/q_i$, has a basis over $\mathrm{Frac}(A)$ consisting of the following vectors in $C^n$:
\[
e_i - \sum_{\substack{1 \le j \le \dim A - 1 \\ t_j = t_i}} w_j^{(i)} e_j, \qquad i = \dim A, \dots, n.
\]
By Lemma~\ref{lem:const-log-tangent}(1) applied to $P \subset S_C$ over $\mathcal{O}_C$, the extended ideal $\mathfrak{a} := P C[[x_1, \dots, x_n]]$ is generated by $C$-linear forms in $\log(q_1), \dots, \log(q_n)$.

Consider the $\mathcal O_C$-algebra endomorphism $\Sigma_C^\#:S_C\to S_C$, $q_i\mapsto q_i^\ell$. In logarithmic coordinates it multiplies every linear form by $\ell$, so it preserves $\mathfrak a$. By Lemma~\ref{lem:ext-contract-identity}(3), $\mathfrak a\cap S_C=P$. Hence $\Sigma_C^\#(P)\subset P$. Now $C$ is algebraically closed in $\mathrm{Frac}(S_C/P)$, so Lemma~\ref{lem:equal-weight-binomial} with $t=1$ gives a primitive sublattice $M\subset\mathbf Z_\ell^n$ with $P=I_M S_C$.

Faithful flatness of $\mathcal O_C$ over $\Lambda$ gives $\mathfrak p=P\cap R=I_M$. We recover $M$ as $\{m\in\mathbf Z_\ell^n\mid q^m-1\in\mathfrak p\}$, which proves uniqueness. For $m\in M$, stability gives $q^{Um}-1=\Phi^\#(q^m-1)\in\mathfrak p$. Thus $Um\in M$.
\end{proof}

By Lemma~\ref{lem:primitive-subtori}(2), the identity $\mathfrak p=I_M$ identifies $\operatorname{Spf}(R/\mathfrak p)$ with the formal subtorus $T_M$. The inclusion $U(M)\subset M$ makes $T_M$ stable under $\Phi$.

\begin{corollary}\label{cor:fixed-field}
In the setting of Lemma~\ref{lem:fixed-elements-algebraic}, the fixed field is $\mathrm{Frac}(A)^\varphi=K$.
\end{corollary}

\begin{proof}
By Proposition~\ref{prop:stable-primes-binomial}, $\mathfrak p=I_M$ for a primitive sublattice $M\subset\mathbf Z_\ell^n$. Extend a basis of $M$ to a basis of $\mathbf Z_\ell^n$. In these coordinates $A$ is a formal power series ring over $\Lambda$, so the kernel of the continuous differential on $\mathrm{Frac}(A)$ is $K$. In particular, $K$ is algebraically closed in $\mathrm{Frac}(A)$. Lemma~\ref{lem:fixed-elements-algebraic} therefore gives $\mathrm{Frac}(A)^\varphi\subset K$, and $\varphi$ fixes $K$.
\end{proof}

\begin{remark}[A gap in the proof of {\cite[Lemma~2.1(2)]{Jon00}}]\label{rem:const-field-step}
In the last paragraph of the proof of \cite[Lemma~2.1(2)]{Jon00}, de Jong drops the stability hypothesis on $\mathfrak p$ and claims that any $h\in\mathrm{Frac}(A)$ with $dh=0$ is in $K$. Differentiating the minimal polynomial as in the proof of Lemma~\ref{lem:diff-alg-const}(2) shows that $h$ is algebraic over $K$. To conclude that the polynomial is linear, de Jong evaluates $h$ at the identity. The augmentation $A\to\Lambda$ induces a homomorphism $A_{\mathfrak a}\to K$, $f/g\mapsto f(0)/g(0)$, where $\mathfrak a=(x_1,\dots,x_n)A$. The proof does not establish that $h\in A_{\mathfrak a}$.

Without the stability hypothesis, the implication $dh=0\Rightarrow h\in K$ fails for
\[
A_0=\Lambda[[x,y]]/(y^2-\varpi x^2),\qquad h=y/x\in\mathrm{Frac}(A_0).
\]
For $\alpha^2=\varpi$, the substitution $y=\alpha x$ embeds $A_0$ in $\mathcal O_{K(\alpha)}[[x]]$. Thus $A_0$ is a flat complete local domain with augmentation $x,y\mapsto0$. Since $h^2=\varpi$, we have $2h\,dh=0$ and hence $dh=0$. If $h$ were in $(A_0)_{(x,y)}$, its value at the identity would be a square root of $\varpi$ in $K$, which is impossible since $v_K(\varpi)=1$.

Under the stability hypothesis, Proposition~\ref{prop:stable-primes-binomial} identifies $A$ with a formal power series ring over $\Lambda$ in suitable coordinates. Corollary~\ref{cor:fixed-field} then gives $\mathrm{Frac}(A)^\varphi=K$. When the power map acts on $\Lambda$ by an automorphism $\tau$, Corollary~\ref{cor:semilinear-fixed-field} gives the fixed field $K^\tau$.
\end{remark}

\subsection{Formal subschemes preserved by an iterate}\label{subsec:binomial-rigidity}
If a prime ideal $\mathfrak p\subset R$ satisfies $\ell\notin\mathfrak p$ and $(\Phi^\#)^c(\mathfrak p)\subset\mathfrak p$ for some $c\ge1$, Proposition~\ref{prop:stable-primes-binomial} applied to $\Phi^c$ gives $\mathfrak p=I_M$ with $U^c(M)\subset M$.

\begin{lemma}\label{lem:uc-implies-u}
Let $U = \mathrm{diag}(\ell^{t_1}, \dots, \ell^{t_n})$ with integers $t_i \ge 1$, and let $M \subset \mathbf{Z}_\ell^n$ be a primitive sublattice. If $U^c(M) \subset M$ for some integer $c \ge 1$, then $U(M) \subset M$.
\end{lemma}

\begin{proof}
Let $W := M \otimes_{\mathbf{Z}_\ell} \mathbf{Q}_\ell \subset \mathbf{Q}_\ell^n$. The inclusion $U^c(M) \subset M$ implies $U^c(W) \subset W$, and since $U^c$ is an invertible linear operator on the finite dimensional vector space $\mathbf{Q}_\ell^n$, we have $U^c(W) = W$. For each distinct weight $t\in\{t_1,\dots,t_n\}$, the spectral projection onto $\ker(U^c-\ell^{ct})$ is a polynomial in $U^c$ over $\mathbf Q_\ell$. These projections preserve $W$, so
\[
W = \bigoplus_{t\in\{t_1,\dots,t_n\}} \bigl(W\cap\ker(U^c-\ell^{ct})\bigr).
\]
On $\ker(U^c-\ell^{ct})$, the diagonal matrix $U$ acts by $\ell^t$. Thus $U(W)=W$. Primitivity gives $M=W\cap\mathbf Z_\ell^n$, so $U(M) \subset W \cap \mathbf{Z}_\ell^n = M$.
\end{proof}

\begin{remark}
The conclusion of Lemma~\ref{lem:uc-implies-u} can fail when $U$ is not diagonal. Take $U = \ell \bigl(\begin{smallmatrix} 0 & 1 \\ 1 & 0 \end{smallmatrix}\bigr)$ and $M = \mathbf{Z}_\ell e_1$. Then $U^2(M) = \ell^2 M \subset M$ and $U(M) = \ell \mathbf{Z}_\ell e_2 \not\subset M$.
\end{remark}

\begin{lemma}\label{lem:irred-comp-finite-map}
Let $f:\mathfrak X\to\mathfrak X$ be a finite morphism, let $\mathfrak Z\subset\mathfrak X$ be a reduced equidimensional closed formal subscheme, and let $\mathfrak Y=f(\mathfrak Z)$ be its scheme theoretic image. For each irreducible component $\mathfrak Z_\alpha$ of $\mathfrak Z$ the following hold:
\begin{enumerate}
\item[(1)] The scheme theoretic image $f(\mathfrak Z_\alpha)$ is an irreducible component of $\mathfrak Y$.
\item[(2)] The subscheme $\mathfrak Z_\alpha$ is an irreducible component of $(f^{-1}(\mathfrak Y))_{\mathrm{red}}$.
\end{enumerate}
\end{lemma}

\begin{proof}
(1) Let $d$ be the common dimension of the components of $\mathfrak Z$. Each $\mathfrak Z_\alpha$ maps finitely and dominantly onto its scheme theoretic image $\mathfrak Y_\alpha=f(\mathfrak Z_\alpha)$. Thus $\mathfrak Y_\alpha$ is irreducible of dimension $d$. Contraction of ideals commutes with intersections, so the scheme theoretic union of the $\mathfrak Y_\alpha$ is $\mathfrak Y$. Since they have the same dimension, the distinct images are exactly its irreducible components.

(2) The reduced subscheme $\mathfrak Z$ lies in $(f^{-1}(\mathfrak Y))_{\mathrm{red}}$. Choose an irreducible component $\mathfrak W$ of this inverse image containing $\mathfrak Z_\alpha$. Its image in $\mathfrak Y$ is irreducible and contains $\mathfrak Y_\alpha$, so (1) gives $f(\mathfrak W)=\mathfrak Y_\alpha$. The finite dominant morphism $\mathfrak W\to\mathfrak Y_\alpha$ gives an integral inclusion of coordinate rings. Both the zero ideal of $\mathcal O(\mathfrak W)$ and the prime defining $\mathfrak Z_\alpha$ contract to zero in $\mathcal O(\mathfrak Y_\alpha)$. By incomparability they coincide. Hence $\mathfrak Z_\alpha=\mathfrak W$.
\end{proof}

\begin{lemma}\label{lem:base-change-components}
Let $D$ be a flat, reduced, equidimensional Noetherian $\Lambda$-algebra. For a finite extension $L/K$ with valuation ring $\mathcal O_L$, let $D'=D\otimes_\Lambda\mathcal O_L$.
\begin{enumerate}
\item[(1)] The algebra $D'$ is flat over $\mathcal O_L$, reduced, and equidimensional.
\item[(2)] For each minimal prime $\mathfrak p$ of $D$, the ideal $\mathfrak pD'$ is radical and has minimal prime decomposition
\[
\mathfrak pD'=\bigcap_{\substack{\mathfrak q'\in\mathrm{Min}(D')\\ \mathfrak q'\cap D=\mathfrak p}}\mathfrak q'.
\]
\end{enumerate}
\end{lemma}

\begin{proof}
(1) Since $\mathcal O_L$ is finite free over $\Lambda$, $D'$ is finite faithfully flat over $D$. The flatness of $D$ over $\Lambda$ gives flatness of $D'$ over $\mathcal O_L$. In particular, $D'$ embeds into $D'[1/\ell]=D[1/\ell]\otimes_K L$. The localization $D'[1/\ell]$ is reduced because $D[1/\ell]$ is reduced and $L/K$ is finite separable.

By going down and faithful flatness, the minimal primes of $D'$ contract to the minimal primes of $D$, and every minimal prime of $D$ occurs. If $\mathfrak q'\cap D=\mathfrak p$ is such a pair, the inclusion $D/\mathfrak p\hookrightarrow D'/\mathfrak q'$ is finite and integral, so the two rings have the same dimension. Equidimensionality of $D$ therefore gives equidimensionality of $D'$.

(2) Let $\mathfrak p$ be a minimal prime of $D$. Since $D$ is flat over $\Lambda$, $\ell\notin\mathfrak p$ and $D/\mathfrak p$ is flat over $\Lambda$. The quotient $D'/\mathfrak pD'$ embeds into the reduced ring $(D/\mathfrak p)[1/\ell]\otimes_K L$, so $D'/\mathfrak pD'$ is reduced. By integrality and going down, the minimal primes over $\mathfrak pD'$ are exactly the minimal primes of $D'$ contracting to $\mathfrak p$. Their intersection is therefore $\mathfrak pD'$.
\end{proof}

\begin{theorem}[Quasilinearity from invariance under an iterate]\label{thm:binomial-rigidity}
Let $I\subset R$ be a radical ideal such that $R/I$ is flat over $\Lambda$ and equidimensional. Assume that $(\Phi^\#)^c(I)\subset I$ for some $c\ge1$, and let $\mathfrak Y=V(I)$.
\begin{enumerate}
\item[(1)] The subscheme $\mathfrak Y$ is $\Phi$-stable quasilinear. More precisely, over a finite cyclotomic extension $L=K(\mu_{\ell^s})$ with $s\ge0$ it has a scheme theoretic decomposition
\[
\mathfrak Y_{\mathcal O_L}
=\overset{\mathrm{sch}}{\bigcup_{\alpha=1}^N}
 t_{\zeta_\alpha}\bigl((T_{M_\alpha})_{\mathcal O_L}\bigr),
\]
where $\zeta_\alpha\in\mathfrak X(\mathcal O_L)$ are torsion points and $M_\alpha\subset\mathbf Z_\ell^n$ are primitive $U$-stable sublattices. The union is empty when $N=0$.
\item[(2)] If $\mathfrak Y$ is irreducible, then $\mathfrak Y=T_M$ for a primitive $U$-stable sublattice $M$.
\end{enumerate}
\end{theorem}

\begin{proof}
(1) If $I=R$, take the empty union. Otherwise let $\mathfrak Y_1,\dots,\mathfrak Y_r$ be the irreducible components of $\mathfrak Y$, with defining primes $\mathfrak p_1,\dots,\mathfrak p_r$. The hypothesis $(\Phi^\#)^c(I)\subset I$ gives $\Phi^c(\mathfrak Y)\subset\mathfrak Y$. Since $\Phi^c$ is finite, the image of each component has the same dimension as that component. Equidimensionality therefore makes each image a component of $\mathfrak Y$. The forward orbit of each component under $\Phi^c$ eventually reaches a periodic component.

Consider first a periodic component $\mathfrak Y_\gamma$ of period $k$. Since $\Phi^{kc}(\mathfrak Y_\gamma)=\mathfrak Y_\gamma$, \eqref{eq:image-inverse-ideals} gives $(\Phi^\#)^{kc}(\mathfrak p_\gamma)\subseteq\mathfrak p_\gamma$. Flatness of $R/I$ over $\Lambda$ gives $\ell\notin\mathfrak p_\gamma$. Proposition~\ref{prop:stable-primes-binomial} applied to this iterate gives $\mathfrak p_\gamma=I_{M_\gamma}$ for a primitive sublattice satisfying $U^{kc}(M_\gamma)\subset M_\gamma$. Lemma~\ref{lem:uc-implies-u} gives $U(M_\gamma)\subset M_\gamma$. Thus every periodic component is a $\Phi$-stable formal subtorus.

For any component $\mathfrak Y_\alpha$, choose $k\ge1$ such that $\Phi^{kc}(\mathfrak Y_\alpha)=T_M$ is periodic. Apply Lemma~\ref{lem:irred-comp-finite-map}(2) to $f=\Phi^{kc}$ and the irreducible subscheme $\mathfrak Y_\alpha$, whose image is $T_M$. It follows that $\mathfrak Y_\alpha$ is an irreducible component of $((\Phi^{kc})^{-1}(T_M))_{\mathrm{red}}$. By Proposition~\ref{prop:images-subtori}(2), $(\Phi^{kc})^{-1}(T_M)=V(I_{M'_\alpha})$ with $M'_\alpha=U^{kc}(M)$. Since $U(M)\subset M$, both $M'_\alpha$ and its saturation $(M'_\alpha)^{\mathrm{sat}}$ are $U$-stable, with $U((M'_\alpha)^{\mathrm{sat}})\subset(M'_\alpha)^{\mathrm{sat}}$ by Remark~\ref{rem:sat-u-stable}.

Choose $s\ge0$ at least as large as every exponent $a_i$ occurring in \eqref{eq:smith-L} for the finitely many $M'_\alpha$, and let $L=K(\mu_{\ell^s})$. Lemma~\ref{lem:nonprimitive-torsion-union} identifies the components of $V(I_{M'_\alpha})_{\mathcal O_L}$ as torsion translates of $T_{(M'_\alpha)^{\mathrm{sat}}}$. The prime defining $\mathfrak Y_\alpha$ in this inverse image is $\mathfrak p_\alpha/I_{M'_\alpha}$. The ring $R/I_{M'_\alpha}$ is flat over $\Lambda$, reduced, and equidimensional by Lemma~\ref{lem:binomial-model-properties}(2). Lemma~\ref{lem:base-change-components}(2) therefore identifies $(\mathfrak Y_\alpha)_{\mathcal O_L}$ with the scheme theoretic union of the translates lying over that prime. Finally, flat base change preserves the finite intersection $I=\bigcap_\alpha\mathfrak p_\alpha$. Taking the union over $\alpha$ gives the stated decomposition of $\mathfrak Y_{\mathcal O_L}$.

(2) If $\mathfrak Y$ is irreducible, then the closed image $\Phi^c(\mathfrak Y)\subset\mathfrak Y$ has the same dimension as $\mathfrak Y$, since $\Phi^c$ is finite. Hence $\Phi^c(\mathfrak Y)=\mathfrak Y$. The classification of periodic components in (1) gives $\mathfrak Y=T_M$.
\end{proof}

\subsection{The formal dynamical criterion}\label{subsec:dynamical-criterion}
\begin{proposition}[$\Phi$-stable formal quasilinearity criterion]\label{prop:formal-dynamical-criterion}
Let $\Phi$ be the diagonal power endomorphism in \eqref{eq:diag-phi}, and let $\mathfrak Z\subset\mathfrak X$ be a closed formal subscheme that is flat over $\Lambda$, reduced, and equidimensional.
\begin{enumerate}
\item[(1)] The subscheme $\mathfrak Z$ is $\Phi$-stable quasilinear if and only if $\Phi^a(\mathfrak Z)\subset\Phi^b(\mathfrak Z)$ as closed formal subschemes for some integers $a>b\ge0$.
\item[(2)] Moreover, $\mathfrak Z$ is a $\Phi$-stable formal subtorus if and only if it is irreducible and $\Phi^c(\mathfrak Z)\subset\mathfrak Z$ for some integer $c\ge1$.
\end{enumerate}
\end{proposition}

\begin{lemma}\label{lem:image-preserves-flat-red-equidim}
Let $\mathfrak{Z} \subset \mathfrak{X}$ be a closed formal subscheme that is flat over $\Lambda$, reduced, and equidimensional. For every integer $b \ge 0$, the scheme theoretic image $\mathfrak{Y} := \Phi^b(\mathfrak{Z})$ is also flat over $\Lambda$, reduced, and equidimensional.
\end{lemma}

\begin{proof}
By \eqref{eq:image-inverse-ideals}, $(\Phi^\#)^b$ induces an injection $R/I_{\mathfrak Y}\hookrightarrow R/I_{\mathfrak Z}$ of $\Lambda$-algebras. Hence $R/I_{\mathfrak Y}$ is reduced and $\ell$-torsion free, and therefore flat over $\Lambda$.

By Remark~\ref{rem:finite-locally-free}, $\Phi^b$ is finite. Each irreducible component of its image is the image of a component of $\mathfrak Z$, and the finite dominant map between these components preserves dimension. Their common dimension gives equidimensionality of $\mathfrak Y$.
\end{proof}

\begin{lemma}\label{lem:pullback-torsion-translate}
Let $L/K$ be a finite extension, let $T\subset\mathfrak X_{\mathcal O_L}$ be a split formal subtorus, and let $\zeta\in\mathfrak X(\mathcal O_L)$ be a torsion point.
\begin{enumerate}
\item[(1)] Over a finite extension $L'/L$ there is a torsion point $\xi\in\mathfrak X(\mathcal O_{L'})$ such that $\Phi(\xi)=\zeta$.
\item[(2)] For every such $\xi$, the inverse image of the translate is given by
\[
\Phi^{-1}(t_\zeta(T))_{\mathcal O_{L'}}=t_\xi\bigl(\Phi^{-1}(T)_{\mathcal O_{L'}}\bigr).
\]
\end{enumerate}
\end{lemma}

\begin{proof}
For $\eta\in\mathfrak X(\mathcal O_L)$, \eqref{eq:diag-phi} gives
\begin{equation}\label{eq:point-action}
q_i(\Phi(\eta))=q_i(\eta)^{\ell^{t_i}}\qquad(1\le i\le n).
\end{equation}
For each $i$, choose an $\ell^{t_i}$-th root of $q_i(\zeta)$ in an algebraic closure of $L$. These roots have $\ell$-power order and therefore lie in $1+\mathfrak m_{L'}$ for a common finite cyclotomic extension $L'/L$. They define a torsion point $\xi$ with $\Phi(\xi)=\zeta$. The identity $\Phi\circ t_\xi=t_\zeta\circ\Phi$ gives $\Phi^{-1}(t_\zeta(T))_{\mathcal O_{L'}}=t_\xi(\Phi^{-1}(T)_{\mathcal O_{L'}})$.
\end{proof}

\begin{lemma}\label{lem:inverse-preserves-quasilinear}
If $\mathfrak{Z} \subset \mathfrak{X}$ is $\Phi$-stable quasilinear, then $\Phi^{-1}(\mathfrak{Z})$ is $\Phi$-stable quasilinear.
\end{lemma}

\begin{proof}
By Definition~\ref{def:formal-quasilinear} and Lemma~\ref{lem:primitive-subtori}, choose a finite extension $L/K$ and write
\[
\mathfrak Z_{\mathcal O_L}=\bigcup_{i=1}^r t_{\zeta_i}\bigl((T_{M_i})_{\mathcal O_L}\bigr)
\]
as a scheme theoretic union, where $\zeta_i\in\mathfrak X(\mathcal O_L)$ are torsion points and $M_i\subset\mathbf Z_\ell^n$ are primitive sublattices satisfying $U(M_i)\subset M_i$ by Lemma~\ref{lem:subtorus-stability}. Since $\Phi_{\mathcal O_L}$ is finite locally free by Remark~\ref{rem:finite-locally-free}, extending ideals along its pullback preserves finite intersections. It suffices to treat one translate.

By Lemma~\ref{lem:pullback-torsion-translate}, over a finite cyclotomic extension of $L$ the inverse image of $t_{\zeta_i}((T_{M_i})_{\mathcal O_L})$ is a torsion translate of $\Phi^{-1}(T_{M_i})$. Proposition~\ref{prop:images-subtori}(2) identifies $\Phi^{-1}(T_{M_i})$ with $V(I_{U(M_i)})$. By Lemma~\ref{lem:nonprimitive-torsion-union}, after a further finite cyclotomic extension, $V(I_{U(M_i)})$ is a scheme theoretic union of torsion translates of $T_{(U(M_i))^{\mathrm{sat}}}$. The subtorus $T_{(U(M_i))^{\mathrm{sat}}}$ is $\Phi$-stable by Remark~\ref{rem:sat-u-stable}. A common finite extension works for all $i$.
\end{proof}

\begin{proof}[Proof of Proposition~\ref{prop:formal-dynamical-criterion}]
If $\mathfrak Z=\varnothing$, the decomposition is the empty union and $\Phi(\mathfrak Z)=\mathfrak Z$. Assume $\mathfrak Z\ne\varnothing$.

For the forward implication of (1), choose a finite extension $L/K$ over which $\mathfrak Z$ is a finite scheme theoretic union of translates $t_{\zeta_j}(T_j)$ with $T_j$ a $\Phi$-stable formal subtorus. By \eqref{eq:point-action}, some $N_0\ge1$ satisfies $\Phi^{N_0}(\zeta_j)=e$ for all $j$, since the $\zeta_j$ have $\ell$-power order and all $t_i$ are positive. Hence $\Phi^{N_0}(t_{\zeta_j}(T_j))=\Phi^{N_0}(T_j)$. Scheme theoretic images preserve finite scheme theoretic unions because contraction of ideals commutes with intersections. Applying this to $\Phi^{N_0}$ and $\Phi^{N_0+1}$ gives
\[
\Phi^{N_0+1}(\mathfrak Z_{\mathcal O_L})
=\bigcup_j^{\mathrm{sch}}\Phi^{N_0+1}(T_j)
\subset\bigcup_j^{\mathrm{sch}}\Phi^{N_0}(T_j)
=\Phi^{N_0}(\mathfrak Z_{\mathcal O_L}),
\]
where the inclusion uses $\Phi(T_j)\subset T_j$. Scheme theoretic images under $\Phi^{N_0}$ and $\Phi^{N_0+1}$ commute with the flat base change $\Lambda\to\mathcal O_L$. Faithful flatness therefore gives $\Phi^{N_0+1}(\mathfrak Z)\subset\Phi^{N_0}(\mathfrak Z)$.

Conversely, suppose $\Phi^a(\mathfrak Z)\subset\Phi^b(\mathfrak Z)$ for some $a>b\ge0$. Let $\mathfrak Y=\Phi^b(\mathfrak Z)$. Lemma~\ref{lem:image-preserves-flat-red-equidim} gives flatness, reducedness, and equidimensionality of $\mathfrak Y$. Since $\Phi^{a-b}(\mathfrak Y)\subset\mathfrak Y$, Theorem~\ref{thm:binomial-rigidity}(1) gives $\Phi$-stable quasilinearity of $\mathfrak Y$.

Let $\mathfrak H=(\Phi^b)^{-1}(\mathfrak Y)$. Lemma~\ref{lem:inverse-preserves-quasilinear} applied $b$ times gives $\Phi$-stable quasilinearity of $\mathfrak H$. $\mathfrak H$ is flat over $\Lambda$ and equidimensional since $\Phi^b$ is finite locally free by Remark~\ref{rem:finite-locally-free}. After inverting $\ell$, the Jacobian of $\Phi^b$ in the coordinates $q_i$ is diagonal with unit entries $\ell^{b t_i}q_i^{\ell^{b t_i}-1}$. It follows that $\mathcal O(\mathfrak H)[1/\ell]$ is finite étale over the reduced ring $\mathcal O(\mathfrak Y)[1/\ell]$, hence reduced. Flatness over $\Lambda$ gives $\mathcal O(\mathfrak H)\hookrightarrow\mathcal O(\mathfrak H)[1/\ell]$, so $\mathfrak H$ is reduced.

By Lemma~\ref{lem:irred-comp-finite-map}(2), each component of $\mathfrak Z$ is a component of $\mathfrak H$. Let $D=\mathcal O(\mathfrak H)$ and let $S\subset\mathrm{Min}(D)$ consist of the primes defining the components of $\mathfrak Z$. Since $\mathfrak Z$ is reduced, $I_{\mathfrak Z}/I_{\mathfrak H}=\bigcap_{\mathfrak p\in S}\mathfrak p$. Choose a finite extension $L/K$ splitting $\mathfrak H$ into torsion translates of $\Phi$-stable formal subtori, and let $D'=D\otimes_\Lambda\mathcal O_L$. The translates maximal under inclusion are the irreducible components of $\mathfrak H_{\mathcal O_L}$. Flatness and Lemma~\ref{lem:base-change-components}(2) give
\[
\left(\bigcap_{\mathfrak p\in S}\mathfrak p\right)D'
=\bigcap_{\mathfrak p\in S}\mathfrak pD'
=\bigcap_{\substack{\mathfrak q'\in\mathrm{Min}(D')\\ \mathfrak q'\cap D\in S}}\mathfrak q'.
\]
Thus $\mathfrak Z_{\mathcal O_L}$ is the scheme theoretic union of the torsion translates indexed by the minimal primes $\mathfrak q'$ with $\mathfrak q'\cap D\in S$. This proves (1).

For (2), a $\Phi$-stable formal subtorus is irreducible and satisfies the inclusion with $c=1$. The converse follows from Theorem~\ref{thm:binomial-rigidity}(2).
\end{proof}

\begin{remark}\label{rem:de-jong-comparison}
When $t_1=\dots=t_n=t$, every primitive sublattice $M\subset\mathbf Z_\ell^n$ satisfies $U(M)=\ell^t M\subset M$. Hence $\Phi$-stable quasilinearity coincides with quasilinearity. Proposition~\ref{prop:formal-dynamical-criterion} then gives the case of \cite[Proposition~1.2]{Jon00} in which $\Phi$ fixes $\Lambda$. Corollary~\ref{cor:semilinear-formal} below recovers \cite[Proposition~1.2]{Jon00} for an arbitrary automorphism $\tau$ of $\Lambda$.
\end{remark}

\subsection{Diagonal power maps with a coefficient automorphism}\label{subsec:semilinear-formal}

Let $\tau$ be an automorphism of $\Lambda$. It preserves the valuation on $K$ and fixes $\mathbf Z_\ell$ pointwise. Let $\Phi_\tau:\mathfrak X\to\mathfrak X$ be the $\tau$-semilinear endomorphism whose pullback $\Phi_\tau^\#:R\to R$ satisfies
\[
\Phi_\tau^\#|_\Lambda=\tau,
\qquad
\Phi_\tau^\#(q_i)=q_i^{\ell^{t_i}}
\qquad(1\le i\le n).
\]
For every integer $j\ge0$,
\[
(\Phi_\tau^\#)^j|_\Lambda=\tau^j,
\qquad
(\Phi_\tau^\#)^j(q_i)=q_i^{\ell^{jt_i}}
\qquad(1\le i\le n).
\]

\begin{proposition}\label{prop:semilinear-primes}
A prime ideal $\mathfrak p\subset R$ with $\ell\notin\mathfrak p$ satisfies $\Phi_\tau^\#(\mathfrak p)\subset\mathfrak p$ if and only if $\mathfrak p=I_M$ for a unique primitive sublattice $M\subset\mathbf Z_\ell^n$ with $U(M)\subset M$.
\end{proposition}

\begin{proof}
Suppose $\Phi_\tau^\#(\mathfrak p)\subset\mathfrak p$. If $f=c+g\in\mathfrak p$ with $0\ne c\in\Lambda$ and $g\in(x_1,\ldots,x_n)$, choose $r\ge0$ such that $\ell^r/c\in\Lambda$. Multiply $f$ by $\ell^r/c$ to obtain an element of $\mathfrak p$ with constant term $\ell^r$. Since $\tau(\ell)=\ell$ and $\tau$ preserves $(\ell,x_1,\ldots,x_n)$, the estimate in Lemma~\ref{lem:m-adic-contraction} shows that the iterates of $(\ell^r/c)f$ under $\Phi_\tau^\#$ converge $\mathfrak m$-adically to $\ell^r$. Closedness of $\mathfrak p$ would give $\ell^r\in\mathfrak p$. Thus $\mathfrak p\subset(x_1,\ldots,x_n)$.

In the completed local ring at the identity, let $y_i=\log(q_i)$. In these coordinates $\Phi_\tau^\#(a y^\nu)=\tau(a)\ell^d y^\nu$, where $d=\sum_i t_i\nu_i$.

The extension $\mathfrak p K[[y_1,\dots,y_n]]$ is preserved by $\Phi_\tau^\#$. Fix $N\ge1$ and order the monomial basis of $K[[y_1,\dots,y_n]]/(y_1,\dots,y_n)^N$. The induced endomorphism carries the image of the extended ideal onto itself, since it is injective and $\tau$-semilinear on this finite dimensional $K$-space. Apply $\Phi_\tau^\#$ to the rows of the reduced row echelon basis of this image, then divide each row by its new pivot. The pivot columns are unchanged, so uniqueness of the reduced row echelon basis recovers the original rows. Thus an entry $a$ in a column of weight $d'$ in a row with pivot of weight $d$ satisfies
\[
a=\ell^{d'-d}\tau(a).
\]
Since $\tau$ preserves valuations, $a\ne0$ forces $d'=d$. Each row is therefore homogeneous for the weights, and $\Phi^\#$ acts on it by the scalar $\ell^d$. Thus $\Phi^\#$ preserves the image of $\mathfrak p K[[y_1,\dots,y_n]]$ modulo $(y_1,\dots,y_n)^N$ for every $N$. Since the extended ideal is closed, $\Phi^\#$ preserves it. Contraction by Lemma~\ref{lem:ext-contract-identity}(3) gives $\Phi^\#(\mathfrak p)\subset\mathfrak p$. Proposition~\ref{prop:stable-primes-binomial} now gives the unique primitive $M$ with $\mathfrak p=I_M$ and $U(M)\subset M$.

Conversely, $\Phi_\tau^\#(q^m-1)=q^{Um}-1$ gives stability for each such $I_M$.
\end{proof}

\begin{corollary}\label{cor:semilinear-fixed-field}
Let $\mathfrak p\subset R$ be a $\Phi_\tau^\#$-stable prime ideal with $\ell\notin\mathfrak p$. The fixed field of the induced endomorphism of $\mathrm{Frac}(R/\mathfrak p)$ is $K^\tau$.
\end{corollary}

\begin{proof}
By Proposition~\ref{prop:semilinear-primes}, $\mathfrak p=I_M$ for a primitive $U$-stable sublattice $M$. Extend a basis of $M$ to a basis of $\mathbf Z_\ell^n$. In these coordinates $R/I_M$ is a formal power series ring over $\Lambda$. Let $\varphi:\mathrm{Frac}(R/I_M)\to\mathrm{Frac}(R/I_M)$ be the field endomorphism induced by $\Phi_\tau^\#$. The endomorphism $\varphi$ preserves the valuation at the height one prime $(\varpi)$. Indeed, identify $(R/I_M)/\varpi$ with $(\Lambda/(\varpi))[[\mathbf Z_\ell^n/M]]$. On this quotient $\Phi_\tau^\#([m+M])=[Um+M]$ for $m\in\mathbf Z_\ell^n$. The induced endomorphism is injective because $U$ induces an injection with finite cokernel on $\mathbf Z_\ell^n/M$ and $\tau$ induces an automorphism of $\Lambda/(\varpi)$. Also, $\tau(\varpi)$ is a unit times $\varpi$.

Return to the original coordinates in which $U$ is diagonal and take $h$ fixed by $\varphi$. Since $\Phi^\#(I_M)\subset I_M$, Lemma~\ref{lem:diff-alg-const}(1) gives a basis of logarithmic differentials among the $dq_i/q_i$. In the expansion of $dh$, the coefficient $a_i$ of $dq_i/q_i$ satisfies $a_i=\ell^{t_i}\varphi(a_i)$. Preservation of the valuation and $t_i>0$ force $a_i=0$. Thus $dh=0$. The computation of $\ker(d)$ in the proof of Lemma~\ref{lem:diff-alg-const}(2), applied to the formal power series ring $R/I_M$, gives $\ker(d)=K$. It follows that $h\in K^\tau$, and every element of $K^\tau$ is fixed.
\end{proof}

\begin{corollary}\label{cor:semilinear-formal}
Let $Z\subset\mathfrak X$ be a closed formal subscheme that is flat over $\Lambda$, reduced, and equidimensional.
\begin{enumerate}
\item[(1)] The inclusion $\Phi_\tau^a(Z)\subset\Phi_\tau^b(Z)$ holds for some integers $a>b\ge0$ if and only if there is a finite extension $L/K$, with valuation ring $\mathcal O_L$, such that $Z_{\mathcal O_L}$ is a finite scheme theoretic union of torsion translates of the base changes of $\Phi_\tau$-stable formal subtori of $\mathfrak X$.
\item[(2)] Moreover, $Z$ is a $\Phi_\tau$-stable formal subtorus if and only if it is irreducible and $\Phi_\tau^c(Z)\subset Z$ for some integer $c\ge1$.
\end{enumerate}
\end{corollary}

\begin{proof}
(1) Suppose $\Phi_\tau^a(Z)\subset\Phi_\tau^b(Z)$ with $a>b\ge0$. The map $\Phi_\tau$ is finite locally free, being the composite of $\Phi$ with the automorphism acting by $\tau$ on coefficients. Finite images preserve the dimensions of irreducible components, so $\Phi_\tau^{a-b}$ maps each component of $\Phi_\tau^b(Z)$ onto a component. Thus for each component of $Z$ some iterate $\Phi_\tau^j$ maps it onto a periodic component $Y$.

Choose a positive multiple $s$ of $a-b$ with $\Phi_\tau^s(Y)=Y$. The coefficient automorphism and exponent matrix of $\Phi_\tau^s$ are $\tau^s$ and $U^s$, respectively. By Proposition~\ref{prop:semilinear-primes}, $Y=T_M$ for a primitive sublattice $M$ with $U^s(M)\subset M$. Lemma~\ref{lem:uc-implies-u} gives $U(M)\subset M$. By Lemma~\ref{lem:irred-comp-finite-map}(2), the original component of $Z$ is a component of $(\Phi_\tau^j)^{-1}(T_M)=V(I_{U^j(M)})$. Lemmas~\ref{lem:nonprimitive-torsion-union} and~\ref{lem:base-change-components} split all these components over a common finite cyclotomic extension into torsion translates of the corresponding $T_{(U^j(M))^{\mathrm{sat}}}$. Each saturated sublattice is $U$-stable by Remark~\ref{rem:sat-u-stable}. Reducedness of $Z$ gives the required scheme theoretic union.

Conversely, choose a finite extension $L/K$ over which $Z_{\mathcal O_L}$ has the stated decomposition. Over $\mathcal O_L$, \eqref{eq:point-action} shows that $\Phi^j$ sends every translating point to the identity for all sufficiently large $j$. For each subtorus in the decomposition we have $\Phi(T_M)=T_M$ by Proposition~\ref{prop:images-subtori}(1), since $U$ preserves $M\otimes_{\mathbf Z_\ell}\mathbf Q_\ell$ bijectively and $M$ is primitive. Thus $\Phi^j(Z)$ is the union of the underlying subtori over the splitting field, and faithful flatness descends this equality to $\Lambda$.

Extend the coefficient automorphism $\tau$ to $R$ by fixing $x_1,\ldots,x_n$. It commutes with $\Phi^\#$, so \eqref{eq:image-inverse-ideals} gives $I_{\Phi_\tau^j(Z)}=\tau^{-j}(I_{\Phi^j(Z)})$. For all sufficiently large $j$, $I_{\Phi^j(Z)}$ is an intersection of ideals $I_M$ generated by binomials over $\mathbf Z_\ell$. Thus $\tau$ preserves the intersection, and $\Phi_\tau^j(Z)=\Phi^j(Z)$ for every such $j$. In particular, $\Phi_\tau^{j+1}(Z)=\Phi_\tau^j(Z)$.

For (2), suppose $Z$ is irreducible and $\Phi_\tau^c(Z)\subset Z$ for $c\ge1$. Applying Proposition~\ref{prop:semilinear-primes} to $\Phi_\tau^c$ and then Lemma~\ref{lem:uc-implies-u} gives $Z=T_M$ with $U(M)\subset M$. Every such $T_M$ satisfies $\Phi_\tau(T_M)\subset T_M$.
\end{proof}

\section{Mixed weight rigid linearity}\label{sec:rigid-spaces}

Let $K$ be the fraction field of $\Lambda$, with $|\ell|=\ell^{-1}$. We use Berthelot's rigid generic fibre functor $\mathfrak Y\mapsto\mathfrak Y_\eta$ from locally Noetherian adic formal $\Lambda$-schemes whose reduced schemes are locally of finite type over the residue field of $\Lambda$ to rigid analytic spaces over $K$; see \cite[Definition~7.0.1, Sections~7.1.1--7.1.3 and~7.2.3]{Jon95}.

\subsection{Affinoid charts on character tori}\label{subsec:affinoid-charts}

\begin{definition}\label{def:rigid-character-space}
For a finite free $\mathbf Z_\ell$-module $\pi$, the \emph{rigid character torus} $X_\pi=(\mathfrak X_\pi)_\eta$ is the rigid analytic generic fibre of $\mathfrak X_\pi$.
\end{definition}

Choose a $\mathbf Z_\ell$-basis $e_1,\ldots,e_n$ of $\pi$ and write $x_i=[e_i]-1$ as in Subsection~\ref{subsec:lattice-dictionary}.

\begin{proposition}\label{prop:generic-fibre-functoriality}
\leavevmode
\begin{enumerate}
\item[(1)] The functor $(-)_\eta$ commutes with fibre products. In particular, for a formal scheme $\mathfrak Y$ over $\Lambda$ and a finite extension $L/K$, there is a canonical isomorphism of rigid analytic $L$-spaces
\[
\bigl(\mathfrak Y\times_{\mathrm{Spf}(\Lambda)}\mathrm{Spf}(\mathcal O_L)\bigr)_\eta
\simeq
\mathfrak Y_\eta\times_{\mathrm{Spm}(K)}\mathrm{Spm}(L).
\]
\item[(2)] The functor $(-)_\eta$ preserves closed immersions and finite morphisms.
\end{enumerate}
\end{proposition}

\begin{proof}
These properties are established in \cite[Proposition~7.2.4(e), (g) and Section~7.2.6]{Jon95}.
\end{proof}

For each integer $N\ge1$, let
\[
A_{N,0}=\frac{\Lambda\langle x_1,\ldots,x_n,T_1,\ldots,T_n\rangle}
                 {(x_i^N-\varpi T_i\mid 1\le i\le n)},\qquad
A_N=A_{N,0}[1/\varpi],
\]
and write $V_{\pi,N}=\mathrm{Spm}(A_N)$.

A covering $X=\bigcup_{N\ge1}U_N$ by affinoid subdomains with $U_N\subset U_{N+1}$ is \emph{admissible} if every morphism from an affinoid $K$-space to $X$ factors through some $U_N$; see \cite[Section~9.1.4, Proposition~2]{BGR84}.

\begin{proposition}\label{prop:explicit-charts}
\leavevmode
\begin{enumerate}
\item[(1)] The algebra $A_{N,0}$ is finite free over $\Lambda\langle T_1,\ldots,T_n\rangle$ with basis $x_1^{a_1}\cdots x_n^{a_n}$ for $0\le a_i<N$. In particular, $A_{N,0}$ is topologically of finite type and flat over $\Lambda$, and $A_N$ is an affinoid $K$-algebra.
\item[(2)] For every complete valued extension $L/K$,
\[
V_{\pi,N}(L)=\{x\in X_\pi(L)\mid |x_i(x)|^N\le |\varpi|\text{ for }1\le i\le n\}.
\]
\item[(3)] For $M\ge N$, there is an inclusion of affinoid subdomains $V_{\pi,N}\hookrightarrow V_{\pi,M}$. The family $(V_{\pi,N})_{N\ge1}$ gives an admissible covering
\[
X_\pi=\bigcup_{N\ge1}V_{\pi,N}.
\]
\item[(4)] The coordinates $x_1,\ldots,x_n$ identify $X_\pi$ with the open unit polydisc $\{z\in(\mathbf A_K^n)^{\mathrm{an}}\mid |z_i|<1\text{ for all }i\}$. In these coordinates the group law satisfies $(z\cdot w)_i=z_i+w_i+z_iw_i$.
\end{enumerate}
\end{proposition}

\begin{proof}
Consider the homomorphism
\[
R_\pi[T_1,\ldots,T_n]/(x_i^N-\varpi T_i\mid 1\le i\le n)
\longrightarrow R_\pi[1/\varpi],\qquad T_i\longmapsto x_i^N/\varpi.
\]
Its image is the dilatation ring $D_N=R_\pi[x_1^N/\varpi,\ldots,x_n^N/\varpi]$, and it becomes an isomorphism after inverting $\varpi$. Its source is $\varpi$-torsion free. Indeed, the reductions $x_1^N,\ldots,x_n^N$ of the defining relations form a regular sequence modulo $\varpi$. At each successive quotient, if $\varpi a=(x_i^N-\varpi T_i)b$, reduction modulo $\varpi$ gives $b=\varpi c$, and cancelling $\varpi$ gives $a=(x_i^N-\varpi T_i)c$. Thus the homomorphism is injective.

In $D_N/\varpi^cD_N$, the relation $x_i^{Nc}=\varpi^cT_i^c=0$ truncates every power series in the $x_i$ to a polynomial. Reduction modulo $\varpi^c$ also turns the restricted series defining $A_{N,0}$ into polynomials. Both quotients therefore have the same presentation over $(\Lambda/\varpi^c)[T_1,\ldots,T_n]$ with relations $x_i^N-\varpi T_i$. Division by these monic polynomials gives the basis in (1). Taking inverse limits over $c\ge1$ identifies $A_{N,0}$ with the $\varpi$-adic completion $\widehat D_N$ of $D_N$ and gives the asserted basis over $\Lambda\langle T_1,\ldots,T_n\rangle$.

Let $\mathfrak m=(\varpi,x_1,\ldots,x_n)$, $I_N=(\varpi,x_1^N,\ldots,x_n^N)$, and $d=(n+1)N$. Every monomial of degree $d$ in the generators of $\mathfrak m$ is divisible by $\varpi$ or by some $x_i^N$. Hence $\mathfrak m^d\subseteq I_N\subseteq(\varpi,\mathfrak m^N)$. Since $D_N=R_\pi[I_N/\varpi]$, the inclusions $\mathfrak m^dD_N\subseteq\varpi D_N\subseteq\mathfrak mD_N$ identify the two adic topologies on $D_N$. If $U_N$ denotes de Jong's chart for $(\varpi,\mathfrak m^N)$, the same ideal inclusions give $U_N\subseteq V_{\pi,N}\subseteq U_d$. Thus the affinoids $V_{\pi,N}$ are cofinal with the charts defining $X_\pi$ in \cite[Sections~7.1.1--7.1.3]{Jon95}. The relations $x_i^N=\varpi T_i$ uniquely determine $T_i=x_i^N/\varpi$, and the conditions $|T_i|\le1$ are equivalent to $|x_i|^N\le|\varpi|$. Projection onto the coordinates $x_1,\ldots,x_n$ therefore identifies $V_{\pi,N}$ with the affinoid subdomain defined by these inequalities, proving (2). For $M\ge N$, the identity $x_i^M/\varpi=x_i^{M-N}(x_i^N/\varpi)$ gives the transition homomorphism $A_M\to A_N$ over $R_\pi$. The radii $|\varpi|^{1/N}$ increase to $1$, giving the inclusions and union in (3) and the isomorphism in (4).

For admissibility in (3), let $f:Y\to X_\pi$ be a morphism from an affinoid $K$-space. The maximum modulus principle gives a common bound $\rho<1$ for the functions $f^*x_i$ on $Y$. Choose $N$ such that $\rho^N\le|\varpi|$. Then $f$ factors through $V_{\pi,N}$; see also \cite[Section~9.1.4, proof of Proposition~5]{BGR84}. Finally, the group law in (4) follows from $q_i\mapsto q_i\otimes q_i$ and $x_i=q_i-1$.
\end{proof}

\begin{remark}\label{rem:chart-basis-indep}
\leavevmode
\begin{enumerate}
\item[(1)] Each $V_{\pi,N}$ is geometrically integral. Over the finite extension $K(\alpha)/K$ with $\alpha^N=\varpi$, there is an isomorphism $A_N\otimes_K K(\alpha)\xrightarrow{\sim}K(\alpha)\langle S_1,\ldots,S_n\rangle$, $x_i\mapsto\alpha S_i$, $T_i\mapsto S_i^N$.
\item[(2)] The domain $V_{\pi,N}$ is independent of the chosen basis of $\pi$. For $n>0$, a change of basis $e_j'=\sum_i a_{ij}e_i$ with $(a_{ij})\in\mathrm{GL}_n(\mathbf Z_\ell)$ gives $x_j'=\prod_i(1+x_i)^{a_{ij}}-1$. The binomial coefficients lie in $\mathbf Z_\ell$, so $|x_j'(\xi)|\le\max_i|x_i(\xi)|$ for every point $\xi\in X_\pi$. Applying the same estimate to the inverse change of basis gives $\max_j|x_j'(\xi)|=\max_i|x_i(\xi)|$. Proposition~\ref{prop:explicit-charts}(2) therefore gives the same affinoid subdomain for both bases.
\end{enumerate}
\end{remark}

\begin{remark}\label{rem:points-character-torus}
Let $L/K$ be a complete nonarchimedean extension and let $\mathfrak m_L=\{a\in\mathcal O_L\mid |a|<1\}$.
\begin{enumerate}
\item[(1)] Fix a basis $e_1,\ldots,e_n$ of $\pi$. For a point $\xi\in X_\pi(L)$, let $u_i=x_i(\xi)\in\mathfrak m_L$. These coordinate values determine a continuous character
\[
\chi_\xi:\pi\longrightarrow1+\mathfrak m_L,
\qquad
\chi_\xi\!\left(\sum_i a_i e_i\right)=\prod_i(1+u_i)^{a_i},
\]
where the $\ell$-adic binomial series converges because $a_i\in\mathbf Z_\ell$ and $|u_i|<1$. Conversely, let $\chi:\pi\to L^\times$ be continuous. Since $\ell^m e_i\to0$ in $\pi$, one has $\chi(e_i)^{\ell^m}\to1$ and hence $|\chi(e_i)|=1$. Reducing modulo $\mathfrak m_L$, the element $\overline{\chi(e_i)}$ has $\ell$-power order. A field of characteristic $\ell$ has no nontrivial $\ell$-power roots of unity, so $\overline{\chi(e_i)}=1$. Hence $\chi(e_i)-1\in\mathfrak m_L$. The assignments $\xi\mapsto\chi_\xi$ and $\chi\mapsto\xi_\chi$, where $x_i(\xi_\chi)=\chi(e_i)-1$, are inverse and give the group isomorphism
\[
X_\pi(L)\xrightarrow{\sim}\mathrm{Hom}_{\mathrm{cont}}(\pi,1+\mathfrak m_L)
=\mathrm{Hom}_{\mathrm{cont}}(\pi,L^\times).
\]
\item[(2)] Every torsion point of $X_\pi(L)$ has $\ell$-power order, since $1+\mathfrak m_L$ contains no nontrivial root of unity of order prime to $\ell$.
\end{enumerate}
\end{remark}

\begin{proposition}\label{prop:closed-subspace-charts}
For any ideal $I \subset R_\pi$, the rigid generic fibre $(\mathrm{Spf}(R_\pi/I))_\eta$ is a closed rigid subspace of $X_\pi$ whose intersection with each affinoid chart $V_{\pi,N}$ is
\[
\bigl(\mathrm{Spf}(R_\pi/I)\bigr)_\eta \cap V_{\pi,N} = \mathrm{Spm}(A_N / I A_N).
\]
\end{proposition}

\begin{proof}
The closed immersion $\mathrm{Spf}(R_\pi/I)\hookrightarrow\mathfrak X_\pi$ induces a closed immersion on generic fibres by Proposition~\ref{prop:generic-fibre-functoriality}(2). On the chart $V_{\pi,N}$, the dilatation ring for $R_\pi/I$ is the image of $D_N$ in $(R_\pi/I)[1/\varpi]$. The kernel $D_N\cap I R_\pi[1/\varpi]$ consists of the elements $h\in D_N$ such that $\varpi^m h\in ID_N$ for some $m\ge0$. Since $D_N$ is Noetherian, its $\varpi$-adic completion preserves this quotient. After inverting $\varpi$, the kernel and $ID_N$ generate the same ideal $IA_N$ in $A_N$, so the chart algebra is $A_N/IA_N$.
\end{proof}

\subsection{Standard integral models}\label{subsec:standard-models}

\begin{definition}\label{def:standard-subspace}
A closed rigid subspace of $X_\pi$ is \emph{standard} if it has the form $Z_J=(\mathrm{Spf}(R_\pi/J))_\eta$ for a radical $\ell$-saturated ideal $J\subset R_\pi$. Its integral model $\mathrm{Spf}(R_\pi/J)$ is \emph{equidimensional} if $R_\pi/J$ is equidimensional.
\end{definition}

\begin{lemma}\label{lem:chart-flatness}
Let $N\ge1$.
\begin{enumerate}
\item[(1)] The canonical homomorphism $R_\pi\to A_N$ is flat.
\item[(2)] For any finite family of ideals $I_\alpha\subset R_\pi$, one has
\[
\left(\bigcap_\alpha I_\alpha\right)A_N=\bigcap_\alpha I_\alpha A_N.
\]
\end{enumerate}
\end{lemma}

\begin{proof}
The dilatation ring $D_N$ is Noetherian, since it is finitely generated over $R_\pi$. Its $\varpi$-adic completion $D_N\to A_{N,0}$ is therefore flat. Inverting $\varpi$ gives a flat homomorphism $R_\pi[1/\varpi]=D_N[1/\varpi]\to A_N$. Precomposing with the flat localization $R_\pi\to R_\pi[1/\varpi]$ proves (1). For (2), tensor the injection $R_\pi/(I_1\cap I_2)\hookrightarrow R_\pi/I_1\oplus R_\pi/I_2$ with $A_N$ and use induction on the number of ideals.
\end{proof}

\begin{lemma}\label{lem:standard-reduced}
For a standard closed rigid subspace $Z_J\subset X_\pi$, each chart algebra $A_N/JA_N$ is reduced. Hence $Z_J$ is reduced.
\end{lemma}

\begin{proof}
The quotient $R_\pi/J$ is reduced and $\varpi$-torsion free. The image of the dilatation ring $D_N$ in $(R_\pi/J)[1/\varpi]$ is reduced and is of finite type over the excellent ring $R_\pi/J$. Its $\varpi$-adic completion is therefore reduced, since the completion map of an excellent Noetherian ring is regular \cite[Tags~0AH2 and~07QK]{Stacks}. Inverting $\varpi$ in this $\varpi$-adic completion gives $A_N/JA_N$ by Proposition~\ref{prop:closed-subspace-charts}. Thus $A_N/JA_N$ is reduced.
\end{proof}

\begin{lemma}\label{lem:jacobson}
Let $I\subset R_\pi$ be an ideal.
\begin{enumerate}
\item[(1)] The ring $(R_\pi/I)[1/\varpi]$ is Jacobson.
\item[(2)] Every maximal ideal $\mathfrak m$ of $(R_\pi/I)[1/\varpi]$ has residue field $\kappa(\mathfrak m)$ finite over $K$. Write $\bar x_i$ for the image of $x_i$ in $\kappa(\mathfrak m)$. Then $|\bar x_i|<1$ for all $i$, and $\mathfrak m$ is the kernel of the continuous homomorphism $(R_\pi/I)[1/\varpi]\to\kappa(\mathfrak m)$ obtained by evaluation at $(\bar x_1,\dots,\bar x_n)$.
\end{enumerate}
\end{lemma}

\begin{proof}
The ring $R_\pi[1/\varpi]=K\otimes_\Lambda\Lambda[[x_1,\ldots,x_n]]$ is the ring of separated power series in $n$ open polydisc variables of \cite[Definition~2.1.1]{Lip00}. By the Nullstellensatz \cite[Theorem~4.1.1]{Lip00}, its radical ideals are intersections of maximal ideals, and the residue fields of its maximal ideals are finite extensions of $K$. In particular the ring and its quotients are Jacobson.

The maximal ideals of $R_\pi[1/\varpi]$ are exactly the kernels of evaluation at rigid points of the open unit polydisc \cite[Proposition~4.1.2(i)]{Lip00}, \cite[Lemma~7.1.9]{Jon95}. Such an evaluation factors through $(R_\pi/I)[1/\varpi]$ precisely when every element of $I$ vanishes at the point. This proves (2).
\end{proof}

\begin{corollary}\label{cor:standard-point-ideals}
Let $J\subset R_\pi$ be a radical $\ell$-saturated ideal.
\begin{enumerate}
\item[(1)] If $J\ne R_\pi$, then $Z_J(L)\ne\varnothing$ for some finite extension $L/K$.
\item[(2)] The ideal $J[1/\varpi]$ is the intersection of the evaluation kernels in $R_\pi[1/\varpi]$ at all points of $Z_J$ over finite extensions of $K$.
\end{enumerate}
\end{corollary}

\begin{proof}
For $J=R_\pi$, the empty intersection in (2) is $R_\pi[1/\varpi]$. Assume $J$ is proper. By $\ell$-saturation, $(R_\pi/J)[1/\varpi]$ is nonzero. Choose a maximal ideal $\mathfrak m$ of $(R_\pi/J)[1/\varpi]$. By Lemma~\ref{lem:jacobson}(2), $\kappa(\mathfrak m)/K$ is finite and the images of $x_1,\ldots,x_n$ in $\kappa(\mathfrak m)$ define a point of $Z_J(\kappa(\mathfrak m))$. This proves (1). Since the ring is reduced and Jacobson by Lemma~\ref{lem:jacobson}(1), the intersection of its maximal ideals is zero. Their inverse images in $R_\pi[1/\varpi]$ are precisely the evaluation kernels in (2).
\end{proof}

For finitely many closed rigid subspaces $Z_1,\dots,Z_r$ of a rigid space, with coherent defining ideals $\mathcal{I}_1,\dots,\mathcal{I}_r$, their \emph{closed analytic union}, denoted by $\bigcup_{\nu=1}^r{}^{\mathrm{an}}Z_\nu$, is the closed rigid subspace defined by $\bigcap_{\nu=1}^r\mathcal{I}_\nu$.

\begin{lemma}\label{lem:analytic-union-components}
Let $J\subset R_\pi$ be a radical $\ell$-saturated ideal with minimal prime decomposition $J=\bigcap_{\nu=1}^{r}\mathfrak p_\nu$.
\begin{enumerate}
\item[(1)] Each $\mathfrak p_\nu$ is $\ell$-saturated and defines a nonempty standard closed rigid subspace $Z_{\mathfrak p_\nu}$.
\item[(2)] The closed analytic union of these subspaces is $Z_J=\bigcup_{\nu=1}^{r}{}^{\mathrm{an}} Z_{\mathfrak p_\nu}$.
\end{enumerate}
\end{lemma}

\begin{proof}
Since $R_\pi/J$ is reduced and Noetherian, its minimal primes are associated primes. The nonzerodivisor $\ell$ thus lies in none of the $\mathfrak p_\nu$. Each $\mathfrak p_\nu$ is $\ell$-saturated, and Corollary~\ref{cor:standard-point-ideals}(1) gives a point of $Z_{\mathfrak p_\nu}$. This proves (1). For (2), Lemma~\ref{lem:chart-flatness}(2) gives $JA_N=\bigcap_{\nu=1}^r\mathfrak p_\nu A_N$ on every affinoid chart. The right hand side defines the closed analytic union on that chart.
\end{proof}

\begin{lemma}\label{lem:faithfulness}
Let $I,J\subset R_\pi$ be radical $\ell$-saturated ideals. Then
\[
Z_J\subset Z_I\quad\Longleftrightarrow\quad I\subset J.
\]
In particular, every standard closed rigid subspace of $X_\pi$ has a unique radical $\ell$-saturated defining ideal in $R_\pi$.
\end{lemma}

\begin{proof}
If $I\subset J$, then $IA_N\subset JA_N$ for every $N$, so $Z_J\subset Z_I$ by Proposition~\ref{prop:closed-subspace-charts}. Conversely, if $Z_J\subset Z_I$, every element of $I$ vanishes at the points of $Z_J$ over finite extensions of $K$. Corollary~\ref{cor:standard-point-ideals}(2) gives $I\subset J[1/\varpi]\cap R_\pi=J$, where the equality uses $\ell$-saturation. Taking $Z_J=Z_I$ gives uniqueness, and Definition~\ref{def:standard-subspace} gives existence.
\end{proof}

Fix a $\mathbf Z_\ell$-basis $e_1,\ldots,e_n$ of $\pi$ and positive integers $t_1,\ldots,t_n$. Let $U:\pi\to\pi$ be the linear map $e_i\mapsto\ell^{t_i}e_i$. Precomposition with $U$ induces a formal endomorphism $\Phi:\mathfrak X_\pi\to\mathfrak X_\pi$. The induced endomorphism of $X_\pi$, also denoted by $\Phi$, and the formal pullback are
\[
\begin{aligned}
\Phi &: X_\pi\longrightarrow X_\pi,
& \chi&\longmapsto\chi\circ U,\\
\Phi^\# &: R_\pi\longrightarrow R_\pi,
& q_i&\longmapsto q_i^{\ell^{t_i}}.
\end{aligned}
\]
For a finite $R_\pi$-algebra $S$, the inverse image of $V_{\pi,N}$ in $(\mathrm{Spf}(S))_\eta$ is $\mathrm{Spm}(A_N\otimes_{R_\pi}S)$ \cite[Section~7.1.11 and Lemma~7.2.2]{Jon95}. For $S=R_\pi$ with structure map $\Phi^\#$, the algebra $A_N\otimes_{R_\pi,\Phi^\#}R_\pi$ is free of rank $\ell^{t_1+\cdots+t_n}$ over $A_N$ by Remark~\ref{rem:finite-locally-free}. Thus the induced endomorphism of $X_\pi$ is finite locally free.

Recall that the \emph{scheme theoretic image} of a finite morphism of rigid spaces $g:Y\to X$ is the closed rigid subspace defined by $\ker(\mathcal O_X\to g_*\mathcal O_Y)$.

\begin{lemma}\label{lem:rigid-images}
Let $J\subset R_\pi$ be radical and $\ell$-saturated. For $a\ge0$ let $J_a=((\Phi^\#)^a)^{-1}(J)$.
\begin{enumerate}
\item[(1)] The ideal $J_a$ is radical and $\ell$-saturated.
\item[(2)] The rigid and formal scheme theoretic images satisfy
\[
\Phi^a(Z_J)=Z_{J_a}=\bigl(\Phi^a(\mathrm{Spf}(R_\pi/J))\bigr)_\eta.
\]
\end{enumerate}
\end{lemma}

\begin{proof}
The composite $R_\pi\xrightarrow{(\Phi^\#)^a}R_\pi\twoheadrightarrow R_\pi/J$ has kernel $J_a$. Thus $R_\pi/J_a$ embeds in the reduced $\ell$-torsion free ring $R_\pi/J$, which proves (1).

For (2), the restriction $\Phi^a|_{Z_J}:Z_J\to X_\pi$ is finite. The inverse image of $V_{\pi,N}=\mathrm{Spm}(A_N)$ under this restriction is therefore affinoid. Its coordinate algebra and the pullback from $A_N$ appear in the bottom row of the base change diagram
\[
\begin{tikzcd}[column sep=large,row sep=normal]
R_\pi \arrow[r, "{h\mapsto(\Phi^\#)^a(h)\bmod J}"] \arrow[d]
  & R_\pi/J \arrow[d] \\
A_N \arrow[r]
  & A_N\otimes_{R_\pi,(\Phi^\#)^a}(R_\pi/J)\mathrlap{.}
\end{tikzcd}
\]
The ordinary tensor product suffices since $R_\pi$ is finite over $(\Phi^\#)^a(R_\pi)$ by Remark~\ref{rem:finite-locally-free}. The upper horizontal map has kernel $J_a$, so flatness of $R_\pi\to A_N$ (Lemma~\ref{lem:chart-flatness}(1)) gives kernel $J_aA_N$ for the lower horizontal map. These kernels define $Z_{J_a}$, and $J_a$ defines the formal image by Lemma~\ref{lem:image-ideals}.
\end{proof}

\subsection{The rigid dynamical criterion}\label{subsec:rigid-dynamical}

A reduced rigid space is called \emph{irreducible} if it is nonempty and cannot be written as the union of two proper closed rigid subspaces.

\begin{definition}\label{def:rigid-subtorus}
A \emph{rigid subtorus} of $X_\pi$ is the generic fibre $H_M=(T_M)_\eta=X_{\pi/M}$ of a formal subtorus, where $M\subset\pi$ is primitive. It is \emph{$\Phi$-stable} if $\Phi(H_M)\subset H_M$.
\end{definition}

\begin{lemma}\label{lem:rigid-subtorus-irred}
The character torus $X_\pi$, every rigid subtorus $H_M \subset X_\pi$, and every translate of $H_M$ over a finite extension of $K$ are geometrically irreducible.
\end{lemma}

\begin{proof}
Over any complete extension of $K$, the charts $V_{\pi,N}$ remain integral by Remark~\ref{rem:chart-basis-indep}(1). Suppose that two closed analytic subspaces cover the character torus after this base change. Each chart lies in one of them. Since the charts are increasing, one of the two subspaces contains charts of arbitrarily large index and hence the whole character torus. This proves geometric irreducibility of $X_\pi$. Since $\pi/M$ is finite free over $\mathbf Z_\ell$, the character torus $H_M=X_{\pi/M}$ is also geometrically irreducible. Translation over its field of definition is an isomorphism, so it preserves geometric irreducibility.
\end{proof}

\begin{remark}\label{rem:rigid-subtorus-stability}
For a primitive submodule $M\subset\pi$, one has
\[
\Phi(H_M)\subset H_M
\quad\Longleftrightarrow\quad
\Phi(T_M)\subset T_M
\quad\Longleftrightarrow\quad
U(M)\subset M.
\]
The first equivalence follows from Lemmas~\ref{lem:faithfulness} and~\ref{lem:rigid-images}(2), and the second from Lemma~\ref{lem:subtorus-stability}.
\end{remark}

\begin{lemma}\label{lem:rigid-torsion-translate}
Let $L/K$ be a finite extension with valuation ring $\mathcal O_L$, let $M\subset\pi$ be primitive, and let $\zeta\in X_\pi(L)$ be a torsion point. The translate $t_\zeta((H_M)_L)\subset X_{\pi,L}$ is a standard closed rigid subspace over $L$, defined by the prime $\ell$-saturated ideal
\[
I_{\zeta,M} = ([\gamma] - \zeta(\gamma) \mid \gamma \in M) \subset \mathcal{O}_L[[\pi]].
\]
\end{lemma}

\begin{proof}
The values $\zeta(\gamma)$ lie in $1+\mathfrak m_L$ by Remark~\ref{rem:points-character-torus}(1), so $\zeta$ is an $\mathcal O_L$-point of $\mathfrak X_\pi$. Lemma~\ref{lem:torsion-translate-ideal} gives $I_{\zeta,M}$ as the ideal of the formal translate $t_\zeta((T_M)_{\mathcal O_L})$.

Extend a basis $e_1,\ldots,e_d$ of $M$ to a basis $e_1,\ldots,e_n$ of $\pi$. In the resulting additive coordinates $x_i=[e_i]-1$, the equations are $x_i=\zeta(e_i)-1$ for $1\le i\le d$. Thus
\[
\mathcal O_L[[\pi]]/I_{\zeta,M}\simeq\mathcal O_L[[x_{d+1},\ldots,x_n]],
\]
which is an $\ell$-torsion free domain. The formal translate has rigid generic fibre $\bigl(t_\zeta((T_M)_{\mathcal O_L})\bigr)_\eta=t_\zeta((H_M)_L)$.
\end{proof}

\Needspace{4\baselineskip}
\begin{definition}\label{def:rigid-quasilinear}
Let $Z\subset X_\pi$ be a reduced closed rigid subspace.
\begin{enumerate}
\item[(1)] The subspace $Z$ is called \emph{rigid quasilinear} if, after a finite extension $L/K$, it admits a decomposition into a finite closed analytic union
\[
Z_L={\bigcup_{\alpha=1}^r}\!^{\mathrm{an}}
t_{\zeta_\alpha}\bigl((H_{M_\alpha})_L\bigr),
\]
where $\zeta_\alpha\in X_\pi(L)$ are torsion points and $M_\alpha\subset\pi$ are primitive submodules.
\item[(2)] The subspace $Z$ is called \emph{$\Phi$-stable rigid quasilinear} if the decomposition in (1) can be chosen so that $\Phi(H_{M_\alpha})\subseteq H_{M_\alpha}$ for every $\alpha$.
\end{enumerate}
\end{definition}

\begin{lemma}\label{lem:rigid-coeff-ext}
Let $L/K$ be a finite extension with valuation ring $\mathcal O_L$. If $Z_J\subset X_\pi$ is standard, then $(Z_J)_L$ is standard over $\mathcal O_L$, defined by $J\mathcal O_L[[\pi]]$, which is radical and $\ell$-saturated.
\end{lemma}

\begin{proof}
Since $\mathcal O_L$ is finite free over $\Lambda$, there is an isomorphism
\[
\mathcal O_L[[\pi]]/J\mathcal O_L[[\pi]]
\cong (R_\pi/J)\otimes_\Lambda\mathcal O_L.
\]
This quotient is flat over $\mathcal O_L$, hence embeds in its localization $(R_\pi/J)[1/\varpi]\otimes_K L$ at $\varpi$.

The reduced ring $(R_\pi/J)[1/\varpi]$ embeds in the product of the fraction fields of its minimal prime quotients. Tensoring with the finite separable extension $L/K$ preserves the injection and gives a product of fields on the right. Thus $(R_\pi/J)[1/\varpi]\otimes_K L$ is reduced.

The quotient $\mathcal O_L[[\pi]]/J\mathcal O_L[[\pi]]$ is therefore reduced and flat over $\mathcal O_L$, so $J\mathcal O_L[[\pi]]$ is radical and $\ell$-saturated. The generic fibre of $\operatorname{Spf}(\mathcal O_L[[\pi]]/J\mathcal O_L[[\pi]])$ is $(Z_J)_L$ by Proposition~\ref{prop:generic-fibre-functoriality}(1).
\end{proof}

\begin{lemma}\label{lem:formal-rigid-quasilinear-equiv}
For a standard closed rigid subspace $Z_J \subset X_\pi$ with formal model $\mathfrak Z_J = \mathrm{Spf}(R_\pi/J)$, $Z_J$ is $\Phi$-stable rigid quasilinear if and only if $\mathfrak Z_J$ is $\Phi$-stable quasilinear.
\end{lemma}

\begin{proof}
Let $L/K$ be a finite extension, let $\zeta_1,\ldots,\zeta_r\in X_\pi(L)$ be torsion points, and let $M_1,\ldots,M_r\subset\pi$ be primitive. The scheme theoretic union of the translates $t_{\zeta_\alpha}\bigl((T_{M_\alpha})_{\mathcal O_L}\bigr)$ equals $(\mathfrak Z_J)_{\mathcal O_L}$ precisely when
\[
J\mathcal O_L[[\pi]]=\bigcap_{\alpha=1}^r I_{\zeta_\alpha,M_\alpha}.
\]
By Lemma~\ref{lem:rigid-coeff-ext}, the ideal on the left defines $(Z_J)_L$. By Lemmas~\ref{lem:rigid-torsion-translate} and~\ref{lem:chart-flatness}(2), the intersection on the right defines the closed analytic union of the rigid translates. Both ideals are radical and $\ell$-saturated, so Lemma~\ref{lem:faithfulness} identifies the displayed ideal equality with the equality of closed rigid subspaces
\[
(Z_J)_L=\bigcup_{\alpha=1}^r{}^{\mathrm{an}}
t_{\zeta_\alpha}\bigl((H_{M_\alpha})_L\bigr).
\]
Finally $T_{M_\alpha}$ is $\Phi$-stable if and only if $H_{M_\alpha}$ is $\Phi$-stable by Remark~\ref{rem:rigid-subtorus-stability}.
\end{proof}

\begin{theorem}[Rigid dynamical criterion]\label{thm:rigid-dynamical}
Let $Z_J\subset X_\pi$ be a standard closed rigid subspace with $R_\pi/J$ equidimensional.
\begin{enumerate}
\item[(1)] $Z_J$ is $\Phi$-stable rigid quasilinear if and only if $\Phi^a(Z_J)\subset\Phi^b(Z_J)$ as closed rigid subspaces for some integers $a>b\ge0$.
\item[(2)] Moreover, $Z_J$ is a $\Phi$-stable rigid subtorus if and only if it is irreducible and $\Phi^c(Z_J)\subset Z_J$ for some integer $c\ge1$.
\end{enumerate}
\end{theorem}

\begin{proof}
For (1), let $\mathfrak Z_J=\mathrm{Spf}(R_\pi/J)$. Since $J$ is radical and $\ell$-saturated, $\mathfrak Z_J$ is reduced and flat over $\Lambda$. It is equidimensional by hypothesis. By Lemma~\ref{lem:formal-rigid-quasilinear-equiv}, $Z_J$ is $\Phi$-stable rigid quasilinear if and only if $\mathfrak Z_J$ is $\Phi$-stable quasilinear. By Lemmas~\ref{lem:rigid-images}(2) and~\ref{lem:faithfulness}, the inclusion $\Phi^a(Z_J)\subset\Phi^b(Z_J)$ is equivalent to $\Phi^a(\mathfrak Z_J)\subset\Phi^b(\mathfrak Z_J)$. Proposition~\ref{prop:formal-dynamical-criterion}(1) applied to $\mathfrak Z_J$ now gives (1).

For (2), a $\Phi$-stable rigid subtorus is irreducible by Lemma~\ref{lem:rigid-subtorus-irred} and satisfies the inclusion with $c=1$. Conversely, suppose that $Z_J$ is irreducible and $\Phi^c(Z_J)\subset Z_J$. Write $J=\bigcap_{\nu=1}^r\mathfrak p_\nu$ for its minimal prime decomposition. By Lemma~\ref{lem:analytic-union-components}(2), the closed subspaces $Z_{\mathfrak p_\nu}$ cover $Z_J$. Since $Z_J$ is irreducible, $Z_J=Z_{\mathfrak p_\nu}$ for some $\nu$. Lemma~\ref{lem:faithfulness} then gives $J=\mathfrak p_\nu$. Thus $\mathfrak Z_J$ is irreducible. The comparison of images in (1) gives $\Phi^c(\mathfrak Z_J)\subset\mathfrak Z_J$. By Proposition~\ref{prop:formal-dynamical-criterion}(2), $\mathfrak Z_J$ is a $\Phi$-stable formal subtorus. Its generic fibre is the required $\Phi$-stable rigid subtorus.
\end{proof}

For a standard closed rigid subspace $Z_J\subset X_\pi$, let $J=\bigcap_{\nu=1}^r\mathfrak p_\nu$ be the minimal prime decomposition. For each $d\ge0$, write
\[
J^{(d)}=\bigcap_{\dim(R_\pi/\mathfrak p_\nu)=d}\mathfrak p_\nu.
\]
An empty intersection is the unit ideal and gives $Z_{J^{(d)}}=\varnothing$.

\begin{corollary}\label{cor:fixed-dim-components}
Let $Z_J\subset X_\pi$ be standard. Suppose that for every $d$ with $Z_{J^{(d)}}\ne\varnothing$ there are integers $a_d>b_d\ge0$ such that
\[
\Phi^{a_d}(Z_{J^{(d)}})\subset\Phi^{b_d}(Z_{J^{(d)}}).
\]
Then $Z_J$ is $\Phi$-stable rigid quasilinear.
\end{corollary}

\begin{proof}
For each nonempty $Z_{J^{(d)}}$, Lemma~\ref{lem:analytic-union-components}(1) shows that $J^{(d)}$ is radical and $\ell$-saturated. The quotient $R_\pi/J^{(d)}$ is equidimensional of dimension $d$ by construction. Thus Theorem~\ref{thm:rigid-dynamical}(1) gives a quasilinear decomposition with $\Phi$-stable subtori for $Z_{J^{(d)}}$. The finitely many nonempty $Z_{J^{(d)}}$ have closed analytic union $Z_J$, since $J=\bigcap_d J^{(d)}$ and Lemma~\ref{lem:chart-flatness}(2) preserves this intersection on every chart. Taking the union of their decompositions over a common finite extension of $K$ shows that $Z_J$ is $\Phi$-stable rigid quasilinear.
\end{proof}

\subsection{Descent along finite isogenies}\label{subsec:finite-isogenies}

\begin{definition}\label{def:finite-isogeny}
A homomorphism of rigid character tori is called a \emph{finite isogeny} if it is finite and surjective.
\end{definition}

Let $\pi\subset\pi'$ be full $\mathbf Z_\ell$-lattices in a finite dimensional $\mathbf Q_\ell$-vector space. Restriction of characters along $\pi\subset\pi'$ defines $f:X_{\pi'}\to X_\pi$.

\begin{lemma}\label{lem:isogeny-properties}
\begin{enumerate}
\item[(1)] The algebra $R_{\pi'}$ is finite free of positive rank over $R_\pi$, and $f$ is a finite étale isogeny.
\item[(2)] The kernel $\ker(f)$ is a finite étale rigid group over $K$ with torsion geometric points. There is a finite extension $L/K$ such that $\ker(f)_L$ is a constant finite rigid group over $L$.
\end{enumerate}
\end{lemma}

\begin{proof}
Choose bases in Smith normal form, so that $e_i=\ell^{a_i}e_i'$ with $e_i\in\pi$, $e_i'\in\pi'$ and $a_i\ge0$. In the coordinates $x_i=[e_i]-1$ and $x_i'=[e_i']-1$, the pullback $f^\#:R_\pi\to R_{\pi'}$ satisfies $f^\#(1+x_i)=(1+x_i')^{\ell^{a_i}}$. Successive Weierstrass division makes $R_{\pi'}$ free of rank $\prod_i\ell^{a_i}$ over $R_\pi$. The coordinate algebra of $f^{-1}(V_{\pi,N})$ is $A_N\otimes_{R_\pi}R_{\pi'}$, which is free of rank $[\pi':\pi]$ over $A_N$. Thus $f$ is finite and faithfully flat. After inverting $\ell$, each derivative $\ell^{a_i}(1+x_i')^{\ell^{a_i}-1}$ is a unit, so $R_\pi[1/\ell]\to R_{\pi'}[1/\ell]$ is finite étale. This proves (1).

In the same coordinates the kernel is $\prod_i\mu_{\ell^{a_i}}$. It is finite étale over $K$ and becomes constant after adjoining the indicated roots of unity. This proves (2).
\end{proof}

\begin{lemma}\label{lem:isogeny-subtorus-images}
Let $M'\subset\pi'$ be primitive.
\begin{enumerate}
\item[(1)] The image $f(H_{M'})$ is the rigid subtorus $H_{\pi\cap M'}$ of $X_\pi$.
\item[(2)] For every finite extension $L/K$ and torsion point $\zeta\in X_{\pi'}(L)$, one has $f(t_\zeta(H_{M'}))=t_{f(\zeta)}(H_{\pi\cap M'})$ over $L$.
\end{enumerate}
\end{lemma}

\begin{proof}
Let $M=\pi\cap M'$. The inclusion $\pi/M\hookrightarrow\pi'/M'$ has finite index. Since $\pi'/M'$ is free, $M$ is primitive. Restriction of characters along $\pi/M\hookrightarrow\pi'/M'$ gives the left vertical arrow in the commutative diagram
\[
\begin{tikzcd}[column sep=1.8em,row sep=normal]
H_{M'} \arrow[r, hook] \arrow[d, "{f|_{H_{M'}}}"']
  & X_{\pi'} \arrow[d, "f"] \\
H_M \arrow[r, hook] & X_\pi\mathrlap{.}
\end{tikzcd}
\]
By Lemma~\ref{lem:isogeny-properties}(1), the map $H_{M'}\to H_M$ is a finite étale isogeny, so $f(H_{M'})=H_M$. After base change to $L$, the identity $f\circ t_\zeta=t_{f(\zeta)}\circ f$ gives (2).
\end{proof}

\begin{lemma}\label{lem:isogeny-standard-pullback}
For every standard closed rigid subspace $Z_J\subset X_\pi$, the diagram
\[
\begin{tikzcd}[column sep=1.8em,row sep=normal]
Z_{JR_{\pi'}} \arrow[r, hook] \arrow[d]
  & X_{\pi'} \arrow[d, "f"] \\
Z_J \arrow[r, hook] & X_\pi
\end{tikzcd}
\]
is Cartesian, and $Z_{JR_{\pi'}}$ is reduced and standard.
\end{lemma}

\begin{proof}
By Lemma~\ref{lem:isogeny-properties}(1), the ring $R_{\pi'}/JR_{\pi'}$ is finite free over $R_\pi/J$, hence $\ell$-torsion free. After inverting $\ell$ it is finite étale over the reduced ring $(R_\pi/J)[1/\ell]$, so it is reduced. The injection into this localization shows that $R_{\pi'}/JR_{\pi'}$ is reduced as well. Thus $JR_{\pi'}$ is radical and $\ell$-saturated. The formal inverse image of $\mathrm{Spf}(R_\pi/J)$ in $\mathfrak X_{\pi'}$ is $\mathrm{Spf}(R_{\pi'}/JR_{\pi'})$. Proposition~\ref{prop:generic-fibre-functoriality}(1) gives the Cartesian square on rigid generic fibres. Reducedness follows from Lemma~\ref{lem:standard-reduced}.
\end{proof}

\begin{lemma}\label{lem:descent-quasilinear}
Let $f : X_{\pi'} \to X_\pi$ be the finite isogeny induced by an inclusion of lattices $\pi \subset \pi'$ with $[\pi':\pi]<\infty$, and let $Z \subset X_\pi$ be a reduced closed rigid subspace. If $(f^{-1}(Z))_{\mathrm{red}}$ is rigid quasilinear, then $Z$ is rigid quasilinear.
\end{lemma}

\begin{proof}
Choose a finite extension $L/K$ over which $((f^{-1}(Z))_{\mathrm{red}})_L$ is a closed analytic union of torsion translates $t_{\zeta_\alpha}\bigl((H_{M_\alpha'})_L\bigr)$, and let $f_L:X_{\pi',L}\to X_{\pi,L}$ be the base change of $f$. Lemma~\ref{lem:isogeny-subtorus-images}(2) gives
\[
f_L\bigl(t_{\zeta_\alpha}((H_{M_\alpha'})_L)\bigr)
=t_{f_L(\zeta_\alpha)}\bigl((H_{\pi\cap M_\alpha'})_L\bigr).
\]
Since $L/K$ is finite separable, reduction commutes with this extension, so $((f^{-1}(Z))_{\mathrm{red}})_L=(f_L^{-1}(Z_L))_{\mathrm{red}}$. The base change $f_L^{-1}(Z_L)\to Z_L$ of the finite isogeny $f_L$ is surjective, and reduction preserves geometric points. Hence the closed analytic union of the images above has the same geometric points as $Z_L$. This union is reduced because its defining ideal on each chart is an intersection of radical ideals. The space $Z_L$ is reduced since $Z$ is reduced and $L/K$ is separable. The affinoid Nullstellensatz \cite[Section~7.1.2]{BGR84} therefore gives equality on every chart $V_{\pi,N,L}$. This is a rigid quasilinear decomposition of $Z_L$.
\end{proof}

\begin{corollary}\label{cor:isogenous-rigid-criterion}
Let $\pi\subset\pi'$ be full lattices, and let $f:X_{\pi'}\to X_\pi$ be the induced finite isogeny. Fix a basis $e_1',\ldots,e_n'$ of $\pi'$ and positive integers $t_1',\ldots,t_n'$. Let $\Phi:X_{\pi'}\to X_{\pi'}$ be the power endomorphism induced by
\[
\Phi^\#:R_{\pi'}\longrightarrow R_{\pi'},
\qquad
[e_i']\longmapsto[e_i']^{\ell^{t_i'}}.
\]
Let $Z_J\subset X_\pi$ be a standard closed rigid subspace and let $Y=f^{-1}(Z_J)=Z_{JR_{\pi'}}\subset X_{\pi'}$. Suppose that for every $d$ with $Z_{(JR_{\pi'})^{(d)}}\ne\varnothing$ there are integers $a_d>b_d\ge0$ such that
\[
\Phi^{a_d}\bigl(Z_{(JR_{\pi'})^{(d)}}\bigr)
\subset
\Phi^{b_d}\bigl(Z_{(JR_{\pi'})^{(d)}}\bigr).
\]
Then $Y$ is $\Phi$-stable rigid quasilinear, and $Z_J$ is rigid quasilinear.
\end{corollary}

\begin{proof}
By Lemma~\ref{lem:isogeny-standard-pullback}, $Y=f^{-1}(Z_J)=Z_{JR_{\pi'}}$ is a standard closed rigid subspace of $X_{\pi'}$. Corollary~\ref{cor:fixed-dim-components} applied to $Y$ shows that $Y$ is $\Phi$-stable rigid quasilinear, and Lemma~\ref{lem:descent-quasilinear} then implies that $Z_J$ is rigid quasilinear.
\end{proof}

\section{Frobenius stable character loci}\label{sec:frobenius-loci}

This section establishes the Frobenius invariance criterion and constructs the finite isogeny $X_{\pi'}\to X_\pi$ used in Theorem~\ref{thm:main-fm} in Section~\ref{sec:jump-loci}.

In this section, let $\Lambda=\mathcal O_E$ for a finite extension $E/\mathbf Q_\ell$ in $\overline{\mathbf Q}_\ell$. For a finite free $\mathbf Z_\ell$-module $\pi$, write $R_\pi=\mathcal O_E[[\pi]]$, $\mathfrak X_\pi=\mathrm{Spf}(R_\pi)$, and $X_\pi=(\mathfrak X_\pi)_\eta$. Fix a field embedding $\iota:\overline{\mathbf Q}_\ell\hookrightarrow\mathbf C$.

\subsection{From Frobenius invariance to \texorpdfstring{$\Phi_w$}{Phi\_w} stability}\label{subsec:local-frob-criterion}
Let $V_\pi=\pi\otimes_{\mathbf Z_\ell}\mathbf Q_\ell$.

\begin{assumption}\label{ass:frobenius-weights}
\leavevmode
\begin{enumerate}
\item[(1)] The automorphism $\sigma\in\mathrm{GL}_{\mathbf Q_\ell}(V_\pi)$ is semisimple and satisfies $\sigma(\pi)=\pi$.
\item[(2)] There is a $\mathbf Q_\ell$-linear decomposition $V_\pi=\bigoplus_{r\in W}V_r$, indexed by a finite set $W\subset\mathbf Z_{>0}$, with $\sigma(V_r)=V_r$ for every $r\in W$. For some real number $c>0$ with $c\ne1$, every eigenvalue $\alpha$ of $\sigma$ on $V_r$ satisfies $|\iota(\alpha)|=c^r$.
\item[(3)] The lattice $\pi$ is the direct sum of the lattices $\pi_r=\pi\cap V_r$ for $r\in W$.
\end{enumerate}
\end{assumption}

Under Assumption~\ref{ass:frobenius-weights}, the action of $\sigma$ on $\pi$ induces the $\mathcal O_E$-algebra automorphism $\sigma^\#:R_\pi\to R_\pi$, $[\gamma]\mapsto[\sigma(\gamma)]$. We also write $\sigma$ for the induced automorphism of $\mathfrak X_\pi$ and for its rigid generic fibre $\sigma:X_\pi\to X_\pi$, $\chi\mapsto\chi\circ\sigma$.

Choose a $\mathbf Z_\ell$-basis $e_1,\dots,e_b$ of $\pi$ by taking a basis of each $\pi_r$. Write $w=(w_1,\dots,w_b)$, with $e_i\in\pi_{w_i}$, and let $q_i=[e_i]$ and $x_i=q_i-1$. The mixed power endomorphism $\Phi_w:\mathfrak X_\pi\to\mathfrak X_\pi$ has pullback
\[
\Phi_w^\#:R_\pi\longrightarrow R_\pi,
\qquad
q_i\longmapsto q_i^{\ell^{w_i}}.
\]
Equivalently, $\Phi_w^\#([\gamma])=[\ell^r\gamma]=[\gamma]^{\ell^r}$ for $\gamma\in\pi_r$. Since $\sigma$ preserves the summands $\pi_r$, it commutes with $\Phi_w$.

Let $l_1,\dots,l_b$ be formal variables. We use the series $\log(1+u)=\sum_{n\ge1}(-1)^{n+1}u^n/n$ and $\exp(u)=\sum_{n\ge0}u^n/n!$.

\begin{proposition}\label{prop:log-coords}
\leavevmode
\begin{enumerate}
\item[(1)] There are mutually inverse continuous $E$-algebra isomorphisms
\[
\begin{aligned}
\mathrm{Exp}_E^\# &: E[[x_1,\dots,x_b]]\xrightarrow{\sim}E[[l_1,\dots,l_b]],
& x_i&\mapsto\exp(l_i)-1,\\
\mathrm{Log}_E^\# &: E[[l_1,\dots,l_b]]\xrightarrow{\sim}E[[x_1,\dots,x_b]],
& l_i&\mapsto\log(1+x_i).
\end{aligned}
\]
\item[(2)] For positive integers $a_1,\dots,a_b$, the following diagram of continuous $E$-algebra homomorphisms commutes:
\[
\begin{CD}
E[[x_1,\dots,x_b]] @>{q_i\mapsto q_i^{a_i}}>> E[[x_1,\dots,x_b]] \\
@V{\mathrm{Exp}_E^\#}V{\rotatebox[origin=c]{90}{$\scriptstyle\sim$}}V @V{\mathrm{Exp}_E^\#}V{\rotatebox[origin=c]{90}{$\scriptstyle\sim$}}V \\
E[[l_1,\dots,l_b]] @>{l_i\mapsto a_i l_i}>> E[[l_1,\dots,l_b]]\mathrlap{.}
\end{CD}
\]
\end{enumerate}
\end{proposition}

\begin{proof}
The series $\exp(l_i)-1$ and $\log(1+x_i)$ have zero constant term, so the substitutions in (1) are continuous. The identities $\log(\exp(l_i))=l_i$ and $\exp(\log(1+x_i))=1+x_i$ give (1), and $\exp(a_i l_i)=\exp(l_i)^{a_i}$ gives (2).
\end{proof}

Over a complete valued extension of $\mathbf Q_\ell$, the series $\log(1+u)$ converges for $|u|<1$ and the series $\exp(u)$ converges for $|u|<\ell^{-1/(\ell-1)}$.

\begin{proposition}\label{prop:analytic-log-coords}
Let $\lambda\in E^\times$ with $\rho=|\lambda|<\ell^{-1/(\ell-1)}$. There are mutually inverse bounded $E$-algebra isomorphisms
\[
\begin{aligned}
\mathrm{Exp}_\rho^\# &: E\langle\lambda^{-1}x_1,\dots,\lambda^{-1}x_b\rangle
\xrightarrow{\sim} E\langle\lambda^{-1}l_1,\dots,\lambda^{-1}l_b\rangle,
& \lambda^{-1}x_i&\mapsto\frac{\exp(l_i)-1}{\lambda},\\
\mathrm{Log}_\rho^\# &: E\langle\lambda^{-1}l_1,\dots,\lambda^{-1}l_b\rangle
\xrightarrow{\sim} E\langle\lambda^{-1}x_1,\dots,\lambda^{-1}x_b\rangle,
& \lambda^{-1}l_i&\mapsto\frac{\log(1+x_i)}{\lambda}.
\end{aligned}
\]
In particular,
\[
\mathrm{Exp}_\rho^\#\bigl(\mathcal O_E\langle\lambda^{-1}x_1,\dots,\lambda^{-1}x_b\rangle\bigr)
=\mathcal O_E\langle\lambda^{-1}l_1,\dots,\lambda^{-1}l_b\rangle.
\]
\end{proposition}

\begin{proof}
Normalize $v_\ell$ by $v_\ell(\ell)=1$. In the scaled coordinates, the two substitutions are
\[
\frac{\exp(l_i)-1}{\lambda}
=\sum_{n\ge1}\frac{\lambda^{n-1}}{n!}(\lambda^{-1}l_i)^n,
\qquad
\frac{\log(1+x_i)}{\lambda}
=\sum_{n\ge1}(-1)^{n+1}\frac{\lambda^{n-1}}{n}(\lambda^{-1}x_i)^n.
\]
Since both $v_\ell(n!)$ and $v_\ell(n)$ are at most $(n-1)/(\ell-1)$, the coefficients in degree $n$ have valuation at least $(n-1)(v_\ell(\lambda)-1/(\ell-1))$. This bound is nonnegative and tends to infinity. Both substitutions thus preserve the integral Tate algebras. They are inverse by Proposition~\ref{prop:log-coords}(1). Inverting $\ell$ gives the bounded $E$-algebra isomorphisms.
\end{proof}

Write $l_i=\log(q_i)$ and $l=(l_1,\dots,l_b)^{\mathrm t}$. If $\sigma(e_i)=\sum_j a_{ji}e_j$, then $\sigma^\#(q_i)=\prod_j q_j^{a_{ji}}$, so taking formal logarithms gives $\sigma^\#(l_i)=\sum_j a_{ji}l_j$. Thus $e_i\mapsto l_i$ identifies $V_\pi\otimes_{\mathbf Q_\ell}E$ with $\bigoplus_i El_i$ as a representation of $\sigma$. Let $E'/E$ be a finite extension containing all eigenvalues of $\sigma$. We also write $\sigma$ and $\Phi_w$ for their base changes to $X_{\pi,E'}$. Choose a block diagonal matrix $C\in\mathrm{GL}_b(E')$, with blocks indexed by $r\in W$, such that
\[
y=(y_1,\dots,y_b)^{\mathrm t}=Cl,
\qquad \sigma^\#(y_i)=\alpha_i y_i.
\]
Since $\sigma(\pi)=\pi$, each $\alpha_i$ is an $\ell$-adic unit. By Assumption~\ref{ass:frobenius-weights}(2) and the block diagonal choice of $C$, we have $|\iota(\alpha_i)|=c^{w_i}$.
Taking complex absolute values in $\alpha_1^{m_1}\cdots\alpha_b^{m_b}=1$ and using $c\ne1$ gives
\begin{equation}\label{eq:eigenvalue-weight-relation}
\prod_{i=1}^b \alpha_i^{m_i} = 1 \quad \Longrightarrow \quad \sum_{i=1}^b w_i m_i = 0 \qquad (m_1, \dots, m_b \in \mathbf{Z}).
\end{equation}

For the supremum norm on $(E')^b$, let $\|C^{-1}\|=\max_{i,j}|(C^{-1})_{ij}|$, with $\|C^{-1}\|=1$ if $b=0$. Choose $\lambda\in(E')^\times$ with $\rho=|\lambda|$ satisfying $\max(1,\|C^{-1}\|)\rho<\ell^{-1/(\ell-1)}$, and let
\[
A_\rho=E'\langle\lambda^{-1}y_1,\dots,\lambda^{-1}y_b\rangle,
\quad \mathfrak m_0=(y_1,\dots,y_b)\subset A_\rho,
\quad \mathbb B_{E'}^b(\rho)=\mathrm{Spm}(A_\rho),
\]
so that the $\mathfrak m_0$-adic completion of $A_\rho$ is $\widehat A_\rho=E'[[y_1,\dots,y_b]]$. The formulas $\sigma^\#(\lambda^{-1}y_i)=\alpha_i\lambda^{-1}y_i$ define an automorphism of $A_\rho$ and hence an automorphism $\sigma$ of $\mathbb B_{E'}^b(\rho)$.

\begin{proposition}\label{prop:diagonal-analytic-coords}
Under Assumption~\ref{ass:frobenius-weights}, let $D_w^\#$ be the continuous $E'$-algebra endomorphism of $E'[[y_1,\dots,y_b]]$ with $D_w^\#(y_i)=\ell^{w_i}y_i$.
\begin{enumerate}
\item[(1)] Composing $\mathrm{Exp}_{E'}^\#$ with $l=C^{-1}y$ gives a continuous $E'$-algebra isomorphism
\[
\mathrm{Exp}_C^\#:E'[[x_1,\dots,x_b]]\xrightarrow{\sim}E'[[y_1,\dots,y_b]],
\qquad x_i\mapsto\exp((C^{-1}y)_i)-1,
\]
satisfying
\begin{equation}\label{eq:expC-conjugacy}
\mathrm{Exp}_C^\# \circ \Phi_w^\# = D_w^\# \circ \mathrm{Exp}_C^\#.
\end{equation}
\item[(2)] There are a $\sigma$-stable affinoid neighbourhood $\mathcal U_\rho\subset X_{\pi,E'}$ of the identity and a $\sigma$-equivariant isomorphism of rigid analytic $E'$-spaces
\[
\mathrm{Exp}_{C,\rho}:\mathbb B_{E'}^b(\rho)\xrightarrow{\sim}\mathcal U_\rho,
\qquad y\mapsto(\exp((C^{-1}y)_i)-1)_i,
\]
with inverse $x\mapsto C(\log(1+x_1),\dots,\log(1+x_b))^{\mathrm t}$.
\end{enumerate}
\end{proposition}

\begin{proof}
(1) Composing the isomorphism $\mathrm{Exp}_{E'}^\#$ of Proposition~\ref{prop:log-coords}(1) with the substitution $l=C^{-1}y$ gives $\mathrm{Exp}_C^\#$. Since $C$ commutes with $\mathrm{diag}(\ell^{w_1},\dots,\ell^{w_b})$, Proposition~\ref{prop:log-coords}(2) gives~\eqref{eq:expC-conjugacy}.

(2) Applying Proposition~\ref{prop:analytic-log-coords} over $E'$ at radius $\|C^{-1}\|\rho<\ell^{-1/(\ell-1)}$ and substituting $l=C^{-1}y$ identifies $\mathbb B_{E'}^b(\rho)$ with the Weierstrass domain $\mathcal U_\rho$ defined in the closed polydisc $|x_i|\le\|C^{-1}\|\rho$ $(1\le i\le b)$ by
\[
\bigl|\bigl(C(\log(1+x_1),\dots,\log(1+x_b))^{\mathrm t}\bigr)_i\bigr|\le\rho
\qquad(1\le i\le b).
\]
These inequalities hold on a sufficiently small closed polydisc centred at the identity, so $\mathcal U_\rho$ is an affinoid neighbourhood of the identity. The identities $\sigma^\#(y_i)=\alpha_i y_i$ show that $\mathrm{Exp}_{C,\rho}$ is $\sigma$-equivariant, so $\sigma$ and $\sigma^{-1}$ preserve $\mathcal U_\rho$.
\end{proof}

\begin{lemma}\label{lem:tate-origin}
Under Assumption~\ref{ass:frobenius-weights}, let $I \subset A_\rho$ be a proper ideal. If $\sigma^\#(I) = I$, then $I \subset \mathfrak{m}_0$.
\end{lemma}

\begin{proof}
Since all weights $w_i \in W$ are strictly positive integers, the relation \eqref{eq:eigenvalue-weight-relation} implies that $\prod_{i=1}^b \alpha_i^{n_i} \ne 1$ for every nonzero tuple $n = (n_1, \dots, n_b) \in \mathbf{N}^b \setminus \{0\}$. By \cite[Proposition~2.1]{EK21}, every proper $\sigma^\#$-stable ideal of $A_\rho$ is contained in the maximal ideal $\mathfrak m_0$ of the origin.
\end{proof}

For any formal power series $f = \sum_{n \in \mathbf{N}^b} a_n y^n \in \widehat{A}_\rho = E'[[y_1, \dots, y_b]]$, write its expansion into weighted homogeneous components with respect to $w = (w_1, \dots, w_b)$ as
\[
f = \sum_{m \ge 0} f^{(m)}, \qquad f^{(m)} := \sum_{\sum_{i=1}^b w_i n_i=m} a_n y^n.
\]

\begin{definition}\label{def:w-homogeneous}
An ideal $\mathfrak{a} \subset \widehat{A}_\rho$ is called \emph{$w$-homogeneous} if $f \in \mathfrak{a}$ implies $f^{(m)} \in \mathfrak{a}$ for every integer $m \ge 0$.
\end{definition}

\begin{proposition}\label{prop:frob-weighted-homog}
Under Assumption~\ref{ass:frobenius-weights}, let $I\subset A_\rho$ be a proper ideal satisfying $\sigma^\#(I)=I$. Then $I\widehat A_\rho$ is $w$-homogeneous.
\end{proposition}

\begin{proof}
Let $\widehat I=I\widehat A_\rho$, which is proper by Lemma~\ref{lem:tate-origin}. Consider the cocharacter $\mu:\mathbf G_{m,E'}\to\mathbf G_{m,E'}^b$, $a\mapsto(a^{w_1},\dots,a^{w_b})$. The image of $\mu$ lies in the Zariski closure of the cyclic subgroup generated by $(\alpha_1,\dots,\alpha_b)$. Indeed, a character of the diagonal torus that is trivial on this cyclic subgroup has exponent vector $m$ satisfying $\prod_i\alpha_i^{m_i}=1$. By~\eqref{eq:eigenvalue-weight-relation}, $\sum_iw_i m_i=0$, so this character is also trivial on the image of $\mu$. The character correspondence for closed subgroups of a split torus \cite[Theorem~12.9]{Mil17} then shows that the image of $\mu$ lies in the Zariski closure of $\langle(\alpha_1,\dots,\alpha_b)\rangle$.

Fix $f\in\widehat I$ and an integer $m\ge0$. For each $N\ge1$, the stabilizer of the image of $\widehat I$ in $\widehat A_\rho/\mathfrak m_0^N\widehat A_\rho$ is Zariski closed in the diagonal torus and contains $(\alpha_1,\ldots,\alpha_b)$, hence also the image of $\mu$. For $n\in\mathbf N^b$, the monomial $y^n$ has weight $\sum_i w_i n_i$ under $\mu$. Thus projection onto the weight $m$ summand sends the class of $f$ to the class of $f^{(m)}$. Hence $f^{(m)}\in\widehat I+\mathfrak m_0^N\widehat A_\rho$ for every $N\ge1$. Since $\widehat I$ is closed in the $\mathfrak m_0$-adic topology, $f^{(m)}\in\widehat I$.
\end{proof}

\begin{lemma}\label{lem:w-homog-phi-stable}
Under Assumption~\ref{ass:frobenius-weights}, let $\mathfrak a\subset\widehat A_\rho$ be a $w$-homogeneous ideal. Then $(\mathrm{Exp}_C^\#)^{-1}(\sqrt{\mathfrak a})$ is preserved by $\Phi_w^\#$.
\end{lemma}

\begin{proof}
First, $\sqrt{\mathfrak a}$ is $w$-homogeneous. For $0\ne f\in\sqrt{\mathfrak a}$ choose $N\ge1$ with $f^N\in\mathfrak a$, and let $d$ be the least weight such that $f^{(d)}\ne0$. The component of $f^N$ of weight $Nd$ is $(f^{(d)})^N$. It is in $\mathfrak a$ by homogeneity, so $f^{(d)}\in\sqrt{\mathfrak a}$. Subtracting $f^{(d)}$ and proceeding by induction on the weight gives $f^{(m)}\in\sqrt{\mathfrak a}$ for every $m$.

For $f\in\sqrt{\mathfrak a}$ every partial sum of $D_w^\#(f)=\sum_{m\ge0}\ell^m f^{(m)}$ is thus in $\sqrt{\mathfrak a}$. The series converges $\mathfrak m_0$-adically and the ideal is closed, so $D_w^\#(f)\in\sqrt{\mathfrak a}$. Conjugating by $\mathrm{Exp}_C^\#$ via~\eqref{eq:expC-conjugacy} shows that $(\mathrm{Exp}_C^\#)^{-1}(\sqrt{\mathfrak a})$ is preserved by $\Phi_w^\#$.
\end{proof}

For a radical $\ell$-saturated ideal $J\subset R_\pi$, recall that $J^{(d)}$ is the intersection of the minimal primes $\mathfrak p$ of $J$ for which $\dim(R_\pi/\mathfrak p)=d$, as in Corollary~\ref{cor:fixed-dim-components}.

\begin{proposition}[Frobenius stability implies the rigid dynamical criterion]\label{prop:frob-gives-phi-criterion}
Under Assumption~\ref{ass:frobenius-weights}, let $J\subset R_\pi$ be a radical $\ell$-saturated ideal satisfying $\sigma^\#(J)=J$. For each $d\ge0$ with $Z_{J^{(d)}}\ne\varnothing$, there exist integers $a_d>b_d\ge0$ such that
\[
\Phi_w^{a_d}(Z_{J^{(d)}})\subset\Phi_w^{b_d}(Z_{J^{(d)}}).
\]
\end{proposition}

\begin{proof}
If $\mathrm{rank}(\pi)=0$, then $X_\pi$ is a point and we can take $a_d=1$ and $b_d=0$. Assume $\mathrm{rank}(\pi)\ge1$. Let $R'=\mathcal O_{E'}[[\pi]]$ and $J'=JR'$. For each minimal prime $\mathfrak p$ of $J$, the domain $R_\pi/\mathfrak p$ is flat over $\mathcal O_E$ because $J$ is $\ell$-saturated. By Lemma~\ref{lem:base-change-components}, $\mathfrak pR'$ is radical and every minimal prime over $\mathfrak pR'$ has quotient of dimension $\dim(R_\pi/\mathfrak p)$. Since $R_\pi\to R'$ is finite faithfully flat, these primes, as $\mathfrak p$ varies, are precisely the minimal primes of $J'$. Flat base change preserves finite intersections, so $J^{(d)}R'=(J')^{(d)}$. Lemma~\ref{lem:rigid-coeff-ext} now gives $(Z_{J^{(d)}})_{E'}=Z_{(J')^{(d)}}$.
For each integer $k\ge0$, Lemma~\ref{lem:rigid-images}(2) shows that the ideal defining $\Phi_w^k(Z_{J^{(d)}})$ is the kernel of the composite $R_\pi\xrightarrow{(\Phi_w^\#)^k}R_\pi\twoheadrightarrow R_\pi/J^{(d)}$. Since $\mathcal O_{E'}$ is flat over $\mathcal O_E$, extending this kernel to $R'$ gives the ideal defining $\Phi_w^k(Z_{(J')^{(d)}})$.

Since $\sigma^\#$ permutes the finitely many minimal primes of $J'$, choose $N\ge 1$ such that $(\sigma^\#)^N$ fixes each of them. In the coordinates $y=Cl$, $(\sigma^\#)^N(y_i)=\alpha_i^N y_i$ and $|\iota(\alpha_i^N)|=(c^N)^{w_i}$. Thus Assumption~\ref{ass:frobenius-weights} holds with $\sigma^N$ and $c^N$ in place of $\sigma$ and $c$, and with the same weights.

Fix a minimal prime ideal $\mathfrak p$ of $(J')^{(d)}$ in $R'$. By Corollary~\ref{cor:standard-point-ideals}(1), $Z_{\mathfrak p}\subset X_{\pi,E'}$ has a rigid point $z$ whose residue field $\kappa(z)$ is a finite extension of $E'$. For $1\le i\le b$, we have $q_i(z)\in1+\mathfrak m_{\kappa(z)}$ and $w_i>0$. Hence $q_i(\Phi_w^k(z))=q_i(z)^{\ell^{kw_i}}\to1$ as $k\to\infty$, so $\Phi_w^k(z)$ tends to the identity. Choose $M=M(\mathfrak p)\ge 0$ such that $\Phi_w^M(z)$ is in the affinoid neighbourhood $\mathcal U_\rho$ of Proposition~\ref{prop:diagonal-analytic-coords}(2).

Let $\mathfrak p_M:=((\Phi_w^\#)^M)^{-1}(\mathfrak p)$. By Lemma~\ref{lem:rigid-images}, $\mathfrak p_M$ is a prime $\ell$-saturated ideal of $R'$ defining $Z_{\mathfrak p_M}=\Phi_w^M(Z_{\mathfrak p})$. Since $\Phi_w$ commutes with $\sigma^N$ and $(\sigma^\#)^N(\mathfrak p)=\mathfrak p$, we have $(\sigma^\#)^N(\mathfrak p_M)=\mathfrak p_M$.

Let $I\subset A_\rho$ be the defining ideal of the closed affinoid subspace $\mathrm{Exp}_{C,\rho}^{-1}(Z_{\mathfrak p_M}\cap\mathcal U_\rho)$ of $\mathbb B_{E'}^b(\rho)$, with the induced closed subspace structure. The inverse image of $\Phi_w^M(z)$ under $\mathrm{Exp}_{C,\rho}$ is a zero of $I$, so $I$ is proper. Equivariance of $\mathrm{Exp}_{C,\rho}$ under $\sigma^N$ gives $(\sigma^\#)^N(I)=I$. Applying Lemma~\ref{lem:tate-origin} to $\sigma^N$ gives $I\subset\mathfrak m_0$. Since $\mathrm{Exp}_{C,\rho}$ sends the origin to the identity, $Z_{\mathfrak p_M}$ contains the identity. Hence $\mathfrak p_M\subset (x_1,\ldots,x_b)R'$.

By Lemma~\ref{lem:ext-contract-identity}, the extension of $\mathfrak p_M$ to $E'[[x_1,\ldots,x_b]]$ is radical and contracts to $\mathfrak p_M$ in $R'$. The restriction and completion maps give a commutative diagram
\[
\begin{tikzcd}[column sep=normal]
R' \arrow[r] \arrow[d,hook] & A_\rho \arrow[d] \\
E'[[x_1,\ldots,x_b]] \arrow[r,"\mathrm{Exp}_C^\#","\sim"'] & E'[[y_1,\ldots,y_b]]\mathrlap{{}=\widehat A_\rho.}
\end{tikzcd}
\]
The upper arrow is restriction to $\mathcal U_\rho$ followed by pullback along $\mathrm{Exp}_{C,\rho}$. The right vertical arrow is the completion map $A_\rho\to\widehat A_\rho$ at $\mathfrak m_0$. Since $I=\mathfrak p_M A_\rho$, we obtain
\begin{equation}\label{eq:exp-preimage-pM}
(\mathrm{Exp}_C^\#)^{-1}\bigl(\sqrt{I\widehat A_\rho}\bigr)
=\mathfrak p_M E'[[x_1,\ldots,x_b]].
\end{equation}
Proposition~\ref{prop:frob-weighted-homog} applied to $\sigma^N$ shows that $I\widehat A_\rho$ is $w$-homogeneous. By Lemma~\ref{lem:w-homog-phi-stable} and \eqref{eq:exp-preimage-pM}, $\Phi_w^\#$ preserves $\mathfrak p_M E'[[x_1,\ldots,x_b]]$. Contracting to $R'$ gives $\Phi_w^\#(\mathfrak p_M)\subset\mathfrak p_M$. By Lemmas~\ref{lem:faithfulness} and~\ref{lem:rigid-images}, this inclusion is equivalent to $\Phi_w(Z_{\mathfrak{p}_M}) \subset Z_{\mathfrak{p}_M}$, or
\[
\Phi_w^{M(\mathfrak{p})+1}(Z_{\mathfrak{p}}) \subset \Phi_w^{M(\mathfrak{p})}(Z_{\mathfrak{p}}).
\]

Let $\mathcal P_d$ be the set of minimal primes of $(J')^{(d)}$ in $R'$ and let $M_d=\max_{\mathfrak p\in\mathcal P_d} M(\mathfrak p)$. Applying $\Phi_w^{M_d-M(\mathfrak p)}$ to the last inclusion gives $\Phi_w^{M_d+1}(Z_{\mathfrak p})\subset\Phi_w^{M_d}(Z_{\mathfrak p})$ for every $\mathfrak p\in\mathcal P_d$. Since inverse images of ideals commute with finite intersections, Lemmas~\ref{lem:analytic-union-components} and~\ref{lem:rigid-images}(2) give
\[
\Phi_w^k(Z_{(J')^{(d)}})
=\bigcup_{\mathfrak p\in\mathcal P_d}{}^{\mathrm{an}}\Phi_w^k(Z_{\mathfrak p})
\qquad(k=M_d,M_d+1).
\]
Thus $\Phi_w^{M_d+1}(Z_{(J')^{(d)}})\subset\Phi_w^{M_d}(Z_{(J')^{(d)}})$. By faithful flatness of $E'/E$, we obtain $\Phi_w^{M_d+1}(Z_{J^{(d)}})\subset\Phi_w^{M_d}(Z_{J^{(d)}})$, as required.
\end{proof}

\begin{corollary}\label{cor:frob-stable-quasilinear}
Under Assumption~\ref{ass:frobenius-weights}, let $J\subset R_\pi$ be a radical $\ell$-saturated ideal satisfying $\sigma^\#(J)=J$. Then $Z_J\subset X_\pi$ is $\Phi_w$-stable rigid quasilinear.
\end{corollary}

\begin{proof}
Apply Corollary~\ref{cor:fixed-dim-components} to the inclusions supplied by Proposition~\ref{prop:frob-gives-phi-criterion}.
\end{proof}

\subsection{Arithmetic actions and weight filtrations}\label{subsec:arithmetic-weights}
Let $F$ be an algebraically closed field of characteristic different from $\ell$, and let $X/F$ be a separated, connected scheme of finite type. Let $\Gamma_\ell(X)$ be the maximal torsion free abelian pro-$\ell$ quotient of $\pi_1(X)$.

For a finitely generated subfield $F_0\subset F$, let $\overline F_0$ denote its algebraic closure in $F$. For a model $X_{F_0}$ of $X$ over $F_0$, write $X_{\overline F_0}=X_{F_0}\times_{F_0}\overline F_0$. For $h\in\mathrm{Gal}(\overline F_0/F_0)$, let $h^*$ denote pullback along the automorphism $\mathrm{id}_{X_{F_0}}\times\mathrm{Spec}(h)$ of $X_{\overline F_0}$.

\begin{definition}\label{def:arithmetic}
\leavevmode
\begin{enumerate}
\item[(1)] A constructible complex $\mathcal F\in D_c^b(X,\overline{\mathbf Q}_\ell)$ is \emph{arithmetic} if there exist a finitely generated subfield $F_0\subset F$, a model $X_{F_0}$, and a complex $\mathcal F_{\overline F_0}\in D_c^b(X_{\overline F_0},\overline{\mathbf Q}_\ell)$ whose pullback to $X$ is isomorphic to $\mathcal F$ and which satisfies $h^*\mathcal F_{\overline F_0}\simeq\mathcal F_{\overline F_0}$ for every $h\in \mathrm{Gal}(\overline F_0/F_0)$ \cite[Definition~5.1]{EK21}.
\item[(2)] A torsion free quotient $u:\Gamma_\ell(X)\twoheadrightarrow\pi$, with $\pi\simeq\mathbf Z_\ell^b$, is \emph{arithmetic} if there exist a finitely generated subfield $F_0\subset F$ and a model $X_{F_0}$ for which an open subgroup of $\mathrm{Gal}(\overline F_0/F_0)$ preserves $\ker(u)$ under the induced Galois action on $\Gamma_\ell(X)$.
\end{enumerate}
\end{definition}

Fix a model $X_{F_0}$ of $X$ and a torsion free quotient $u:\Gamma_\ell(X)\twoheadrightarrow\pi$. Identify $H^1(X_{\overline F_0},\mathbf Q_\ell)$ with the space of continuous $\mathbf Q_\ell$-valued homomorphisms on $\Gamma_\ell(X)$. The injection
\[
V_\pi^\vee=\mathrm{Hom}_{\mathbf Z_\ell}(\pi,\mathbf Q_\ell)\hookrightarrow H^1(X_{\overline F_0},\mathbf Q_\ell),
\qquad \lambda\mapsto\lambda\circ u,
\]
identifies $V_\pi^\vee$ with the annihilator of $\ker(u)$.

Since $\pi$ is torsion free, the functionals in $V_\pi^\vee$ separate its points. Hence $\ker(u)=\bigcap_{\lambda\in V_\pi^\vee}\ker(\lambda\circ u)$. A subgroup of $\mathrm{Gal}(\overline F_0/F_0)$ therefore preserves $\ker(u)$ if and only if it preserves $V_\pi^\vee$. In particular, $u$ is arithmetic relative to the model $X_{F_0}$ if and only if $V_\pi^\vee$ has an open stabilizer. If $u$ is arithmetic relative to $X_{F_0}$, passing to the finite separable extension of $F_0$ fixed by the stabilizer of $V_\pi^\vee$ makes the injection Galois equivariant.

\begin{remark}\label{rem:common-descent}
An arithmetic quotient $u$ and an arithmetic complex $\mathcal F$ admit a common model. Choose a finitely generated subfield $K\subset F$ containing the fields of definition of both models of $X$ and large enough that their identification over $F$ descends to $K$. Write $\mathcal F_{\overline K}$ for the pullback of the descended complex to $X_{\overline K}$. In $\mathrm{Gal}(\overline K/K)$, the stabilizers of $\ker(u)$ and of the isomorphism class of $\mathcal F_{\overline K}$ are open, since they contain the inverse images of the respective open stabilizers over the original fields. Let $F_0$ be the fixed field of their intersection. Base change to this finite separable extension of $K$ gives a model $X_{F_0}$ on which every $h\in \mathrm{Gal}(\overline F_0/F_0)$ satisfies $h^*\mathcal F_{\overline F_0}\simeq\mathcal F_{\overline F_0}$ and $h(\ker(u))=\ker(u)$.
\end{remark}

\begin{remark}\label{rem:finite-field-case}
Suppose $F=\overline{\mathbf F}_p$ with $p\ne\ell$ and choose a model over a finite subfield $\mathbf F_q$. The absolute Galois group $\mathrm{Gal}(F/\mathbf F_q)\simeq\widehat{\mathbf Z}$ is topologically generated by the geometric Frobenius automorphism $\mathrm{Frob}_q:F\to F$, $x\mapsto x^{1/q}$. The open subgroups of $\mathrm{Gal}(F/\mathbf F_q)$ are
\[
\mathrm{Gal}(F/\mathbf F_{q^n})=\overline{\langle\mathrm{Frob}_q^n\rangle},
\qquad n\ge1.
\]
Definition~\ref{def:arithmetic} therefore takes the following form:
\begin{enumerate}
\item[(1)] A complex $\mathcal F$ is arithmetic if and only if $(\mathrm{Frob}_q^n)^*\mathcal F\simeq\mathcal F$ for some $n\ge1$.
\item[(2)] A quotient $u$ is arithmetic if and only if $\mathrm{Frob}_q^n(\ker(u))=\ker(u)$ for some $n\ge1$.
\end{enumerate}
\end{remark}

We recall the definitions of pure and mixed Galois representations and their weight filtrations following \cite[Section~2]{Jan10}.

\begin{definition}\label{def:pure-rep}
Let $V$ be a finite dimensional continuous $\mathbf Q_\ell$-representation of $\mathrm{Gal}(\overline F_0/F_0)$. For $n\in\mathbf Z$, the representation $V$ is \emph{pure of weight $n$} if there is a normal integral scheme $S$ of finite type over $\mathbf Z[1/\ell]$, with function field $F_0$, satisfying the following conditions:
\begin{enumerate}
\item[(1)] The representation $V$ is unramified over $S$, meaning that the action factors through $\pi_1(S,\overline\eta)$, where $\overline\eta=\mathrm{Spec}(\overline F_0)$ is the geometric generic point of $S$.
\item[(2)] For every closed point $s\in S$ with residue field $\kappa(s)=\mathbf F_{q_s}$, let $\mathrm{Frob}_s$ denote the conjugacy class in $\pi_1(S,\overline\eta)$ induced by the geometric Frobenius automorphism $\overline{\kappa(s)}\to\overline{\kappa(s)}$, $x\mapsto x^{1/q_s}$. Every eigenvalue $\alpha$ of $\mathrm{Frob}_s$ on $V$ is algebraic over $\mathbf Q$ and satisfies $|\tau(\alpha)|=q_s^{n/2}$ for every embedding $\tau:\mathbf Q(\alpha)\hookrightarrow\mathbf C$.
\end{enumerate}
\end{definition}

\begin{definition}\label{def:mixed-rep}
Let $V$ be a finite dimensional continuous $\mathbf Q_\ell$-representation of $\mathrm{Gal}(\overline F_0/F_0)$. We call $V$ \emph{mixed} if it admits a finite increasing $\mathrm{Gal}(\overline F_0/F_0)$-stable filtration $W'_\bullet V$ such that every nonzero graded piece $\mathrm{Gr}^{W'}_a V=W'_aV/W'_{a-1}V$ is pure of an integer weight $n_a$. We call $V$ \emph{mixed with all weights strictly positive} if the filtration can be chosen with $n_a>0$ for every nonzero graded piece.
\end{definition}

\begin{definition}\label{def:weight-filtration}
A \emph{weight filtration} on a finite dimensional continuous $\mathbf Q_\ell$-representation $V$ of $\mathrm{Gal}(\overline F_0/F_0)$ is a finite increasing $\mathrm{Gal}(\overline F_0/F_0)$-stable filtration $W_\bullet V$ for which every nonzero graded piece $\mathrm{Gr}^W_n V$ is pure of weight $n$. A weight filtration is unique when it exists.
\end{definition}

Suppose that $u$ is arithmetic. Choose a field of definition $F_0$ such that $\mathrm{Gal}(\overline F_0/F_0)$ preserves $\ker(u)$. Then $V_\pi^\vee$ is a Galois-stable subspace of $H^1(X_{\overline F_0},\mathbf Q_\ell)$, and its induced weight filtration is given by
\[
W_n V_\pi^\vee = V_\pi^\vee \cap W_n H^1(X_{\overline{F}_0}, \mathbf{Q}_\ell).
\]
For this intersection formula to give a weight filtration on a subspace of $H^1(X_{\overline F_0},\mathbf Q_\ell)$, the subspace must be Galois-stable.

\begin{lemma}[A geometric Frobenius lift and semisimplicity]\label{lem:frob-lift}
Suppose that $u:\Gamma_\ell(X)\twoheadrightarrow\pi$ is arithmetic, that $V_\pi^\vee$ is mixed with all weights strictly positive, and that $\mathcal F$ is arithmetic. For some common field of definition $F_0$ as in Remark~\ref{rem:common-descent}, there exist
\begin{enumerate}
\item[(1)] a normal integral scheme $S$ with function field $F_0$, of finite type over a finite field in positive characteristic and over $\mathbf Z[1/\ell]$ in characteristic zero,
\item[(2)] a closed point $s\in S$ with residue field $\kappa(s)=\mathbf F_{q_s}$,
\item[(3)] an integer $N>0$, and
\item[(4)] an element $\sigma\in\mathrm{Gal}(\overline F_0/F_0)$ lifting a representative of the conjugacy class of $\mathrm{Frob}_s^N$ in $\pi_1(S,\overline\eta)$.
\end{enumerate}
The lift $\sigma$ satisfies $\sigma^*\mathcal F_{\overline F_0}\simeq\mathcal F_{\overline F_0}$ and $\sigma(\ker(u))=\ker(u)$, and acts semisimply on $V_\pi$. Moreover, $V_\pi$ has a decomposition $V_\pi=\bigoplus_{r\in W}V_r$ into $\sigma$-stable $\mathbf Q_\ell$-subspaces indexed by $W\subset\mathbf Z_{>0}$. Every eigenvalue $\alpha$ of $\sigma$ on $V_r$ satisfies $|\iota(\alpha)|=q_s^{-Nr/2}$.
\end{lemma}

\begin{proof}
Choose $F_0$ and $X_{F_0}$ as in Remark~\ref{rem:common-descent}, so that every element of $\mathrm{Gal}(\overline F_0/F_0)$ preserves $\ker(u)$ and fixes the isomorphism class of $\mathcal F_{\overline F_0}$. If $\mathrm{char}(F)=0$ and $\ell=2$, replace $F_0$ by the finite extension $F_0(\sqrt{-1})\subset F$ so that it has no real places. Choose a normal affine integral scheme $S$ with function field $F_0$ over which $X_{F_0}$ extends to a separated morphism of finite type $f:\mathscr X\to S$. Take $S$ of finite type over a finite field in positive characteristic and over $\mathbf Z[1/\ell]$ in characteristic zero (over $\mathbf Z[1/2,\sqrt{-1}]$ when $\mathrm{char}(F)=0$ and $\ell=2$).

The scheme $S$ satisfies condition (LO) of \cite[Section~2.1]{Sun17}. The constant system $\mathbf Z_\ell$ is stratifiable since its reductions have constant simple constituents. By constructibility, after shrinking $S$ to a dense affine open, $R^1f_*\mathbf Q_\ell$ is lisse. After a further shrinking, generic base change \cite[Theorem~2.5(i)]{Sun17} gives $(R^1f_*\mathbf Q_\ell)_{\bar t}\simeq H^1(\mathscr X_{\bar t},\mathbf Q_\ell)$ for every geometric point $\bar t\to S$ with algebraically closed residue field, since the spectrum of that field satisfies (LO). Choose a filtration $W'_\bullet V_\pi^\vee$ with strictly positive weights as in Definition~\ref{def:mixed-rep}. The subrepresentation $V_\pi^\vee\subset H^1(X_{\overline F_0},\mathbf Q_\ell)$ and the steps of this filtration define lisse subsheaves of $R^1 f_*\mathbf Q_\ell$. By Definition~\ref{def:pure-rep}, after a further shrinking, each nonzero graded piece of $W'_\bullet$ is pure of its specified weight at every closed point of $S$.

Choose a closed point $s\in S$ with residue field $\kappa(s)=\mathbf F_{q_s}$ and a geometric point $\overline s$ above $s$. An étale path from $\overline\eta$ to $\overline s$ determines a representative of $\mathrm{Frob}_s$ in $\pi_1(S,\overline\eta)$. Choose a lift $\phi_s\in\mathrm{Gal}(\overline F_0/F_0)$ of this representative. Write $H^1=H^1(\mathscr X_{\overline s},\mathbf Q_\ell)$. Transport along the path identifies the action of $\phi_s$ on $H^1(X_{\overline F_0},\mathbf Q_\ell)$ with the action of $\mathrm{Frob}_s$ on $H^1$.

The Picard $1$-motive of $\mathscr X_s$ is $[L_s\to G_s]$, where $L_s$ is a finite free abelian group and $G_s$ is an extension of an abelian variety $A_s$ by a torus $T_s$. By \cite[Theorems~1.9 and~6.22, Proposition~5.8, and Proposition~6.24]{Man13}, its realization gives a $\mathrm{Frob}_s$-stable filtration with graded pieces
\[
\mathrm{Gr}^W_r H^1(\mathscr X_{\overline s},\mathbf Q_\ell)\simeq
\begin{cases}
V_\ell(T_s)(-1) & r=0,\\
V_\ell(A_s)(-1) & r=1,\\
(L_s\otimes_{\mathbf Z}\mathbf Q_\ell)(-1) & r=2
\end{cases},
\]
where $V_\ell(H)$ denotes the rational $\ell$-adic Tate module of $H$; see \cite[Section~5.7]{Man13}.

The action of $\mathrm{Frob}_s$ on each of $X^*(T_s)$ and $L_s$ has finite order. Choose $N>0$ so that $\mathrm{Frob}_s^N$ acts trivially on both lattices, and let $\sigma=\phi_s^N$. In weight $0$, the twist in $V_\ell(T_s)\simeq\mathrm{Hom}_{\mathbf Z}(X^*(T_s),\mathbf Q_\ell(1))$ cancels the outer twist by $(-1)$, so $\sigma$ acts as the identity. In weight $2$ the twist by $(-1)$ gives multiplication by $q_s^N$. In weight $1$, $\mathrm{Frob}_s$ acts semisimply by \cite[Section~2, p.~138]{Tat66}.

For $r=1,2$, let $P_r(T)\in\mathbf Q_\ell[T]$ be the minimal polynomial of $\sigma$ on $\mathrm{Gr}^W_r H^1$, with $P_r=1$ if this graded piece is zero. Both polynomials are squarefree. The roots of $P_r$ have absolute value $q_s^{Nr/2}$ under $\iota$, so $P_1$ and $P_2$ are coprime. Since $P_2(\sigma)H^1\subset W_1H^1$ and $P_1(\sigma)W_1H^1\subset W_0H^1$, the squarefree polynomial $P_1P_2$ annihilates the action of $\sigma$ on $H^1/W_0H^1$. Thus $\sigma$ acts semisimply on this quotient.

On $W_0H^1$, $\sigma$ acts as the identity. On each nonzero graded piece of $W'_\bullet V_\pi^\vee$ of weight $r>0$, every eigenvalue of $\sigma$ has absolute value $q_s^{Nr/2}>1$ under $\iota$. Hence $V_\pi^\vee\cap W_0H^1=0$, and the projection to $H^1/W_0H^1$ embeds $V_\pi^\vee$ as a $\sigma$-stable subspace. In particular, $\sigma$ acts semisimply on $V_\pi^\vee$.

For $r=1,2$, let $K_r=\ker(P_r(\sigma):V_\pi^\vee\to V_\pi^\vee)$. Since $P_1$ and $P_2$ are coprime and their product annihilates $V_\pi^\vee$, we have $V_\pi^\vee=K_1\oplus K_2$. Let $W=\{r\in\{1,2\}\mid K_r\ne0\}$. For $r\in W$, the dual of the projection onto $K_r$ identifies $V_r:=K_r^\vee$ with a subspace of $V_\pi$, giving $V_\pi=\bigoplus_{r\in W}V_r$. The action on each $V_r$ is semisimple, and its eigenvalues are the inverses of those on the weight $r$ summand of $V_\pi^\vee$, so they have absolute value $q_s^{-Nr/2}$ under $\iota$.
\end{proof}

\begin{proposition}\label{prop:weight-adapted-torus}
Suppose that $u:\Gamma_\ell(X)\twoheadrightarrow\pi$ is arithmetic and that $V_\pi^\vee$ is mixed with all weights strictly positive. Fix a lift $\sigma$ of a power of geometric Frobenius and a decomposition $V_\pi=\bigoplus_{r\in W}V_r$ as in Lemma~\ref{lem:frob-lift}. There exists a $\sigma$-stable lattice $\pi'\subset V_\pi$ with the following properties:
\begin{enumerate}
\item[(1)] The lattice $\pi'$ satisfies
\[
\pi\subset\pi',\qquad
[\pi':\pi]<\infty,\qquad
\pi'=\bigoplus_{r\in W}(\pi'\cap V_r).
\]
\item[(2)] The inclusion $\pi\hookrightarrow\pi'$ induces a finite étale $\sigma$-equivariant isogeny
\[
f:X_{\pi'}\longrightarrow X_\pi,\qquad
\chi\longmapsto\chi|_\pi.
\]
\item[(3)] Multiplication by $\ell^r$ on each summand $\pi'\cap V_r$ induces the mixed power endomorphism $\Phi_w:\mathfrak X_{\pi'}\to\mathfrak X_{\pi'}$, with pullback
\[
\Phi_w^\#:R_{\pi'}\longrightarrow R_{\pi'},\qquad
[\gamma]\longmapsto[\gamma]^{\ell^r}
\quad(\gamma\in\pi'\cap V_r).
\]
\end{enumerate}
\end{proposition}

\begin{proof}
For (1), each $\pi\cap V_r$ is a full $\sigma$-stable lattice in $V_r$. Their direct sum has finite index in $\pi$, so we can choose $M\ge0$ with $\ell^M\pi\subset\bigoplus_r(\pi\cap V_r)$. Let $\pi'=\bigoplus_r\ell^{-M}(\pi\cap V_r)$. Then $\pi'$ contains $\pi$ with finite index and $\pi'\cap V_r=\ell^{-M}(\pi\cap V_r)$. Each summand is $\sigma$-stable.

For (2), Lemma~\ref{lem:isogeny-properties}(1) applied to $\pi\subset\pi'$ gives the finite étale isogeny. Restriction of characters commutes with $\sigma$ since $\sigma$ preserves both lattices. For (3), multiplication by $\ell^r$ preserves the summand $\pi'\cap V_r$. On the completed group ring, $\Phi_w^\#([\gamma])=[\ell^r\gamma]=[\gamma]^{\ell^r}$ for $\gamma\in\pi'\cap V_r$.
\end{proof}

\section{Arithmetic cohomology jump loci}\label{sec:jump-loci}

Let $F$ be an algebraically closed field of characteristic different from $\ell$ and let $X/F$ be a separated, connected scheme of finite type. Fix a quotient $u:\Gamma_\ell(X)\twoheadrightarrow\pi\simeq\mathbf Z_\ell^b$. Let $E/\mathbf Q_\ell$ be a finite extension in $\overline{\mathbf Q}_\ell$. Write $V_\pi=\pi\otimes_{\mathbf Z_\ell}\mathbf Q_\ell$, $R_\pi=\mathcal O_E[[\pi]]$, and $X_\pi=(\mathrm{Spf}(R_\pi))_\eta$.

\subsection{Fourier--Mellin transforms}\label{subsec:fm-complexes}
Let $(A,\mathfrak m_A)$ be a complete Noetherian local $\mathcal O_E$-algebra with finite residue field, and let $\varrho:\pi_1(X)\to A^\times$ be a continuous character. For every $n\ge1$, let $\mathcal L_{A/\mathfrak m_A^n}$ be the rank one étale local system of $A/\mathfrak m_A^n$-modules on $X$ defined by the composite
\[
\pi_1(X)\xrightarrow{\varrho}A^\times\longrightarrow(A/\mathfrak m_A^n)^\times.
\]
We denote the inverse system $(\mathcal L_{A/\mathfrak m_A^n})_{n\ge1}$ by $\mathcal L_A$.

\begin{definition}\label{def:m-adic-cohom}
Following \cite[Section~4.1]{EK21}, for $\mathcal{F}_{\mathcal{O}_E}\in D_c^b(X,\mathcal O_E)$ we define the \emph{$\mathfrak m_A$-adic étale cohomology complexes} by
\[
R\Gamma_{(c)}(X, \mathcal{F}_{\mathcal{O}_E} \otimes_{\mathcal{O}_E}^{\mathbf{L}} \mathcal{L}_A) := R\varprojlim_n R\Gamma_{(c)}(X, \mathcal{F}_{\mathcal{O}_E} \otimes_{\mathcal{O}_E}^{\mathbf{L}} \mathcal{L}_{A/\mathfrak{m}_A^n}).
\]
In formulas with $(c)$, we omit $c$ for ordinary étale cohomology $R\Gamma$ and retain it for étale cohomology with compact support $R\Gamma_c$.
\end{definition}

\begin{proposition}\label{prop:fm-finiteness-base-change}
Suppose that $A$ has finite global dimension and let $\mathcal F_{\mathcal O_E}\in D_c^b(X,\mathcal O_E)$.
\begin{enumerate}
\item[(1)] The complex $R\Gamma_{(c)}(X,\mathcal F_{\mathcal O_E}\otimes_{\mathcal O_E}^{\mathbf L}\mathcal L_A)$ is in $D^b_{\mathrm{coh}}(A)$.
\item[(2)] For every quotient $A\twoheadrightarrow A'$ with its quotient topology, let $\mathcal L_{A'}$ be the local system defined by the composite $\pi_1(X)\xrightarrow{\varrho}A^\times\to(A')^\times$. The base change map
\[
R\Gamma_{(c)}(X,\mathcal F_{\mathcal O_E}\otimes_{\mathcal O_E}^{\mathbf L}\mathcal L_A)\otimes_A^{\mathbf L}A'
\longrightarrow R\Gamma_{(c)}(X,\mathcal F_{\mathcal O_E}\otimes_{\mathcal O_E}^{\mathbf L}\mathcal L_{A'})
\]
is an isomorphism.
\end{enumerate}
\end{proposition}

\begin{proof}
For ordinary cohomology, both assertions follow from \cite[Proposition~4.1(1) and~(2)]{EK21}.

For compactly supported cohomology, let $A_n=A/\mathfrak m_A^n$, $K_n=\mathcal F_{\mathcal O_E}\otimes^{\mathbf L}_{\mathcal O_E}\mathcal L_{A_n}$, and $C_n=R\Gamma_c(X,K_n)$. Since $\mathcal L_{A_{n+1}}$ is locally free of rank one and reduces to $\mathcal L_{A_n}$, change of coefficients gives $C_{n+1}\otimes^{\mathbf L}_{A_{n+1}}A_n\simeq C_n$ \cite[Tag~0GL6]{Stacks}. The complex $C_1$ has bounded finite-dimensional cohomology over the field $A_1$ \cite[Tags~0G2A and~0GLH]{Stacks}, so it is perfect. Since $A=\varprojlim_n A_n$ and the transition maps $A_{n+1}\to A_n$ are surjective with nilpotent kernels, $C=R\varprojlim_n C_n$ is perfect over $A$ and $C\otimes_A^{\mathbf L}A_n\simeq C_n$ \cite[Tag~0CQG]{Stacks}. This proves (1).

For $A'=A/I$, let $B_n=A/(I+\mathfrak m_A^n)$. Change of coefficients identifies $C\otimes_A^{\mathbf L}B_n$ with $R\Gamma_c(X,\mathcal F_{\mathcal O_E}\otimes^{\mathbf L}_{\mathcal O_E}\mathcal L_{B_n})$. Choose a bounded finite free complex $P^\bullet$ representing $C$. The transition maps of $P^\bullet\otimes_A B_n$ are degreewise surjective. Since $A'=\varprojlim_n B_n$, their derived inverse limit is computed termwise and equals $P^\bullet\otimes_A A'$. As $B_n=A'/\mathfrak m_{A'}^n$, Definition~\ref{def:m-adic-cohom} gives the base change isomorphism in (2).
\end{proof}

We now take $A=R_\pi$. The \emph{universal character} $\varrho_{\mathrm{univ}}$ is the composite
\[
\pi_1(X)\to\Gamma_\ell(X)\xrightarrow{u}\pi\xrightarrow{\gamma\mapsto[\gamma]}R_\pi^\times.
\]
Let $\mathcal L_\pi$ be the rank one $R_\pi$-local system on $X$ defined by $\varrho_{\mathrm{univ}}$. Since $R_\pi\simeq\mathcal O_E[[x_1,\ldots,x_b]]$ is regular local of global dimension $b+1$, the cohomology complexes in Proposition~\ref{prop:fm-finiteness-base-change}(1) admit bounded representatives consisting of finite free $R_\pi$-modules and hence are perfect \cite[Tags~00OC and~066Z]{Stacks}.

\Needspace{4\baselineskip}
\begin{definition}\label{def:fm-transform}
\leavevmode
\begin{enumerate}
\item[(1)] The \emph{integral Fourier--Mellin transform} $\mathrm{FM}_\pi^{(c)}(X,-):D_c^b(X,\mathcal O_E)\to D_{\mathrm{coh}}^b(R_\pi)$ is given by
\[
\mathrm{FM}_\pi^{(c)}(X, \mathcal{F}_{\mathcal{O}_E}) := R\Gamma_{(c)}(X, \mathcal{F}_{\mathcal{O}_E} \otimes_{\mathcal{O}_E}^{\mathbf{L}} \mathcal{L}_\pi).
\]
\item[(2)] For $\mathcal F\in D_c^b(X,E)$, choose an integral lattice $\mathcal F_{\mathcal O_E}\in D_c^b(X,\mathcal O_E)$ together with an isomorphism $\mathcal F_{\mathcal O_E}[1/\ell]\simeq\mathcal F$ \cite[Section~1.1.3]{Del80}. The \emph{rational Fourier--Mellin transform} $\mathrm{FM}_\pi^{(c)}(X,-):D_c^b(X,E)\to D_{\mathrm{coh}}^b(R_\pi[1/\ell])$ is given by
\[
\mathrm{FM}_\pi^{(c)}(X, \mathcal{F}) := \mathrm{FM}_\pi^{(c)}(X, \mathcal{F}_{\mathcal{O}_E}) \otimes_{R_\pi} R_\pi[1/\ell].
\]
\end{enumerate}
\end{definition}

\Needspace{4\baselineskip}
\begin{remark}\label{rem:fm-specialization}
Let $\mathcal F\in D_c^b(X,E)$.
\begin{enumerate}
\item[(1)] The rational complex $\mathrm{FM}_\pi^{(c)}(X,\mathcal F)$ is independent, up to isomorphism, of the choice of integral lattice. Given two $\mathcal O_E$-lattices for $\mathcal F$, choose a third lattice mapping to both with cones annihilated by a power of $\ell$. Applying $\mathrm{FM}_\pi^{(c)}(X,-)$ to the two lattice maps gives morphisms whose cones are annihilated by the same power of $\ell$. Hence both morphisms become isomorphisms after inverting $\ell$.
\item[(2)] Let $\chi\in X_\pi$ be a rigid point with residue field $\kappa(\chi)$, a finite extension of $E$. Let $\mathcal L_\chi$ be the rank one $\kappa(\chi)$-local system on $X$ defined by the composite $\pi_1(X)\to\pi\xrightarrow{\chi}\kappa(\chi)^\times$. Specialization at $\chi$ gives a canonical isomorphism
\begin{equation}\label{eq:fm-specialization}
\mathrm{FM}_\pi^{(c)}(X,\mathcal F)\otimes_{R_\pi[1/\ell]}^{\mathbf L}\kappa(\chi)\xrightarrow{\sim}
R\Gamma_{(c)}(X,\mathcal F\otimes_E^{\mathbf L}\mathcal L_\chi).
\end{equation}
Indeed, the evaluation homomorphism $\mathrm{ev}_\chi:R_\pi\to\kappa(\chi)$ satisfies $\mathrm{ev}_\chi([\gamma])=\chi(\gamma)$, so extension of scalars along $\mathrm{ev}_\chi$ carries $\mathcal L_\pi$ to $\mathcal L_\chi$. The quotient $R_\pi/\ker(\mathrm{ev}_\chi)$ is finite over $\mathcal O_E$, and $(R_\pi/\ker(\mathrm{ev}_\chi))[1/\ell]\simeq\kappa(\chi)$. Thus \eqref{eq:fm-specialization} follows by applying Proposition~\ref{prop:fm-finiteness-base-change}(2) to this quotient and inverting $\ell$.
\end{enumerate}
\end{remark}

\begin{remark}\label{rem:alg-closed-invariance}
Let $p:X\to X_{\overline F_0}$ be the projection for a common model as in Remark~\ref{rem:common-descent}.
\begin{enumerate}
\item[(1)] The projection $p$ induces an isomorphism
\[
\Gamma_\ell(X)\xrightarrow{\sim}\Gamma_\ell(X_{\overline F_0})
\]
on maximal torsion free abelian pro-$\ell$ quotients. Indeed, for a connected scheme $Y$ and $n\ge1$, one has
\[
H^1(Y,\mathbf Z/\ell^n)\simeq\mathrm{Hom}_{\mathrm{cont}}(\pi_1(Y),\mathbf Z/\ell^n).
\]
Invariance of these cohomology groups under extension of algebraically closed fields \cite[Arcata, V, Corollary~3.3]{Del77} identifies the maximal abelian pro-$\ell$ quotients by Pontryagin duality, and hence their maximal torsion free quotients. Thus $u$ also defines a quotient $\Gamma_\ell(X_{\overline F_0})\twoheadrightarrow\pi$.
\item[(2)] After a finite extension of $E$ if necessary, choose an $E$-model $\mathcal G$ of $\mathcal F_{\overline F_0}$ \cite[Section~1.1.3]{Del80}. Then $p^*\mathcal G$ is an $E$-model of $\mathcal F$, and $p$ induces an isomorphism
\[
\mathrm{FM}_\pi^{(c)}(X_{\overline F_0},\mathcal G)
\xrightarrow{\sim}
\mathrm{FM}_\pi^{(c)}(X,p^*\mathcal G).
\]
Indeed, choose an integral lattice $K$ in $\mathcal G$. By Remark~\ref{rem:fm-specialization}(1), we can compute the two rational transforms using $K$ and $p^*K$, respectively. Let $\mathfrak m_{R_\pi}$ be the maximal ideal of $R_\pi$. For each $n\ge1$, let $K_n=K\otimes_{\mathcal O_E}^{\mathbf L}\mathcal L_{R_\pi/\mathfrak m_{R_\pi}^n}$ on $X_{\overline F_0}$. By (1), the universal character on $X$ is the pullback of the universal character on $X_{\overline F_0}$, so $p^*K_n$ is the corresponding twist of $p^*K$. Invariance under extension of algebraically closed fields gives compatible isomorphisms \cite[Arcata, IV, Theorem~5.4 and V, Theorem~3.2 and Corollary~3.3]{Del77}
\[
R\Gamma_{(c)}(X_{\overline F_0},K_n)\xrightarrow{\sim}R\Gamma_{(c)}(X,p^*K_n).
\]
Applying $R\varprojlim_n$ gives the isomorphism of integral transforms by Definitions~\ref{def:m-adic-cohom} and~\ref{def:fm-transform}(1). Inverting $\ell$ then gives the asserted isomorphism of rational transforms by Definition~\ref{def:fm-transform}(2).
\end{enumerate}
\end{remark}

\subsection{Standard models of jump loci}\label{subsec:integral-eq-equiv}
Fix $\mathcal F\in D_c^b(X,E)$, $i\in\mathbf Z$, and $m\ge1$.

\begin{definition}\label{def:fm-loci}
\leavevmode
\begin{enumerate}
\item[(1)] The \emph{cohomology jump loci} are
\begin{equation}\label{eq:jump-locus-def}
\Sigma_m^{i,(c)}(X, \mathcal{F}) = \bigl\{\chi \in X_\pi \;\big|\; \dim_{\kappa(\chi)} H_{(c)}^i(X, \mathcal{F} \otimes_E^{\mathbf{L}} \mathcal{L}_\chi) \ge m\bigr\}.
\end{equation}
\item[(2)] Let $M^i := H^i(\mathrm{FM}_\pi(X, \mathcal{F}))$ and $\mathfrak p_\chi:=\ker(R_\pi[1/\ell]\to\kappa(\chi))$ for $\chi\in X_\pi$. The \emph{nonflat locus} is
\begin{equation}\label{eq:nf-locus-def}
\mathrm{NF}^i(X, \mathcal{F}) = \bigl\{\chi \in X_\pi \;\big|\; (M^i)_{\mathfrak{p}_\chi} \text{ is not flat over } (R_\pi[1/\ell])_{\mathfrak{p}_\chi}\bigr\}.
\end{equation}
\end{enumerate}
\end{definition}

Fix a rigid point $\chi\in X_\pi$ with residue field $\kappa(\chi)$. By \eqref{eq:fm-specialization}, $\chi\in\Sigma_m^{i,(c)}(X,\mathcal F)$ if and only if $\dim_{\kappa(\chi)}H^i\!\left(\mathrm{FM}_\pi^{(c)}(X,\mathcal F)\otimes^{\mathbf L}_{R_\pi[1/\ell]}\kappa(\chi)\right)\ge m$.

Flatness can be checked on any affinoid chart $V_{\pi,N}=\mathrm{Spm}(A_N)$ containing $\chi$. Indeed, if $\mathfrak m_\chi\subset A_N$ represents $\chi$, the local map $(R_\pi[1/\ell])_{\mathfrak p_\chi}\to(A_N)_{\mathfrak m_\chi}$ is faithfully flat by Lemma~\ref{lem:chart-flatness}. Thus $(M^i)_{\mathfrak p_\chi}$ is flat if and only if $(M^i\otimes_{R_\pi[1/\ell]} A_N)_{\mathfrak m_\chi}$ is flat. The flat locus of $M^i$ in $\mathrm{Spec}(R_\pi[1/\ell])$ is open because $M^i$ is finite over the Noetherian ring $R_\pi[1/\ell]$ \cite[Tag~00RC]{Stacks}.

For a bounded complex $P^\bullet$ of finite free $R_\pi$-modules, let $r_i=\operatorname{rank}_{R_\pi}P^i$ and write $d^j:P^j\to P^{j+1}$ for the differential in each degree $j\in\mathbf Z$. For a matrix $D$ over $R_\pi$, let $I_t(D)$ be the ideal of its $t\times t$ minors, with $I_t(D)=R_\pi$ for $t\le0$ and $I_t(D)=0$ when $t$ exceeds its number of rows or columns. Following \cite[Definition--Proposition~2.2]{BW15DG}, define
\[
D_m^i(P^\bullet) := I_{r_i - m + 1}\begin{pmatrix} d^{i-1} & 0 \\ 0 & d^i \end{pmatrix} \subset R_\pi.
\]

Choose an integral lattice $\mathcal F_{\mathcal O_E}$ for $\mathcal F$. Let $P^\bullet$ and $P_c^\bullet$ be bounded finite free representatives of $\mathrm{FM}_\pi(X,\mathcal F_{\mathcal O_E})$ and $\mathrm{FM}_\pi^c(X,\mathcal F_{\mathcal O_E})$, respectively.

\begin{definition}\label{def:fm-locus-ideals}
\leavevmode
\begin{enumerate}
\item[(1)] The \emph{reduced integral cohomology jump ideals} are
\[
J_m^i:=\sqrt{D_m^i(P^\bullet)R_\pi[1/\ell]}\cap R_\pi,
\qquad
J_m^{i,c}:=\sqrt{D_m^i(P_c^\bullet)R_\pi[1/\ell]}\cap R_\pi.
\]
\item[(2)] For the module $M^i$ of Definition~\ref{def:fm-loci}(2), let $\mathcal N_i$ be the set of primes $\mathfrak p\in\operatorname{Spec}(R_\pi[1/\ell])$ for which $(M^i)_{\mathfrak p}$ is not flat over $(R_\pi[1/\ell])_{\mathfrak p}$. The reduced integral ideal of the nonflat locus is
\[
 \mathfrak n_i:=\bigcap_{\mathfrak p\in\mathcal N_i}(\mathfrak p\cap R_\pi).
\]
If $\mathcal N_i$ is empty, the intersection is $R_\pi$.
\end{enumerate}
\end{definition}

\Needspace{4\baselineskip}
\begin{proposition}\label{prop:standard-fm-loci}
Fix $i\in\mathbf Z$ and $m\ge1$.
\begin{enumerate}
\item[(1)] The ideals $J_m^i$, $J_m^{i,c}$ and $\mathfrak n_i$ depend only on $\mathcal F$ and $u$, and are radical and $\ell$-saturated.
\item[(2)] The loci $\Sigma_m^i(X,\mathcal F)$, $\Sigma_m^{i,c}(X,\mathcal F)$, and $\mathrm{NF}^i(X,\mathcal F)$ are Zariski closed in $X_\pi$ and, endowed with their reduced closed rigid-subspace structures, satisfy
\[
\Sigma_m^i(X,\mathcal F)_{\mathrm{red}}=Z_{J_m^i},\qquad
\Sigma_m^{i,c}(X,\mathcal F)_{\mathrm{red}}=Z_{J_m^{i,c}},\qquad
\mathrm{NF}^i(X,\mathcal F)_{\mathrm{red}}=Z_{\mathfrak n_i}.
\]
\end{enumerate}
\end{proposition}

\begin{proof}
For $\mathfrak p\in\mathrm{Spec}(R_\pi[1/\ell])$, write $d^j_{\mathfrak p}=d^j\otimes_{R_\pi}\kappa(\mathfrak p)$. Since $\operatorname{im}(d^{i-1}_{\mathfrak p})\subset\ker(d^i_{\mathfrak p})$, we have
\[
\dim_{\kappa(\mathfrak p)}H^i(P^\bullet\otimes_{R_\pi}\kappa(\mathfrak p))
=r_i-\operatorname{rank}(d^{i-1}_{\mathfrak p})-\operatorname{rank}(d^i_{\mathfrak p}).
\]
The sum of the two ranks is the rank of the block diagonal matrix defining $D_m^i(P^\bullet)$. Thus the fiber has cohomology of dimension at least $m$ if and only if $\mathfrak p$ contains $D_m^i(P^\bullet)R_\pi[1/\ell]$.

For (compactly supported) cohomology, the fibre dimensions depend only on $\mathrm{FM}_\pi^{(c)}(X,\mathcal F)$. Hence $J_m^i$ and $J_m^{i,c}$ are independent of the chosen finite free representatives and, by Remark~\ref{rem:fm-specialization}(1), of the integral lattice. The nonflat locus depends only on $M^i=H^i(\mathrm{FM}_\pi(X,\mathcal F))$.

Contracting a radical ideal of $R_\pi[1/\ell]$ to $R_\pi$ produces a radical $\ell$-saturated ideal of $R_\pi$. By Proposition~\ref{prop:closed-subspace-charts} and Lemma~\ref{lem:standard-reduced}, the standard closed rigid subspaces $Z_{J_m^i}$, $Z_{J_m^{i,c}}$, and $Z_{\mathfrak{n}_i}$ are reduced. By \eqref{eq:fm-specialization}, a rigid point $\chi\in X_\pi$ lies in $Z_{J_m^i}$ if and only if $\chi\in\Sigma_m^i(X,\mathcal F)$, and in $Z_{J_m^{i,c}}$ if and only if $\chi\in\Sigma_m^{i,c}(X,\mathcal F)$. By \eqref{eq:nf-locus-def}, $\chi$ lies in $Z_{\mathfrak n_i}$ if and only if $M^i$ fails to be flat at $\mathfrak p_\chi$.
\end{proof}

\begin{proposition}\label{prop:frob-equiv-loci}
Suppose that $u$ and $\mathcal F$ are arithmetic relative to a common descent over $F_0$. Let $\sigma\in\mathrm{Gal}(\overline F_0/F_0)$ preserve $\ker(u)$ and satisfy $\sigma^*\mathcal F_{\overline F_0}\simeq\mathcal F_{\overline F_0}$.
\begin{enumerate}
\item[(1)] The automorphism $\sigma^\#$ fixes each of the ideals $J_m^i$, $J_m^{i,c}$, and $\mathfrak n_i$ in $R_\pi$.
\item[(2)] The automorphism $\sigma$ preserves each of the reduced closed rigid subspaces $\Sigma_m^i(X,\mathcal F)_{\mathrm{red}}$, $\Sigma_m^{i,c}(X,\mathcal F)_{\mathrm{red}}$, and $\mathrm{NF}^i(X,\mathcal F)_{\mathrm{red}}$.
\end{enumerate}
\end{proposition}

\begin{proof}
Choose an isomorphism $\sigma^*\mathcal F_{\overline F_0}\xrightarrow{\sim}\mathcal F_{\overline F_0}$. The complex $\mathcal F_{\overline F_0}$, its identification with $\mathcal F$ after pullback to $X$, and the chosen isomorphism are defined over a finite extension $E'/E$ \cite[Section~1.1.3]{Del80}. Finite free extension of coefficients commutes with both Fourier--Mellin transforms. Since minors commute with base change and flatness descends along faithfully flat local maps, the three reduced ideals over $\mathcal O_{E'}[[\pi]]$ contract to $J_m^i$, $J_m^{i,c}$, and $\mathfrak n_i$, respectively. It therefore suffices to prove (1) after replacing $E$ by $E'$.

By Remark~\ref{rem:alg-closed-invariance}, we identify $\pi$ and the rational Fourier--Mellin transforms on $X$ with their counterparts on $X_{\overline F_0}$. In the following calculation $\mathcal L_\pi$ denotes the universal $R_\pi$-local system on $X_{\overline F_0}$. The actions on $V_\pi^\vee$ and $\pi$ satisfy $(\sigma\lambda)(\sigma\gamma)=\lambda(\gamma)$ for $\lambda\in V_\pi^\vee$ and $\gamma\in\pi$. On the coefficient ring, the induced automorphism is $\sigma^\#:R_\pi\to R_\pi$, $[\gamma]\mapsto[\sigma(\gamma)]$. Pulling back the universal character therefore gives $\sigma^*\mathcal L_\pi\simeq\mathcal L_\pi\otimes_{R_\pi,(\sigma^\#)^{-1}}R_\pi$. Hence $\sigma^*(\mathcal L_\pi\otimes_{R_\pi,\sigma^\#}R_\pi)\simeq\mathcal L_\pi$. Applying $\sigma^*$ to the twisted complex on $X_{\overline F_0}$ gives the first isomorphism in
\[
\mathrm{FM}_\pi^{(c)}(X,\mathcal F)\otimes_{R_\pi[1/\ell],\sigma^\#}^{\mathbf L} R_\pi[1/\ell]
\simeq \mathrm{FM}_\pi^{(c)}(X_{\overline F_0},\sigma^*\mathcal F_{\overline F_0})
\simeq \mathrm{FM}_\pi^{(c)}(X,\mathcal F).
\]
The second uses $\sigma^*\mathcal F_{\overline F_0}\simeq\mathcal F_{\overline F_0}$.
Choose a bounded finite free representative $P^\bullet$ of $\mathrm{FM}_\pi^{(c)}(X,\mathcal F_{\mathcal O_E})$. The scalar twist on the left has differential matrices obtained by applying $\sigma^\#$ to the matrices of $P^\bullet\otimes_{R_\pi} R_\pi[1/\ell]$. Their minors are the images under $\sigma^\#$ of the original minors. By Proposition~\ref{prop:standard-fm-loci}(1), the two isomorphic rational complexes yield the same radical jump ideal. Thus $\sigma^\#$ fixes the radical jump ideal $\sqrt{D_m^i(P^\bullet)R_\pi[1/\ell]}$. It also fixes the nonflat locus of $M^i$ in $\mathrm{Spec}(R_\pi[1/\ell])$, since the isomorphism of ordinary transforms induces a $\sigma^\#$-semilinear automorphism of $M^i$. Contracting these radical ideals to $R_\pi$ gives (1). By Lemma~\ref{lem:faithfulness}, invariance of the defining ideals is equivalent to invariance of the standard subspaces in (2).
\end{proof}

Let $\pi'$ be a lattice containing $\pi$ with finite index as in Proposition~\ref{prop:weight-adapted-torus}(1). Restriction of characters gives the finite étale isogeny $f:X_{\pi'}\to X_\pi$. Let $\mathcal L_{\pi'}$ be the rank one $R_{\pi'}$-local system on $X$ defined by the composite character
\[
\pi_1(X)\longrightarrow\Gamma_\ell(X)\xrightarrow{u}\pi\hookrightarrow\pi'
\xrightarrow{\gamma'\mapsto[\gamma']} R_{\pi'}^\times.
\]
Using $\mathcal L_{\pi'}$ in Definition~\ref{def:fm-transform} gives the integral and rational transforms $\mathrm{FM}_{\pi'}^{(c)}$. We denote the resulting jump and nonflat loci in $X_{\pi'}$ by $\Sigma_{m,\pi'}^{i,(c)}(X,\mathcal F)$ and $\mathrm{NF}_{\pi'}^i(X,\mathcal F)$.

\begin{lemma}\label{lem:fm-isogeny-base-change}
Base change along $R_\pi\to R_{\pi'}$ and $R_\pi[1/\ell]\to R_{\pi'}[1/\ell]$ induces the canonical isomorphisms
\[
\begin{aligned}
\mathrm{FM}_\pi^{(c)}(X,\mathcal F_{\mathcal O_E})\otimes_{R_\pi}^{\mathbf L}R_{\pi'} &\xrightarrow{\sim} \mathrm{FM}_{\pi'}^{(c)}(X,\mathcal F_{\mathcal O_E}),\\
\mathrm{FM}_\pi^{(c)}(X,\mathcal F)\otimes_{R_\pi[1/\ell]}^{\mathbf L}R_{\pi'}[1/\ell] &\xrightarrow{\sim} \mathrm{FM}_{\pi'}^{(c)}(X,\mathcal F).
\end{aligned}
\]
\end{lemma}

\begin{proof}
By Lemma~\ref{lem:isogeny-properties}(1), the local homomorphism $R_\pi\to R_{\pi'}$ is finite free. The quotient $R_{\pi'}/\mathfrak m_{R_\pi}R_{\pi'}$ is Artinian local, with nilpotent maximal ideal. In particular, there is an integer $c\ge1$ such that for every $n\ge1$,
\[
\mathfrak m_{R_{\pi'}}^{cn}\subset\mathfrak m_{R_\pi}^{n}R_{\pi'}\subset\mathfrak m_{R_{\pi'}}^{n}.
\]
The two systems of quotients are therefore cofinal.

The monodromy of $\mathcal L_{\pi'}$ is the composite $\pi_1(X)\xrightarrow{\varrho_{\mathrm{univ}}}R_\pi^\times\to R_{\pi'}^\times$. Hence $\mathcal L_{\pi'}\simeq\mathcal L_\pi\otimes_{R_\pi}R_{\pi'}$. For $n\ge1$, let $A_n=R_\pi/\mathfrak m_{R_\pi}^n$ and $B_n=R_{\pi'}/\mathfrak m_{R_\pi}^nR_{\pi'}$. Since $B_n$ is finite free over $A_n$, the projection formula gives
\[
R\Gamma_{(c)}(X,\mathcal F_{\mathcal O_E}\otimes_{\mathcal O_E}^{\mathbf L}\mathcal L_{A_n})
  \otimes_{A_n}B_n
\xrightarrow{\sim}
R\Gamma_{(c)}(X,\mathcal F_{\mathcal O_E}\otimes_{\mathcal O_E}^{\mathbf L}\mathcal L_{B_n}).
\]
Since $R_{\pi'}$ is finite free over $R_\pi$, tensoring with $R_{\pi'}$ commutes with the derived inverse limit over $n$. By cofinality, the derived inverse limit on the right is $\mathrm{FM}_{\pi'}^{(c)}(X,\mathcal F_{\mathcal O_E})$. Tensoring the integral isomorphism over $R_{\pi'}$ with $R_{\pi'}[1/\ell]$ gives the rational isomorphism.
\end{proof}

\begin{corollary}\label{lem:pullback-fm-isogeny}
For the isogeny $f:X_{\pi'}\to X_\pi$, the jump and nonflat loci satisfy the following equalities of reduced closed rigid subspaces of $X_{\pi'}$:
\[
\begin{aligned}
f^{-1}\bigl(\Sigma_m^{i,(c)}(X,\mathcal F)_{\mathrm{red}}\bigr)&=\Sigma_{m,\pi'}^{i,(c)}(X,\mathcal F)_{\mathrm{red}},\\
f^{-1}\bigl(\mathrm{NF}^i(X,\mathcal F)_{\mathrm{red}}\bigr)&=\mathrm{NF}_{\pi'}^i(X,\mathcal F)_{\mathrm{red}}.
\end{aligned}
\]
\end{corollary}

\begin{proof}
For a rigid point $\chi'\in X_{\pi'}$, let $\chi:=f(\chi')$. The characters $\chi$ and $\chi'$ fit into the commutative diagram
\[
\begin{tikzcd}[column sep=1.8em]
\pi \arrow[r,hook] \arrow[d,"\chi"'] & \pi' \arrow[d,"\chi'"] \\
\kappa(\chi)^\times \arrow[r,hook] & \kappa(\chi')^\times\mathrlap{,}
\end{tikzcd}
\]
where the lower arrow is induced by the residue field inclusion $\kappa(\chi)\hookrightarrow\kappa(\chi')$. Hence $\mathcal{L}_{\chi'}\simeq\mathcal{L}_\chi\otimes_{\kappa(\chi)}\kappa(\chi')$. Taking (compactly supported) cohomology gives
\[
\dim_{\kappa(\chi')}H^i_{(c)}(X,\mathcal F\otimes_E^{\mathbf L}\mathcal L_{\chi'})
=\dim_{\kappa(\chi)}H^i_{(c)}(X,\mathcal F\otimes_E^{\mathbf L}\mathcal L_\chi).
\]
Hence $f^{-1}(\Sigma_m^{i,(c)}(X,\mathcal F))=\Sigma_{m,\pi'}^{i,(c)}(X,\mathcal F)$ as point sets.

Since $R_{\pi'}[1/\ell]$ is flat over $R_\pi[1/\ell]$, Lemma~\ref{lem:fm-isogeny-base-change} gives $H^i(\mathrm{FM}_{\pi'}(X,\mathcal F))\simeq M^i\otimes_{R_\pi[1/\ell]} R_{\pi'}[1/\ell]$. For each prime $\mathfrak q\subset R_{\pi'}[1/\ell]$ above $\mathfrak p\subset R_\pi[1/\ell]$, the local homomorphism $(R_\pi[1/\ell])_{\mathfrak p}\to(R_{\pi'}[1/\ell])_{\mathfrak q}$ is faithfully flat. Thus the base changed module is flat at $\mathfrak q$ if and only if $M^i$ is flat at $\mathfrak p$. This gives $f^{-1}(\mathrm{NF}^i(X,\mathcal F))=\mathrm{NF}_{\pi'}^i(X,\mathcal F)$.

Since the inverse image under $f$ of a standard closed rigid subspace is already reduced and standard by Lemma~\ref{lem:isogeny-standard-pullback}, Corollary~\ref{cor:standard-point-ideals}(2) implies that the reduced closed rigid subspace structures coincide.
\end{proof}

\subsection{Proof of the main theorem}\label{subsec:structure-thm}

\begin{theorem}[Quasilinearity of mixed weight Fourier--Mellin loci]\label{thm:main-fm}
Suppose that $u:\Gamma_\ell(X)\twoheadrightarrow\pi$ and $\mathcal F$ are arithmetic relative to a common descent and that $V_\pi^\vee$ is mixed with all weights strictly positive. For $i\in\mathbf Z$ and $m\ge1$, let $Z$ be any of the three reduced closed rigid subspaces
\[
\Sigma_m^i(X,\mathcal F)_{\mathrm{red}},\qquad
\Sigma_m^{i,c}(X,\mathcal F)_{\mathrm{red}},\qquad
\mathrm{NF}^i(X,\mathcal F)_{\mathrm{red}}.
\]
\begin{enumerate}
\item[(1)] $Z$ is rigid quasilinear on $X_\pi$.
\item[(2)] Let $f:X_{\pi'}\to X_\pi$ be the finite étale isogeny of Proposition~\ref{prop:weight-adapted-torus}(2), and let $\Phi_w:X_{\pi'}\to X_{\pi'}$ be the rigid generic fibre of the mixed power endomorphism of Proposition~\ref{prop:weight-adapted-torus}(3). Then $f^{-1}(Z)$ is $\Phi_w$-stable rigid quasilinear on $X_{\pi'}$.
\end{enumerate}
\end{theorem}

\begin{proof}
Choose the common field of definition $F_0$ and the lift $\sigma\in\mathrm{Gal}(\overline F_0/F_0)$ of $\mathrm{Frob}_s^N$ supplied by Lemma~\ref{lem:frob-lift}. By Propositions~\ref{prop:standard-fm-loci} and~\ref{prop:frob-equiv-loci}(1), $Z = Z_J$ is a standard closed rigid subspace defined by a radical $\ell$-saturated ideal $J \subset R_\pi$ satisfying $\sigma^\#(J) = J$.

Let $\pi\subset\pi'$ be the $\sigma$-stable overlattice of Proposition~\ref{prop:weight-adapted-torus}(1) and let $Y=f^{-1}(Z)$. By Lemma~\ref{lem:isogeny-standard-pullback} and Corollary~\ref{lem:pullback-fm-isogeny}, $Y$ is the corresponding locus on $X_{\pi'}$. Its defining ideal $JR_{\pi'}$ is radical, $\ell$-saturated, and fixed by $\sigma^\#$.

The lattice $\pi'$ and the automorphism $\sigma$ satisfy Assumption~\ref{ass:frobenius-weights} with $c=q_s^{-N/2}\ne1$. For every $d$ with $Z_{(JR_{\pi'})^{(d)}}\ne\varnothing$, Proposition~\ref{prop:frob-gives-phi-criterion} gives integers $a_d>b_d\ge0$ such that
\[
\Phi_w^{a_d}\bigl(Z_{(JR_{\pi'})^{(d)}}\bigr)
\subset
\Phi_w^{b_d}\bigl(Z_{(JR_{\pi'})^{(d)}}\bigr).
\]
By Corollary~\ref{cor:isogenous-rigid-criterion}, $Y$ is $\Phi_w$-stable rigid quasilinear on $X_{\pi'}$ and $Z$ is rigid quasilinear on $X_\pi$.
\end{proof}

\begin{corollary}\label{cor:pure-weight}
Under the assumptions of Theorem~\ref{thm:main-fm}, suppose that $V_\pi^\vee$ is pure of weight $w>0$. Then the reduced rigid subspaces $\Sigma_m^i(X,\mathcal F)_{\mathrm{red}}$, $\Sigma_m^{i,c}(X,\mathcal F)_{\mathrm{red}}$, and $\mathrm{NF}^i(X,\mathcal F)_{\mathrm{red}}$ are $\Phi_w$-stable rigid quasilinear on $X_\pi$.
\end{corollary}

\begin{proof}
Since $V_\pi^\vee$ is pure of weight $w>0$, the weight decomposition of $V_\pi$ has only the summand $V_w$. Thus we take $\pi' = \pi$ and $f = \mathrm{id}_{X_\pi}$ in Theorem~\ref{thm:main-fm}(2).
\end{proof}

\begin{corollary}\label{cor:torsion-detection}
Under the assumptions of Theorem~\ref{thm:main-fm}, the following conditions are equivalent for each $i\in\mathbf Z$ and $m\ge1$:
\begin{enumerate}
\item[(1)] The jump locus $\Sigma_m^{i,(c)}(X,\mathcal F)$ is nonempty.
\item[(2)] There exist a finite extension $L/E$ and a torsion character $\zeta\in X_\pi(L)$ whose associated rank one $L$-local system $\mathcal L_\zeta$ satisfies $\dim_L H^i_{(c)}(X,\mathcal F\otimes_E^{\mathbf L}\mathcal L_\zeta)\ge m$.
\end{enumerate}
\end{corollary}

\begin{proof}
By Theorem~\ref{thm:main-fm}(1), after a finite extension of $E$ the jump locus is a finite union of torsion translates of subtori. Each translate contains its translating torsion point because the subtorus contains the identity. Thus a nonempty jump locus contains a torsion point. The converse follows from the definition of the jump locus.
\end{proof}

\subsection{Examples}\label{subsec:examples}

\begin{example}[Open curves]\label{ex:open-curves}
Let $\overline Y_0$ be a smooth projective geometrically connected curve of genus $g\ge1$ over a finite field $\mathbf F_q$ of characteristic different from $\ell$, and let $D_0$ be a nonempty reduced divisor of degree $r$ on $\overline Y_0$. Write $\overline Y$ and $D$ for their base changes to $\overline{\mathbf F}_q$ and let $Y=\overline Y\setminus D$. Restriction to $Y$ and the residues along $D$ give the Frobenius-equivariant localization sequence
\[
 0\longrightarrow
 \underbrace{H^1(\overline Y,\mathbf Q_\ell)}_{\operatorname{Gr}_1^W H^1(Y,\mathbf Q_\ell)}
 \longrightarrow H^1(Y,\mathbf Q_\ell)
 \longrightarrow
 \underbrace{\ker\!\bigl(H^0(D,\mathbf Q_\ell)\xrightarrow{\sum}\mathbf Q_\ell\bigr)(-1)}
 _{\operatorname{Gr}_2^W H^1(Y,\mathbf Q_\ell)}
 \longrightarrow0.
\]
The weight one piece has dimension $2g$. After a finite extension of $\mathbf F_q$ splitting $D_0$, the weight two piece is isomorphic to $\mathbf Q_\ell(-1)^{\oplus(r-1)}$. Thus $H^1(Y,\mathbf Q_\ell)$ is pure of weight one if $r=1$ and has both weights one and two if $r\ge2$.

Suppose $r\ge2$ and let $u:\Gamma_\ell(Y)\twoheadrightarrow\pi$ be an arithmetic quotient for which $\operatorname{Gr}_1^W V_\pi^\vee$ and $\operatorname{Gr}_2^W V_\pi^\vee$ are both nonzero, as for $\pi=\Gamma_\ell(Y)$. On the character torus $X_{\pi'}$ of Proposition~\ref{prop:weight-adapted-torus}, write $\boldsymbol q_1$ and $\boldsymbol q_2$ for the tuples of multiplicative coordinates corresponding to the two weights of $V_\pi^\vee$. In these coordinates the mixed power endomorphism is $\Phi_w:X_{\pi'}\to X_{\pi'}$, $(\boldsymbol q_1,\boldsymbol q_2)\mapsto(\boldsymbol q_1^\ell,\boldsymbol q_2^{\ell^2})$, with powers taken coordinatewise.
\end{example}

\begin{example}[Passing to the maximal quotient can destroy purity]\label{ex:zero-weight}
Let $C_0$ be the split nodal rational curve obtained by identifying $0$ and $\infty$ in $\mathbf P^1_{\mathbf F_q}$, where $\mathbf F_q$ has characteristic different from $\ell$. Let $F=\overline{\mathbf F}_q$, $C=C_0\times_{\mathbf F_q}F$, and $X=C\times_F\mathbf G_{m,F}$, with basepoint $(1,1)$. The graph loop $\gamma_0$ of $C$ and the tame pro-$\ell$ generator $\gamma_2$ of $\mathbf G_m$ give
\[
\Gamma_\ell(X)=\mathbf Z_\ell\gamma_0\oplus\mathbf Z_\ell\gamma_2,
\qquad \sigma(\gamma_0)=\gamma_0,\qquad \sigma(\gamma_2)=q^{-1}\gamma_2,
\]
where $\sigma$ is induced by geometric Frobenius. Let $e_0,e_2$ be the dual basis of $H^1(X,\mathbf Q_\ell)$.

Consider the quotient $\pi_2=\Gamma_\ell(X)/\mathbf Z_\ell\gamma_0$. Its kernel is fixed by $\sigma$, so the quotient is arithmetic, and $V_{\pi_2}^\vee=\mathbf Q_\ell e_2$ is pure of weight two. Passing to the maximal quotient replaces $V_{\pi_2}^\vee$ by $H^1(X,\mathbf Q_\ell)=\mathbf Q_\ell e_0\oplus\mathbf Q_\ell e_2$, which has weights zero and two.

The change is also visible in the jump loci. Choose a nontorsion unit $a\in1+\ell\mathbf Z_\ell$, with $a\in1+4\mathbf Z_2$ if $\ell=2$. Let $M_a$ be the rank one local system on $C$ defined by the character $\gamma_0\mapsto a$, and let $\mathcal F=i_*M_a$, where $i:C\hookrightarrow X$ is the inclusion of $C\times\{1\}$. Since $\sigma$ fixes $\gamma_0$ and $i$ is defined over $\mathbf F_q$, $\mathcal F$ is arithmetic. For a rigid character $\chi$ of $\Gamma_\ell(X)$, the projection formula gives
\[
H^0(X,\mathcal F\otimes_{\mathbf Q_\ell}^{\mathbf L}\mathcal L_\chi)
\simeq H^0(C,M_a\otimes_{\mathbf Q_\ell}^{\mathbf L}i^*\mathcal L_\chi).
\]
On the normalization $\mathbf P^1_F$, the twisted local system is constant. A section with value $v\in\kappa(\chi)$ descends across the node if and only if $a\chi(\gamma_0)v=v$. Hence the displayed cohomology group is nonzero if and only if $\chi(\gamma_0)=a^{-1}$.

Every character factoring through $\pi_2$ satisfies $\chi(\gamma_0)=1$, so $\Sigma_{1,\pi_2}^0(X,\mathcal F)=\varnothing$. On $X_{\Gamma_\ell(X)}$, let $H=\{\chi\mid\chi(\gamma_0)=1\}$ and let $\chi_a$ be the character with $\chi_a(\gamma_0)=a^{-1}$ and $\chi_a(\gamma_2)=1$. Then $H$ is a codimension one rigid subtorus and
\[
\Sigma_{1,\Gamma_\ell(X)}^0(X,\mathcal F)
=\{\chi\in X_{\Gamma_\ell(X)}\mid\chi(\gamma_0)=a^{-1}\}
=\chi_a H.
\]
Every character in $\chi_aH$ takes the nontorsion value $a^{-1}$ at $\gamma_0$. Thus the jump locus for the maximal quotient is nonempty and contains no torsion point, so it is not rigid quasilinear. This gives a counterexample to the reduction in the proof of \cite[Theorem~6.2, p.~20]{EK21} when only the quotient is changed.
\end{example}

\begin{remark}[Comparison with Esnault and Kerz]\label{rem:point-set-comparison}
To avoid the loss of purity in Example~\ref{ex:zero-weight}, we retain the given arithmetic quotient $u:\Gamma_\ell(X)\twoheadrightarrow\pi$. Lemma~\ref{lem:frob-lift} supplies a lift of a power of geometric Frobenius acting semisimply on $V_\pi$ itself. Esnault and Kerz formulate their conclusions for Zariski closed subsets of $\mathcal G(\overline{\mathbf Q}_\ell)=\mathrm{Hom}_{\mathrm{cont}}(\pi,\overline{\mathbf Q}_\ell^\times)$, identified with $\mathrm{Spm}(B)$ for $B=R_\pi[1/\ell]\otimes_E\overline{\mathbf Q}_\ell$ \cite[Sections~3.1 and~6.1]{EK21}. To compare with Theorem~\ref{thm:main-fm}, we first identify these characters with the points of $X_\pi$ over finite extensions of $E$.

Choose a finite extension $E'/E$ in $\overline{\mathbf Q}_\ell$ over which a decomposition on $X_\pi$ from Theorem~\ref{thm:main-fm}(1) is defined. For a rigid space $Z$ over $E'$, write $Z(\overline{\mathbf Q}_\ell)=\bigcup_L Z(L)$, where $L$ ranges over the finite extensions of $E'$ in $\overline{\mathbf Q}_\ell$. Then for every primitive submodule $M\subset\pi$ we have
\[
X_\pi(\overline{\mathbf Q}_\ell)=\mathcal G(\overline{\mathbf Q}_\ell),\qquad
H_M(\overline{\mathbf Q}_\ell)=\mathrm{Hom}_{\mathrm{cont}}(\pi/M,\overline{\mathbf Q}_\ell^\times).
\]
Indeed, Remark~\ref{rem:points-character-torus}(1) identifies $X_\pi(L)$ with $\mathrm{Hom}_{\mathrm{cont}}(\pi,L^\times)$. Conversely, given $\chi\in\mathcal G(\overline{\mathbf Q}_\ell)$, choose a $\mathbf Z_\ell$-basis $e_1,\ldots,e_b$ of $\pi$ and let $L=E'(\chi(e_1),\ldots,\chi(e_b))$. The $\mathbf Z$-span of the $e_i$ is dense in $\pi$, and $L^\times$ is closed in $\overline{\mathbf Q}_\ell^\times$ because $L$ is complete. Continuity gives $\chi(\pi)\subset L^\times$, proving the first equality. The second follows by applying the first to $\pi/M$, since $H_M=X_{\pi/M}$.

Under these identifications, a torsion translate $\zeta H_M$ gives the coset $\zeta\,\mathrm{Hom}_{\mathrm{cont}}(\pi/M,\overline{\mathbf Q}_\ell^\times)$ of \cite[Definition~3.3]{EK21}. Conversely, the characters of each of the three reduced loci in Proposition~\ref{prop:standard-fm-loci}(2) determine its defining ideal by Corollary~\ref{cor:standard-point-ideals}(2) and $\ell$-saturation.

For the cohomology jump loci, the condition in \cite[Section~6.1]{EK21} is $\dim H^i>j$. Thus for $j\ge0$ their locus $\Sigma_\pi^i(\mathcal F,j)$ consists of the $\overline{\mathbf Q}_\ell$-points of $\Sigma_{j+1}^i(X,\mathcal F)_{\mathrm{red}}$, and for $j<0$ it is all of $\mathcal G(\overline{\mathbf Q}_\ell)$.

For the nonflat locus, the Fourier--Mellin cohomology module over $B$ in \cite[Theorem~6.2(2)]{EK21} is $M^i\otimes_{R_\pi[1/\ell]}B$, where $M^i=H^i(\mathrm{FM}_\pi(X,\mathcal F))$. Let $\mathfrak m_\chi$ be the maximal ideal of $B$ corresponding to $\chi$, and let $\mathfrak p_\chi=\mathfrak m_\chi\cap R_\pi[1/\ell]$. Since $(R_\pi[1/\ell])_{\mathfrak p_\chi}\to B_{\mathfrak m_\chi}$ is faithfully flat, $M^i$ is flat at $\mathfrak p_\chi$ if and only if $M^i\otimes_{R_\pi[1/\ell]}B$ is flat at $\mathfrak m_\chi$. Hence the nonflat locus of \cite[Theorem~6.2(2)]{EK21} consists of the $\overline{\mathbf Q}_\ell$-points of $\mathrm{NF}^i(X,\mathcal F)_{\mathrm{red}}$.

For the given arithmetic quotient, Corollary~\ref{cor:pure-weight} therefore recovers both conclusions of \cite[Theorem~6.2]{EK21} in pure weight and gives the compactly supported analogue for cohomology jump loci. Under the common arithmetic descent hypotheses of Theorem~\ref{thm:main-fm}, the same two conclusions hold when $V_\pi^\vee$ is mixed with strictly positive weights, as conjectured in the first paragraph of \cite[Remark~6.3]{EK21}.
\end{remark}

\section*{Acknowledgements}
Many of the ideas in this paper grew out of my master's thesis, which I wrote at the University of Regensburg as part of the ALGANT program. I thank my thesis advisor, Prof.~Dr.~Moritz Kerz, for suggesting the problem and for contributing many ideas used here.

\end{document}